\documentclass{article}
\usepackage{arxiv}
\usepackage[utf8]{inputenc}
\usepackage[T1]{fontenc}    
\usepackage{csquotes}
\usepackage[backend=bibtex, sorting=none]{biblatex}
\usepackage{amsmath}
\usepackage{amsfonts}
\usepackage{graphicx}
\usepackage{caption} 
\usepackage[margin=10pt,font=small,labelfont=bf,labelsep=endash]{caption}
\usepackage{subcaption}
\usepackage{hyperref}

\title{Adaptive POD Galerkin technique for reservoir simulation and optimization}
\author{ \href{https://orcid.org/0000-0003-1674-4891}{\includegraphics[scale=0.06]{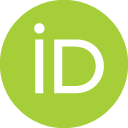}\hspace{1mm}Dmitry Voloskov}\\
        Skolkovo Institute of Science and Technology\\
        3 Nobelya st., Moscow, Russia 121205\\
        \texttt{dmitry.voloskov@skoltech.ru} \\
        \And
        Dimitri Pissarenko \\
        Skolkovo Institute of Science and Technology\\
        3 Nobelya st., Moscow, Russia 121205\\
        \texttt{d.pissarenko@skoltech.ru}\\
        TOTAL S.A. \\
        7 Lesnaya st., Moscow, Russia 124196\\
}
\hypersetup{
pdftitle={Adaptive POD Galerkin technique for reservoir simulation and optimization},
pdfauthor={Dmitry Voloskov, Dimitri Pissarenko},
pdfkeywords={Reservoir Simulation, Reduced Order Model, Proper Orthogonal Decomposition, Production Optimization},
}

\date{August 2020}

\begin{document}

\maketitle
\begin{abstract}
    In this work, a novel method with an adaptive functional basis for reduced order models (ROM) based on proper orthogonal decomposition (POD) is introduced.
    The method is intended to be applied in particular to hydrocarbon reservoir simulations, where a range of varying
    boundary conditions must be explored.
    The proposed method allows us to update the POD functional basis constructed for a specific problem setting in order to match varying boundary conditions,
    such as modified well locations and geometry, without the necessity to recalculate each time the whole set of basis functions.
    Such an adaptive technique allows us to
    significantly reduce the number of snapshots required to calculate the new basis, and hence reduce
    the computational cost of the simulations.
    The proposed method was applied to a two-dimensional immiscible displacement model, and simulations were performed using a high resolution model,
    a classical POD reduced model, and a reduced model whose POD basis was adapted to varying well location and geometry.
    Numerical simulations show that the proposed approach allows us to reduce
    the required number of model snapshots by a few orders of magnitude compared to a classical POD scheme,
    without noticeable loss of accuracy of calculated fluid production rates.
    Such an adaptive POD scheme can therefore provide a significant gain in computational efficiency for problems
    where multiple or iterative simulations with varying boundary conditions are required, such as optimization of well design or production optimization.
\end{abstract}
\section{Introduction}
Fluid flow simulation in subsurface porous media is an essential process
in virtually all reservoir engineering applications. Hydrodynamic reservoir models
are typically
based on a numerical solution of a system of nonlinear partial differential equations \cite{aziz_petroleum_1979} that describe the evolution in time of pressure field and fluid flow. These equations are discretized using such numerical methods as finite volumes \cite{monteagudo_control-volume_2004}
or finite elements \cite{young_finite-element_1981},
and are solved using IMPES (implicit pressure, explicit saturation), sequential, or
fully implicit schemes \cite{coats_compositional_1998,spillette_high-stability_1973,jiang_nonlinear_2019}.
In order to accurately reflect structural and physical complexities of a hydrocarbon reservoir, the numerical model should typically
contain millions of cells. Numerical simulations based on such models inevitably become a very challenging computational task.
Another factor that additionally increases the computational complexity of reservoir models
is the non-linearity of the governing equations. As a result, industry grade reservoir
simulations are a very resource consuming task that usually requires the use of high performance computing hardware. The same reason justifies a significant effort by the reservoir simulation community aimed at improving the computational efficiency of the simulations either through  parallelization or by algorithmic means. Regardless of such improvements, the
computational cost and consequently the duration of the simulation often become prohibitive for
such classes of problems as reservoir optimization and uncertainty quantification: they
typically require to perform thousands of runs of simulation scenarios.

In this paper, an applications of the Reduced Order Modeling (ROM) approach to reservoir simulation problems where a range of boundary conditions needs to be explored, is considered.
A new variant of proper orthogonal decomposition (POD) method, that allows us to considerably decrease the amount of calculations for each particular set of boundary conditions,
and hence to reduce the overall computational cost of the problem, is introduced.
ROM methods and their applications to resource demanding simulations in various areas of science and engineering are actively explored
in the recent years \cite{mehta_methodology_2017,shlizerman_proper_2012,xiao_non-intrusive_2015,yao_reduced-order_2019}.
A class of ROM methods frequently applied to large-scale simulation problems are
based on POD. In these methods, POD is used in order to obtain a
consistent representation of a model in a functional space whose dimension is lower than that of the original model.
To that end, the model's state variables are projected onto a
lower-dimensional POD domain, and Galerkin or least-squares Petrov-Galerkin projection is applied in order
to obtain a reduced system of equations \cite{carlberg_galerkin_2017}.

One of the main difficulties in the realization of efficient POD-based ROMs is related to handling the non-linearities of the model's equations.
In iterative algorithms frequently used for solving non-linear systems of equations such as Newton method,
calculations of the non-linear terms and the Jacobian estimation are required. To do so, one generally needs to project the approximated solution back to the original domain
and calculate the nonlinear functional in the full-scale domain on every iteration of the
algorithm. The cost of such full scale non-linear functional calculation and of the
forth-and-back projections may considerably offset the gain obtained
due to the model reduction.

Several approaches were suggested in order to treat non-linearities in an efficient way.
Chaturantabut and Sorensen introduced the Discrete Empirical Interpolation
Method (DEIM)\cite{chaturantabut_nonlinear_2010} in order to treat non-linearities in POD-based ROMs governed by systems of
time-dependent partial differential equations (PDE).
POD-DEIM has become one of the most widely used ROM methods applied to reservoir simulation problems
\cite{tan_trajectory-based_2019, efendiev_online_2016}.
Carlberg et al. developed the Gauss-Newton
with Approximated Tensors (GNAT) method \cite{carlberg_efficient_2011} which also uses
POD in order to reduce the vector of unknowns, but in contrast to DEIM it operates in a
fully discrete domain.
Jiang and Durlofsky successfully applied the GNAT approach to complex reservoir simulations \cite{jiang_implementation_2019}.
Rewienski and White developed a method called trajectory piecewise linearization (TPWL) \cite{rewienski_trajectory_2003}. In this method, a number of the system's states and Jacobians are first calculated and saved, then new
simulations are obtained as a result of linear expansions around previously saved states. This approach
can also be applied in a reduced subspace such as the one obtained through POD.
A combination of POD and TPWL (POD-TPWL) is now widely used in order to model subsurface flows
\cite{cardoso_linearized_2010,  he_enhanced_2011}. Trehan and Durlofsky
\cite{trehan_trajectory_2016} developed an extension of TPWL called trajectory piece-wise
quadratic extension procedure (TPWQ) and combined it with POD (POD-TPWQ).

In recent years, ROM methods using machine learning (ML) have as well been actively explored.
Kani and Elsheikh developed a deep residual recurrent neural network (DR-RNN)
approach \cite{kani_dr-rnn:_2017} and applied it to modeling of two-phase subsurface flows
\cite{kani_reduced-order_2018}. They used POD to project the original problem onto a reduced
subspace and applied a recurrent neural network (RNN) in order to model the dynamics in
the reduced space. Another group of methods use Variational Autoencoders (VAE) in order
to obtain a reduced representation of model's states. Lee and Carlberg \cite{lee_model_2019}
modified the GNAT approach and used
VAE instead of POD. Temirchev et al. \cite{temirchev_deep_2020} used VAE
combined with RNN in order to mimic the dynamics of subsurface flows. In \cite{temirchev_reduced_2019}, an approach called
Neural Differential Equations based ROM (NDE-b-ROM) was suggested: the authors applied
the Neural Ordinary Differential Equations method \cite{chen_neural_2019} to model
the dynamics in the reduced latent space, while the reduced representation was obtained with the help
of VAE. Fraces et al. applied feedforward neural networks to approximate derivatives in
Buckley-Leverett problem\cite{fraces_physics_2020}.
They used transfer learning approach and Generative Adversarial Networks to obtain
continuos in both time and space approximation of PDE solution.

A number of important reservoir engineering problems involve iterative simulations whose total
computational cost may be particularly high. A typical example of such problems are various
optimization tasks, such as finding the optimal well locations, well geometries, well
completion schemes, and well control schedules.
Numerical solution of such optimization problems requires multiple simulations of essentially the same reservoir unit with varying well parameters and schedules.
Reduced reservoir models including POD-based ROMs, were successfully applied to well control
optimization problems
\cite{jansen_use_2017,cardoso_use_2010,trehan_trajectory_2016,insuasty_tensor-based_2015,cardoso_linearized_2010,he_enhanced_2011}.
However, some reservoir optimization problems require simulations with varying well locations
and well geometries, in addition to different well control sequences. For such problems,
standard POD-based ROMs fail to correctly reproduce the flow dynamics even after a slight change
in well location or well geometry with respect to the model used to construct the basis.
In such cases, classical POD schemes require the construction of a new POD basis after every
change in well locations or geometry. That implies re-calculating from scratch of a new training
data set based on high-resolution model simulations with a significant associated computational overhead.

In this paper, a POD-based adaptive scheme that allows us to account for changes in
well location and well geometry at the expense of a relatively small computational overhead is introduced.
The proposed scheme requires significantly less training data compared to what is necessary for constructing a new POD-basis from scratch.

This paper proceeds as follows. In Sect. \ref{sec: problem_formulation}, the governing
equations for a two-phase immiscible displacement problem are presented together with an overview
of POD-Galerkin ROM method and its application to this problem. In Sect. \ref{sec: horizontal_well_pod},  a test problem related
to the optimization of the geometry of a horizontal production
well is set up, and two different approaches for construction of the POD-basis for such
kind of problems are described. The results of simulations using full resolution and adaptive POD models are compared
and discussed. In Sect. \ref{sec: summary}, advantages and shortcomings of the proposed adaptive POD scheme are discussed,
and directions for further work are suggested.

\section{POD-Galerkin ROM of Two-phase Immiscible Flow}\label{sec: problem_formulation}
\subsection{Two-phase immiscible flow}
The mathematical model of a two-phase immiscible flow is obtained by combining the system of mass conservation (continuity)
equations for fluid phases (oil and water), and Darcy's law for each phase. The continuity equation takes the form
\begin{equation}
    \frac{\partial\phi\rho_{o, w} s_{o, w}}{\partial t} -
    \nabla \cdot (\rho_{o, w}\mathbf{v}_{o, w}) + q_{o,w} = 0 ,
\end{equation}
where subscript $o, w$ denote oil and water phase respectively, $\phi$ is the porosity, $\rho$ - the fluid density,
$s$ - the fluid saturation, and  $\mathbf{v}$ is the Darcy velocity that can be expressed as follows
\begin{equation}\label{eq: darcy}
\mathbf{v}_{o, w} = -\lambda_{o, w}\mathbf{K}\nabla (p_{o, w} - \rho_{o, w}g h) .
\end{equation}
Here $\mathbf{K}$ is the absolute permeability tensor, $\lambda = \frac{k_r}{\mu}$ - the phase mobility,
$k_r$ is the relative permeability of the corresponding phase, and $\mu$ is the viscosity of this phase, $p$ - the fluid pressure,
$g$ - the gravitational
acceleration, $h$ - the depth,  and $q$ is the source or sink term \cite{kani_reduced-order_2018}.
After neglecting the capillary pressure, compressibility and gravitational effects,
mass conservation equation and Darcy's law are combined in order to obtain a system of equations for the global pressure

\begin{equation}\label{eq: pressure_equation}
    \nabla \cdot \mathbf{K}\lambda\nabla p = q ,
\end{equation}
and for the saturation of the water phase
\begin{equation}\label{eq: saturation_equation}
    \phi \frac{\partial s_w}{\partial t} + \nabla \cdot \mathbf{v}_w = \frac{q_w}{\rho_w} ,
\end{equation}
where $p = p_o = p_w$ is the global pressure, $\lambda  = \lambda_w + \lambda_o$ - the total mobility,
$q = q_w + q_o$ - is the source or sink term.
The discrete form of the model can be obtained by dividing the domain into blocks and by applying the
finite volume method to Eqs. \eqref{eq: pressure_equation} and \eqref{eq: saturation_equation}.
Discretized pressure equation takes the form
\begin{equation}\label{eq: discretized_p_equation}
    \mathbf{A}\mathbf{s}_p = \mathbf{b}\;,
\end{equation}
where $\mathbf{A} \in \mathbb{R}^{n\times n}$ is a coefficient matrix,
$\mathbf{s_p} \in \mathbb{R}^n$ is the pressure state vector, and $b \in  \mathbb{R}^n$ is the right hand side of the equation.
Each element of the unknown vector ${s_p}_i$ represents the mean pressure value in the $i$-th grid block.
Saturation equation takes the form
\begin{equation}\label{eq: discretized_s_equation}
    \frac{\mathrm{d} \mathbf{s}_s}{\mathrm{d}t} + \mathbf{B}(\mathbf{v})f_w(\mathbf{s}_s) = \mathbf{d},
\end{equation}
where $ \mathbf{s}_s \in \mathbb{R}^n$ is the saturation state vector, $\mathbf{v}\in \mathbb{R}^n$ is the velocity vector, obtained
from the pressure field, $\mathbf{B} \in \mathbb{R}^{n\times n}$ is a coefficient matrix
depending on the velocity vector, $f_w(\mathbf{s}_s)$  is the nonlinear term
depending on the saturation field, and $\mathbf{d} \in \mathbb{R}^n$ is right hand side of the equation.
Equations \eqref{eq: discretized_p_equation} and \eqref{eq: discretized_s_equation} are coupled
through a dependence of matrix $\mathbf{A}$ on the saturation field and through a dependence of the velocity
$\mathbf{v}$ on the pressure field. There are several ways of constructing the numerical solution of such coupled problems. In the present work,
the IMPES method \cite{fanchi_principles_2018} is used: at each
time step,
the saturation field from the previous time step is used to construct the matrix $\mathbf{A}$. The pressure
equation is solved using an implicit scheme in order to obtain the pressure field. The obtained pressure field is then used to calculate the velocity field
$\mathbf{v}$ and construct the matrix $\mathbf{B}$. After that Eq. \eqref{eq: discretized_s_equation}
is solved explicitly and the saturation field is obtained.

\subsection{POD-Galerkin model}
\subsubsection{POD basis}
The main objective of proper orthogonal decomposition is to obtain an optimal
low dimensional functional basis that is capable to adequately represent high dimensional data. Once constructed,
the POD basis can be used in order
to formulate a reduced order model that corresponds to the original high resolution model.
POD decomposes a given fluctuating field into an orthonormal
system of spatial modes $\mathbf{u_i}(x)$ and the corresponding temporal coefficients $a_i(t)$ \cite{kunisch_galerkin_2003}
\begin{equation}
    \mathbf{u'}(x) = \sum\limits_{i=1}^Na_i \mathbf{u_i}(x) .
\end{equation}
The discrete variant of POD is also known as principal component analysis (PCA), and both
methods are
closely related to singular value decomposition method (SVD) \cite{abdi_principal_2010}.
In order to generate a set
of POD modes, the data set should be presented as a matrix $\mathbf{X}$, where each row represents the variable field at a given instant, such that if the field state is represented by $n$ values, and the
data set consists of $m$ field states, then $ \mathbf{X} \in \mathbb{C}^{m\times n}$.

The optimality of the POD basis means that for any given dimensionality of the basis $r$,
the truncation error is minimal
\begin{equation}
    \int\limits_t\int\limits_x\left(
        u(x, t) - \sum\limits_1^ra_i(t)u_i(x)
    \right)^2 = \min_{\phi, b}\int\limits_t\int\limits_x\left(
        u(x, t) - \sum\limits_1^rb_i(t)\phi_i(x)
    \right)^2 .
\end{equation}
In the discrete case, this equation can be written as
\begin{equation} \label{eq: optimality_}
    \sum\limits_{i=1}^n\sum\limits_{j=1}^m
        \left(
            x_{i, j} - \sum\limits_{k=1}^r(u^k_ia^k_j)
        \right)^2 = \min\limits_{\phi, b}
    \sum\limits_{i=1}^n\sum\limits_{j=1}^m
    \left(
        x_{i, j} - \sum\limits_{k=1}^r(\phi^k_ib^k_j)
    \right)^2\;,
\end{equation}
where $x_{i, j}$ corresponds to the $j$-th value of the $i$-th state representation,
$u_i^k$ is the $i$-th value of the $k$-th basis vector $\mathbf{u}^k$, and
$a^k_j$  is a projection of the $j$-th field state onto the $ k$-th basis vector.
This result can be obtained by the factorization of matrix $X$ using SVD

\begin{equation}\label{eq: svd}
    \mathbf{X=U\Sigma V^*}\;,
\end{equation}
where $\mathbf{U} \in \mathbb{C}^{n\times n}$ is the left singular matrix,
$\mathbf{V} \in \mathbb{C}^{n\times n}$ is the right singular matrix,
$\mathbf{\Sigma} \in \mathbb{R}^{m\times n}$ is the singular matrix (a diagonal matrix with non
negative values $\sigma_i$ on the main diagonal).
In matrix form the factorization \eqref{eq: svd} can be written as
\begin{equation}
    \mathbf{X} =
        \left[
            \begin{array}{ccccc}
                \mathbf{u}_1 &
                \dots
                \mathbf{u}_i &
                \dots
                \mathbf{u}_m &
            \end{array}
        \right]
        \left[
            \begin{array}{ccccc}
                \sigma_1 & & & &\\
                & \ddots & & &\\
                & & \sigma_i & &\\
                & & & \ddots & \\
                & & & & \sigma_n \\
                0 & \dots & \dots & \dots & 0

            \end{array}
        \right]
        \left[
            \begin{array}{c}
                \mathbf{v}^*_1 \\
                \vdots\\
                \mathbf{v}^*_i \\
                \vdots\\
                \mathbf{v}^*_n
            \end{array}
        \right]\;,
\end{equation}
where $\mathbf{u}_i$ is an $m$-dimensional column vector, and
$\mathbf{v}^*_i$ is a $n$-dimensional row vector.
The singular matrix $\mathbf{\sigma}$ is constructed such that
$\sigma_1 \geq \sigma_2 \geq \dots 0 $ \cite{shlizerman_proper_2012}.
The optimal reduced basis is obtained by taking the first $r$ left singular vectors
$\mathbf{u}_i$.

In order to obtain a reduced POD basis for a specific problem,
one needs to construct a "snapshot" matrix $\mathbf{X}$ which is composed of state vectors $s$ obtained from the solution of the full system.
For two- or three-dimensional problems, the model's states are first flattened to vectors.
These vectors, called snapshots, are then stacked to compose the snapshot matrix. The reduced POD basis is obtained
by applying SVD to the snapshot matrix and keeping the first
$r$ columns of the calculated singular matrix.

\subsubsection{POD-Galerkin ROM}
POD-Galerkin ROM of a two-phase immiscible displacement in IMPES formulation can be stated as follows.
During the offline stage (also called the training stage), pressure snapshots are recorded, and snapshot matrix $\mathbf{X}_p$ is constructed.
Then SVD \eqref{eq: svd} is applied to the snapshot matrix, and the reduced basis
$\mathbf{U}^r_p = [\mathbf{u}_1 \dots \mathbf{u}_i \dots \mathbf{u}_r]$ is obtained.
During the online stage, a reduced representation of the pressure equation \eqref{eq: discretized_p_equation} is constructed.
It can be written as
\begin{equation}\label{eq: pressure_equation_reduced}
    \mathbf{A}^r\mathbf{s}_p^r=\mathbf{b}^r\,,
\end{equation}
where $\mathbf{A^r} = {\mathbf{U}^r_p}^\top\mathbf{A}\mathbf{U}^r_p$ is a projection
of the matrix equation onto the reduced subspace,
$\mathbf{s}_p^r = {\mathbf{U}^r_p}^\top\mathbf{s}_p$, and
${\mathbf{b}}^r = {\mathbf{U}^r_p}^\top\mathbf{b}$  are the projections of the state
vector and of the right hand side of the equation onto the reduced subspace.
This equation is solved in the reduced subspace in order to obtain a
new reduced pressure state $\mathbf{s}^r_p$. It is used to form a full representation of the
pressure state $\widetilde{\mathbf{s}}_p$, the velocity field $\widetilde{\mathbf{v}}$, and the coefficient matrix $\mathbf{B}(\mathbf{v})$.
The saturation field $\widetilde{\mathbf{s}}_s$ is then calculated explicitly, and a new
coefficient matrix $\widetilde{\mathbf{A}}$ is formed.
This procedure is repeated for the subsequent time steps.

\section{Well orientation optimization using POD-Galerkin ROM} \label{sec: horizontal_well_pod}
\subsection{Test problem setup}
\begin{figure}[htbp]
    \begin{subfigure}{0.5\textwidth}
        \begin{center}
            \includegraphics[width=0.9\textwidth]{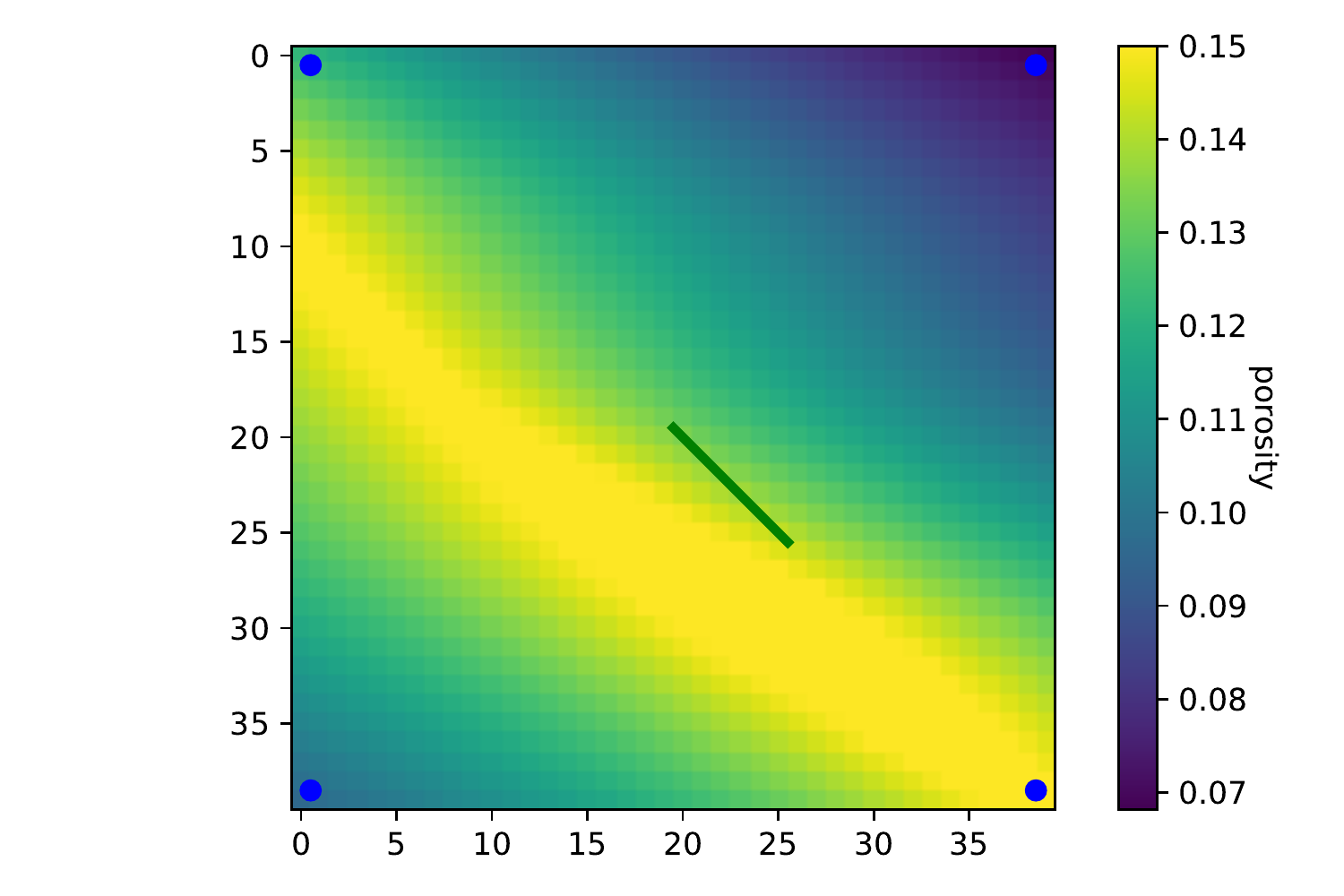}
            \caption{}
            \label{fig: porosity}
        \end{center}
    \end{subfigure}
    \begin{subfigure}{0.5\textwidth}
        \begin{center}
            \includegraphics[width=0.9\textwidth]{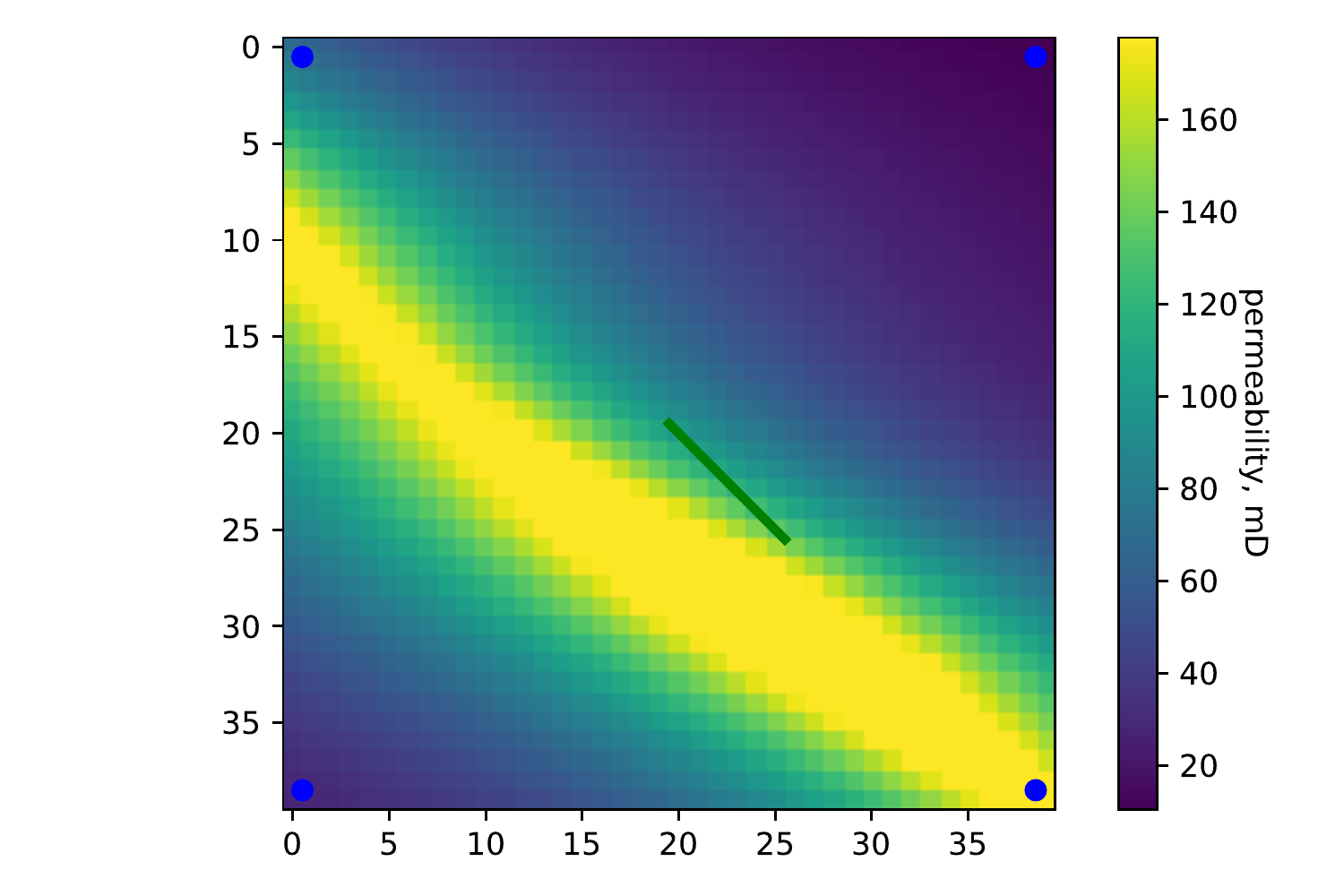}
            \caption{}
            \label{fig: permeability}
        \end{center}
    \end{subfigure}
    \caption{\ref{sub@fig: porosity} - porosity field, \ref{sub@fig: permeability} - permeability field of the model; blue dots in the corners indicate the positions of injector wells, green line shows the location of the producer.
             }
    \label{fig: porosity_permeability}
\end{figure}

Let us consider a simplified production optimization problem in which one needs to
optimize the orientation (azimuth) of a horizontal producing section of a well.
A two-dimensional immiscible displacement problem in a square
domain of the size $1000\times1000\; \text{m}$ is considered, and the domain is divided
into $40\times40$ square cells.
Heterogeneous porosity and permeability fields are generated numerically in order to mimic a
high permeability fluvial channel crossing a less permeable formation (fig. \ref{fig: porosity_permeability}).
In the corners of the square area are placed four injector wells that inject water at a controlled injection pressure. The oil that initially saturates the model is displaced toward the center of the simulation
area where fluids are recovered by a 150 meter long horizontal producer (fig. \ref{fig: porosity_permeability}).

\subsection{Universal POD-basis}
\begin{figure}[htbp]
    \begin{center}
        \includegraphics[width=100mm]{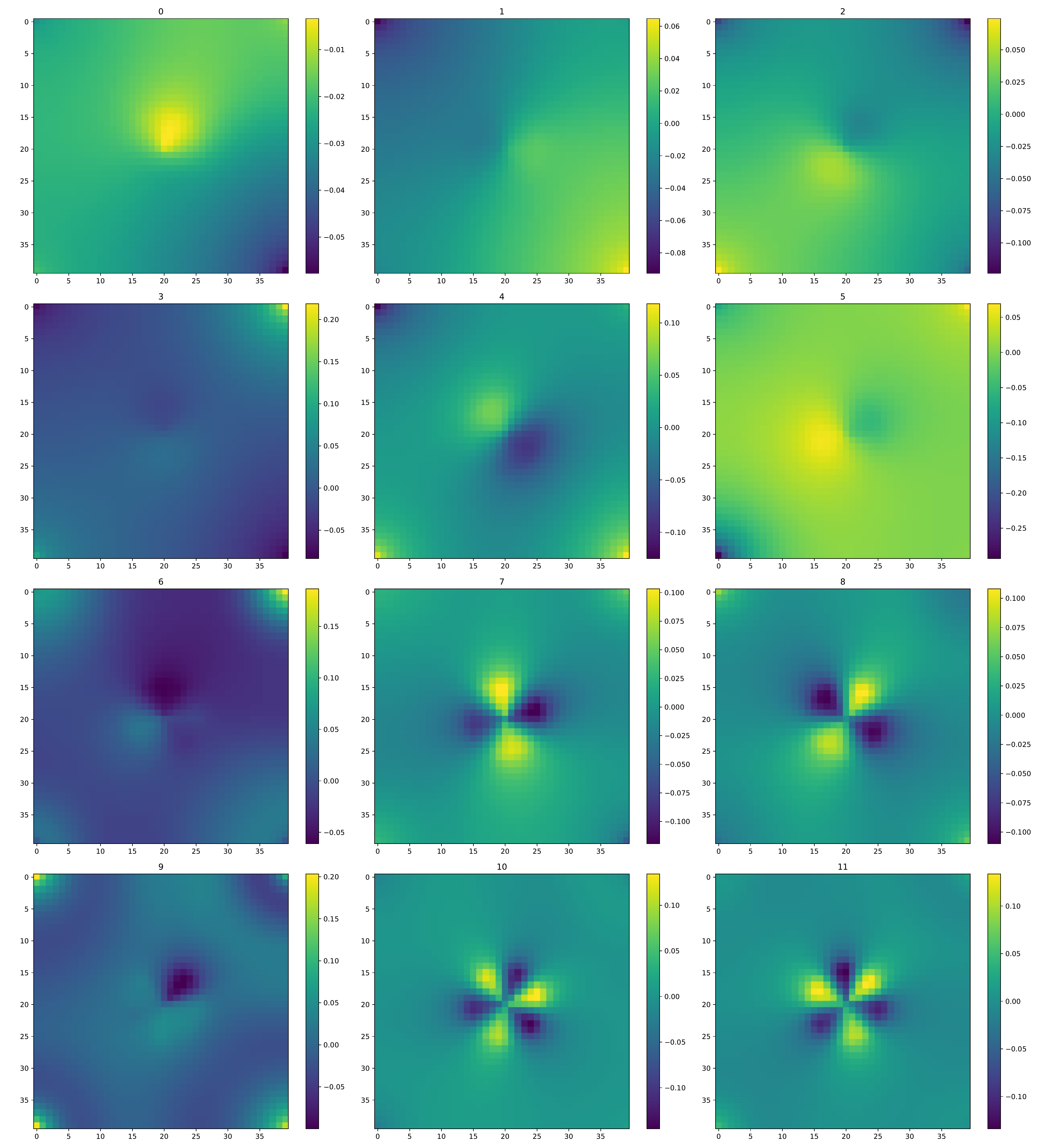}
    \end{center}
    \caption{The first 12 principal components of the universal POD basis.}
    \label{fig: universal_basis_components}
\end{figure}

\begin{figure}[htbp]
    \begin{center}
            \begin{subfigure}{0.45\textwidth}
            \begin{center}
                \includegraphics[width=\textwidth]{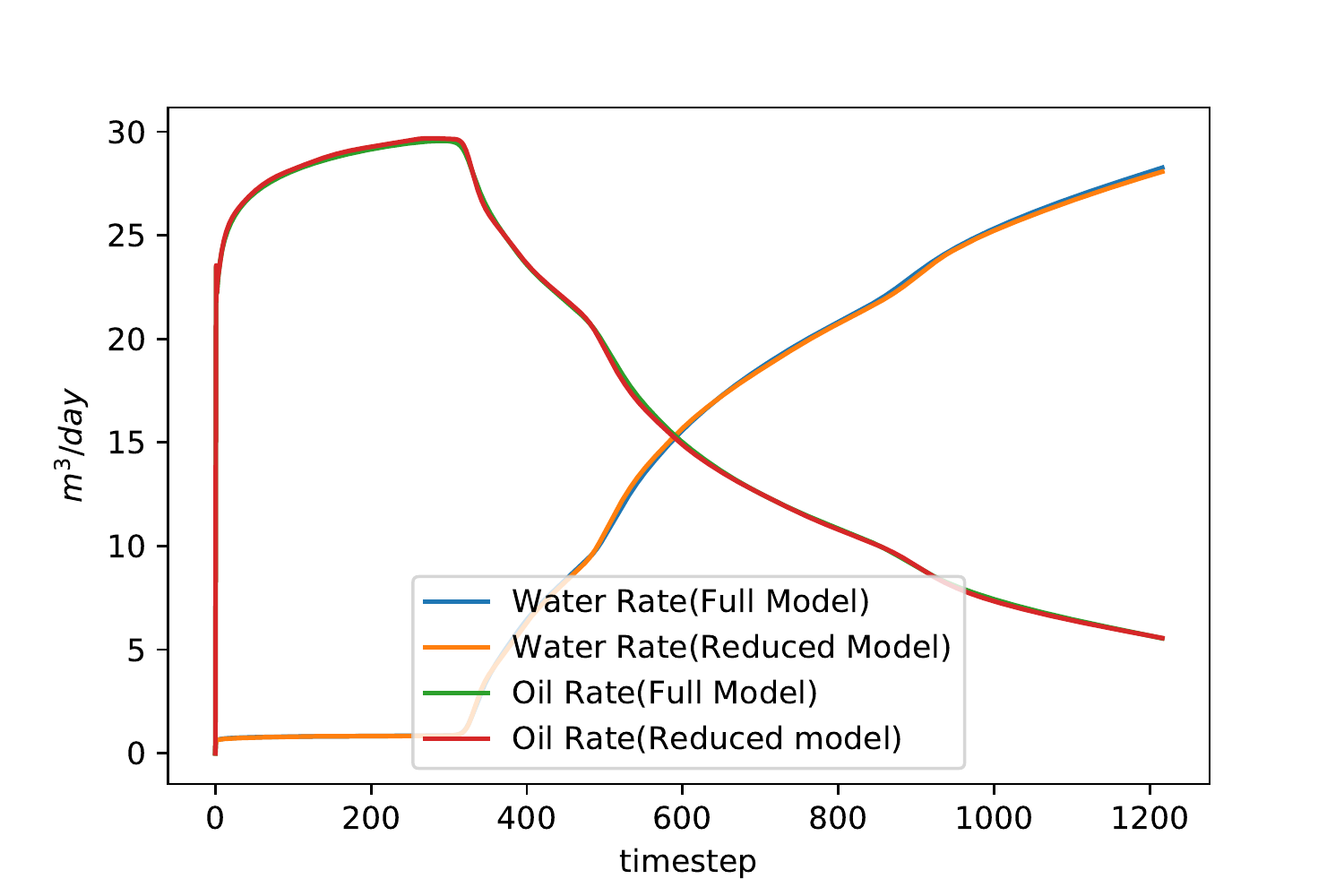}
            \end{center}
            \caption{}
            \label{fig: production_rates_63}
        \end{subfigure}
        \begin{subfigure}{0.45\textwidth}
            \begin{center}
                \includegraphics[width=\textwidth]{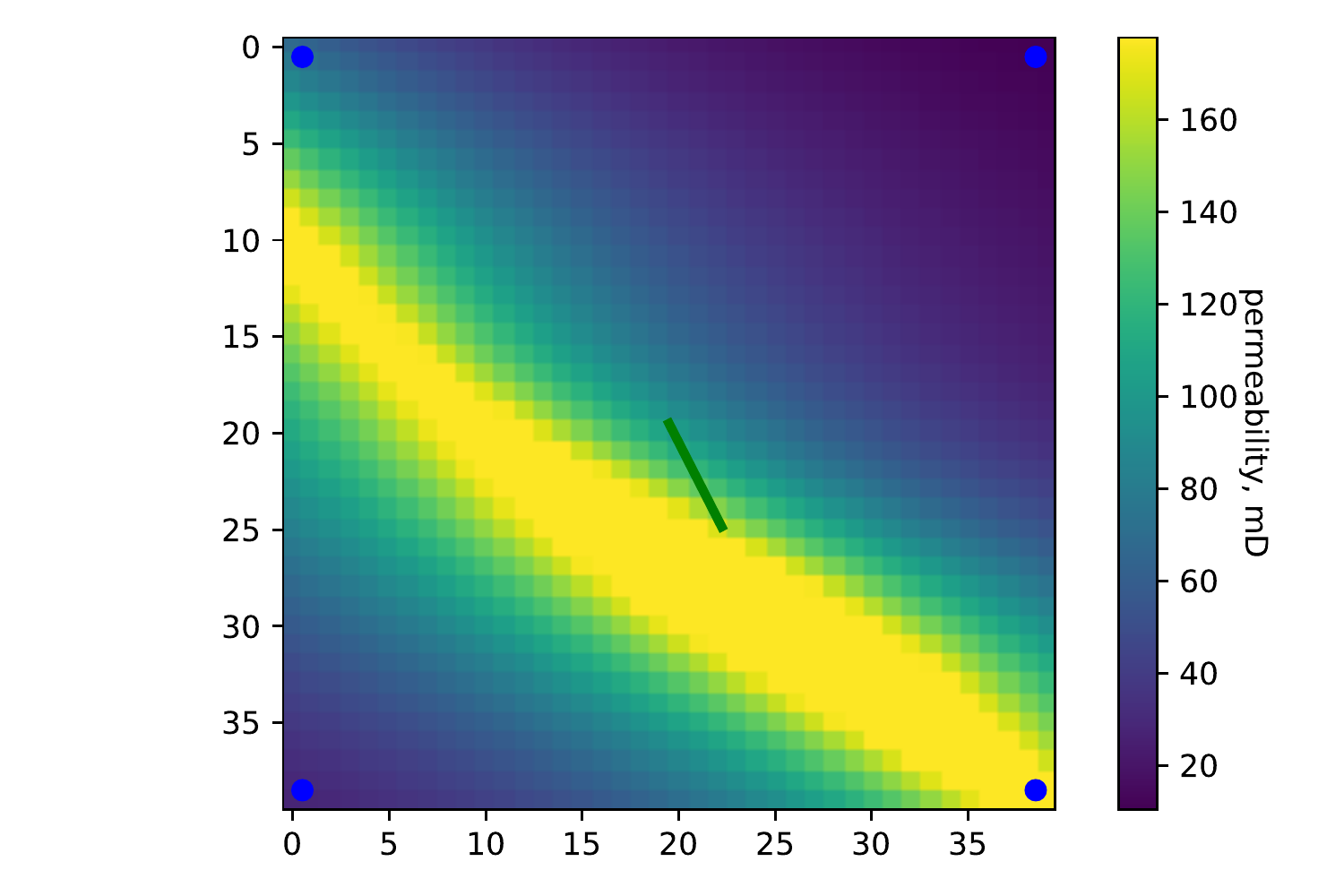}
            \end{center}
            \caption{}
            \label{fig: well_placement_63}
        \end{subfigure}
        \caption{Simulation of the fluid production of the well oriented 63 degrees clockwise
        from the horizontal axis.  \ref{sub@fig: production_rates_63} - production rates,
        \ref{sub@fig: well_placement_63} - well placement scheme}
        \label{fig: production_rates_63_full}
    \end{center}
\end{figure}

\begin{figure}[htbp]
    \begin{subfigure}{0.45\textwidth}
        \begin{center}
            \includegraphics[width=\textwidth]{figures/rates_100_components_175_degree.pdf}
        \end{center}
        \caption{}
        \label{fig: production_rates_175}
    \end{subfigure}
    \begin{subfigure}{0.45\textwidth}
        \begin{center}
            \includegraphics[width=\textwidth]{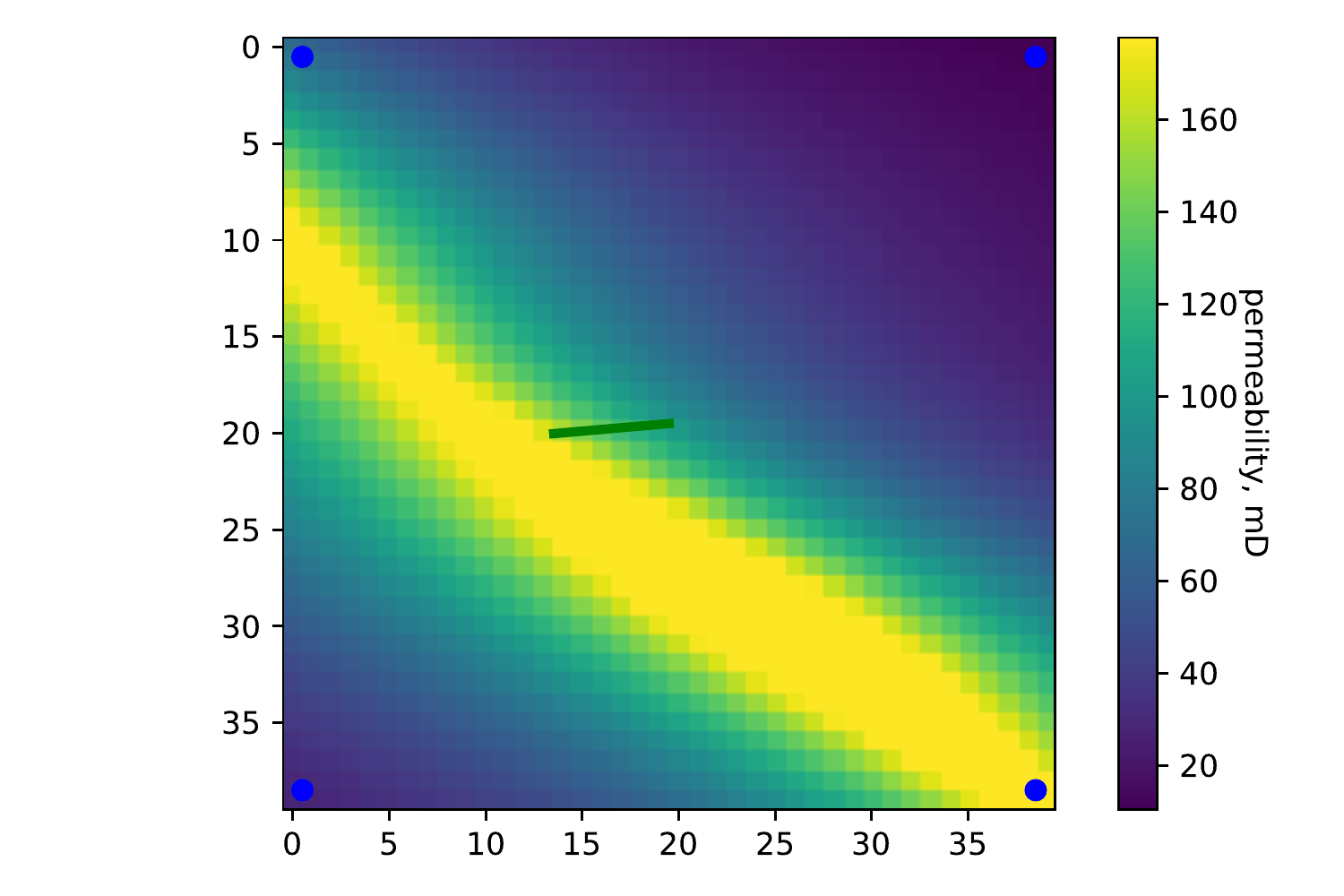}
        \end{center}
        \caption{}
        \label{fig: well_placement_175}
    \end{subfigure}
    \caption{Simulation with the producer oriented 175 degrees clockwise
    from the horizontal axis \ref{sub@fig: production_rates_175} - production rates,
             \ref{sub@fig: well_placement_175} - well placement scheme}
    \label{fig: production_rates_175_full}
\end{figure}

\begin{figure}[htbp]
    \begin{subfigure}{0.5\textwidth}
        \begin{center}
            \includegraphics[width=0.9\textwidth]{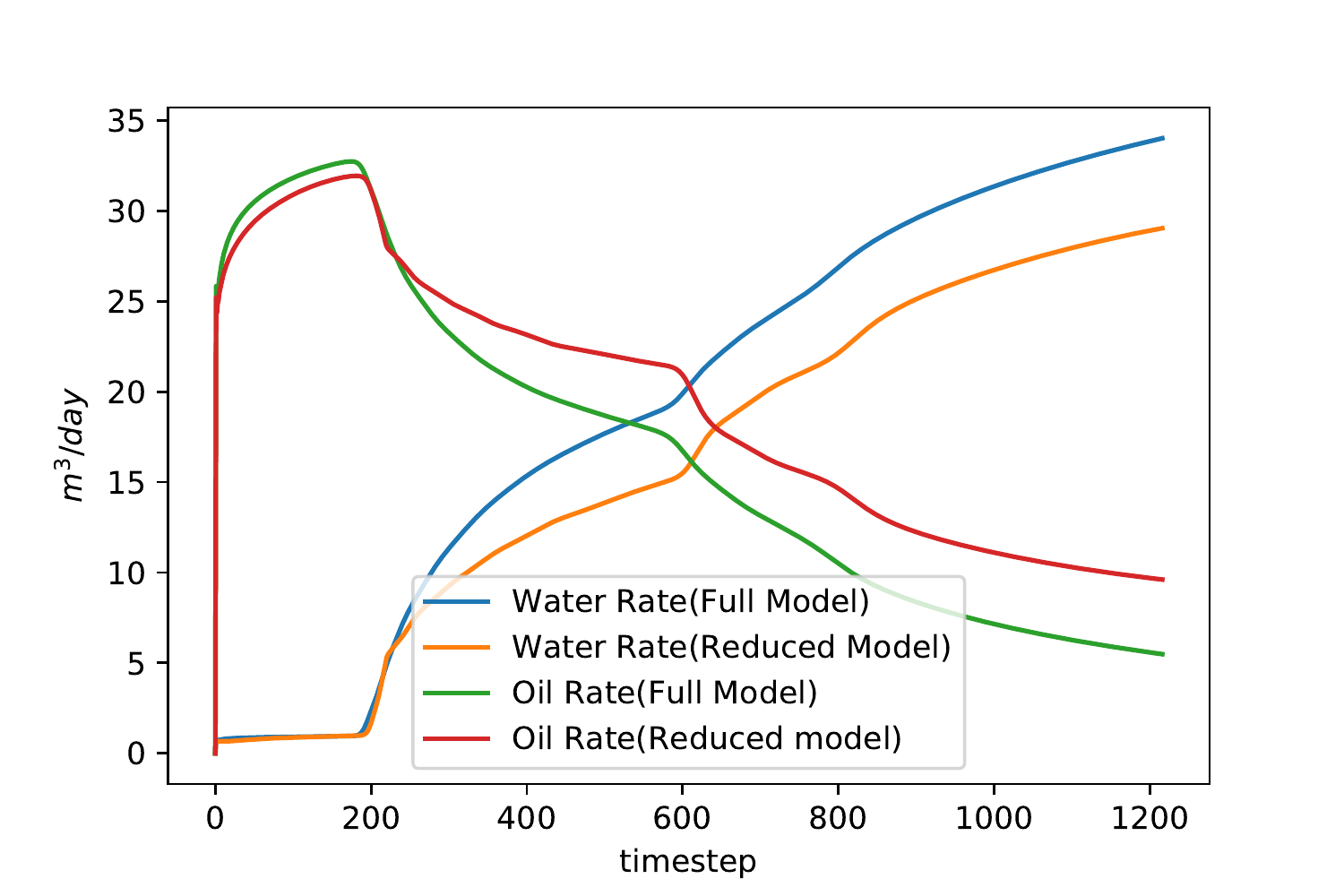}
        \end{center}
        \caption{}
        \label{fig: rates_63_20_components}
    \end{subfigure}
    \begin{subfigure}{0.5\textwidth}
        \begin{center}
            \includegraphics[width=0.9\textwidth]{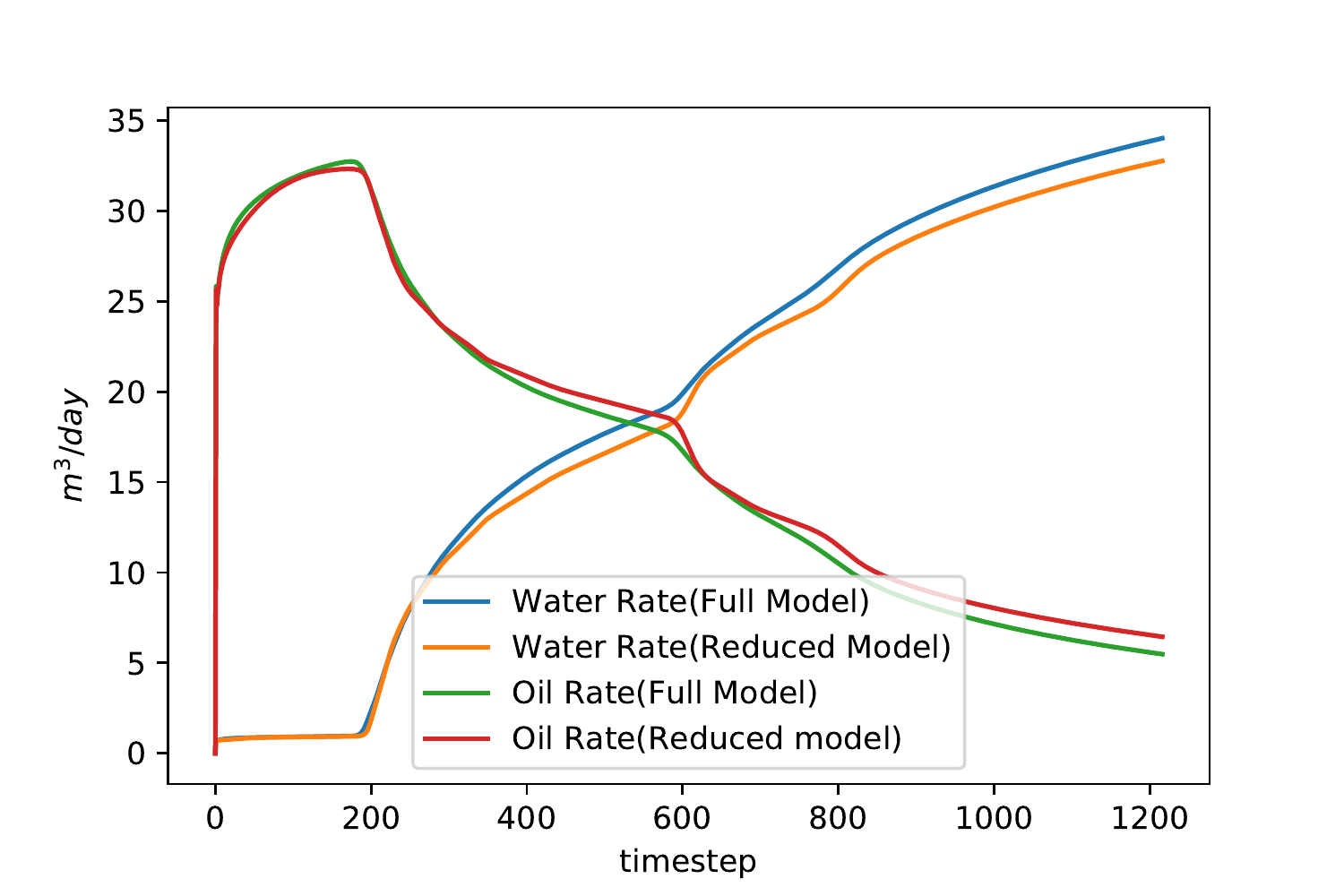}
        \end{center}
        \caption{}
        \label{fig: rates_63_40_components}
    \end{subfigure}
    \begin{subfigure}{0.5\textwidth}
        \begin{center}
            \includegraphics[width=0.9\textwidth]{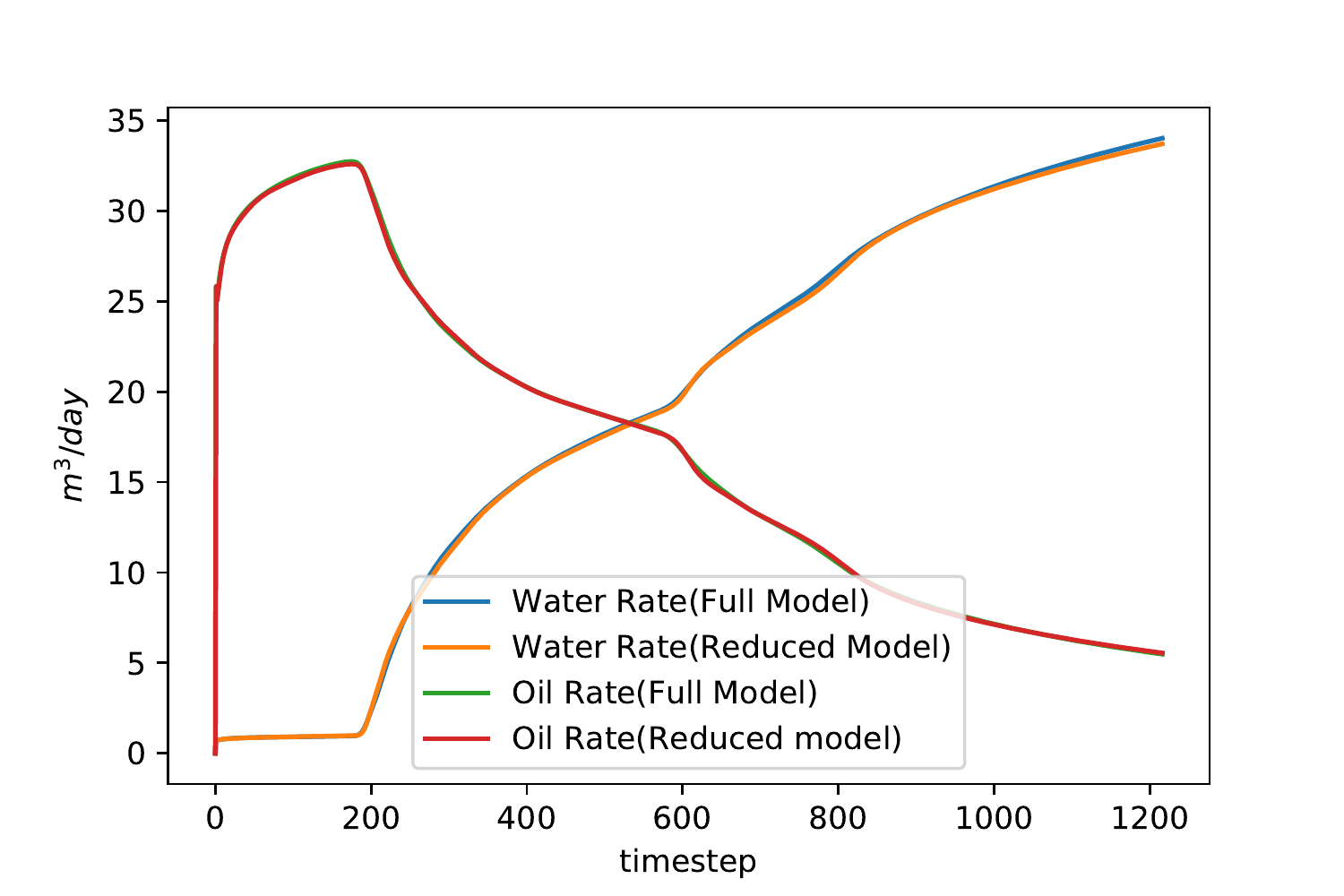}
        \end{center}
        \caption{}
        \label{fig: rates_63_100_components}
    \end{subfigure}
    \begin{subfigure}{0.5\textwidth}
        \begin{center}
            \includegraphics[width=0.9\textwidth]{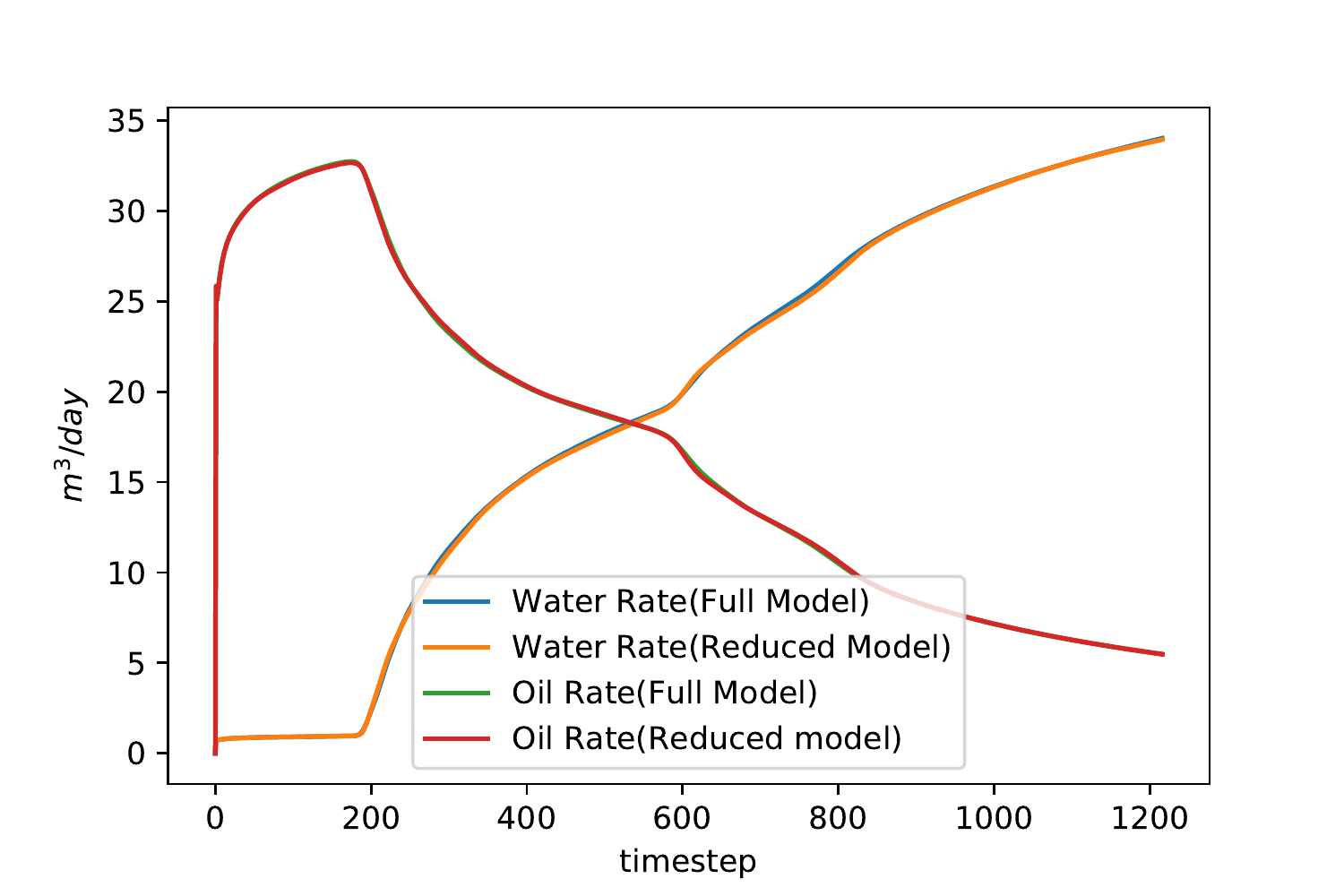}
        \end{center}
        \caption{}
        \label{fig: rates_63_120_components}
    \end{subfigure}
    \caption{Production rates simulated using different numbers of components of the reduced basis. \ref{sub@fig: rates_63_20_components} - 20 components, \ref{sub@fig: rates_63_40_components} - 40 components,
    \ref{sub@fig: rates_63_100_components} - 100 components,
    \ref{sub@fig: rates_63_120_components} - 120 components.}
    \label{fig: rates_various_components}
\end{figure}

In order to address optimal well placement problems, the POD-Galerkin model must be capable of simulating production scenarios with various well locations and geometries. Such a "universal" model
thus needs to capture the key spatio-temporal features of all these scenarios. This can be achieved by
generating a sufficiently long learning data set that contains the maximum amount of information
regarding possible modes of the model, or projections of the solutions on the POD basis.
In order to generate such a data set, a simulation in which the orientation (azimuth) of the horizontal producer
was randomly changed at regular time intervals was recorded. The injection pressure at each injection well was also varied randomly through the simulation.
All the pressure fields obtained from this simulation were flattened into vectors and
stacked together to form a "snapshot" matrix, where each column represents a specific
pressure field state.
Singular value decomposition was performed on this matrix \eqref{eq: svd}, and the first
$r$ columns of the resulting left singular matrix $\mathbf{U}$ were taken in order to obtain the reduced POD-basis.
In figure \ref{fig: universal_basis_components}, the first 12 components of the resulting reduced
basis are shown.
This basis was further used to formulate the reduced problem \eqref{eq: pressure_equation_reduced}.
In Figs. \ref{fig: production_rates_63_full} and
\ref{fig: production_rates_175_full} the simulated fluid production rates
for two particular producer orientations (63 and 175 degrees clockwise
from the horizontal axis) are shown.
One can observe that in both cases the production curves simulated by the reduced model
demonstrate a very close match with the solutions obtained by the full model.
However, due to the complex spatio-temporal structure of the solutions used for constructing the
POD-basis, a relatively large number of components (80-100) of the reduced basis
have to be kept in order to achieve such an accuracy. In Fig.
\ref{fig: rates_various_components} the results of modeling of the production rates with different
numbers of the reduced basis components are presented.

\begin{figure}[htbp]
    \begin{center}
        \includegraphics[width=0.6\textwidth]{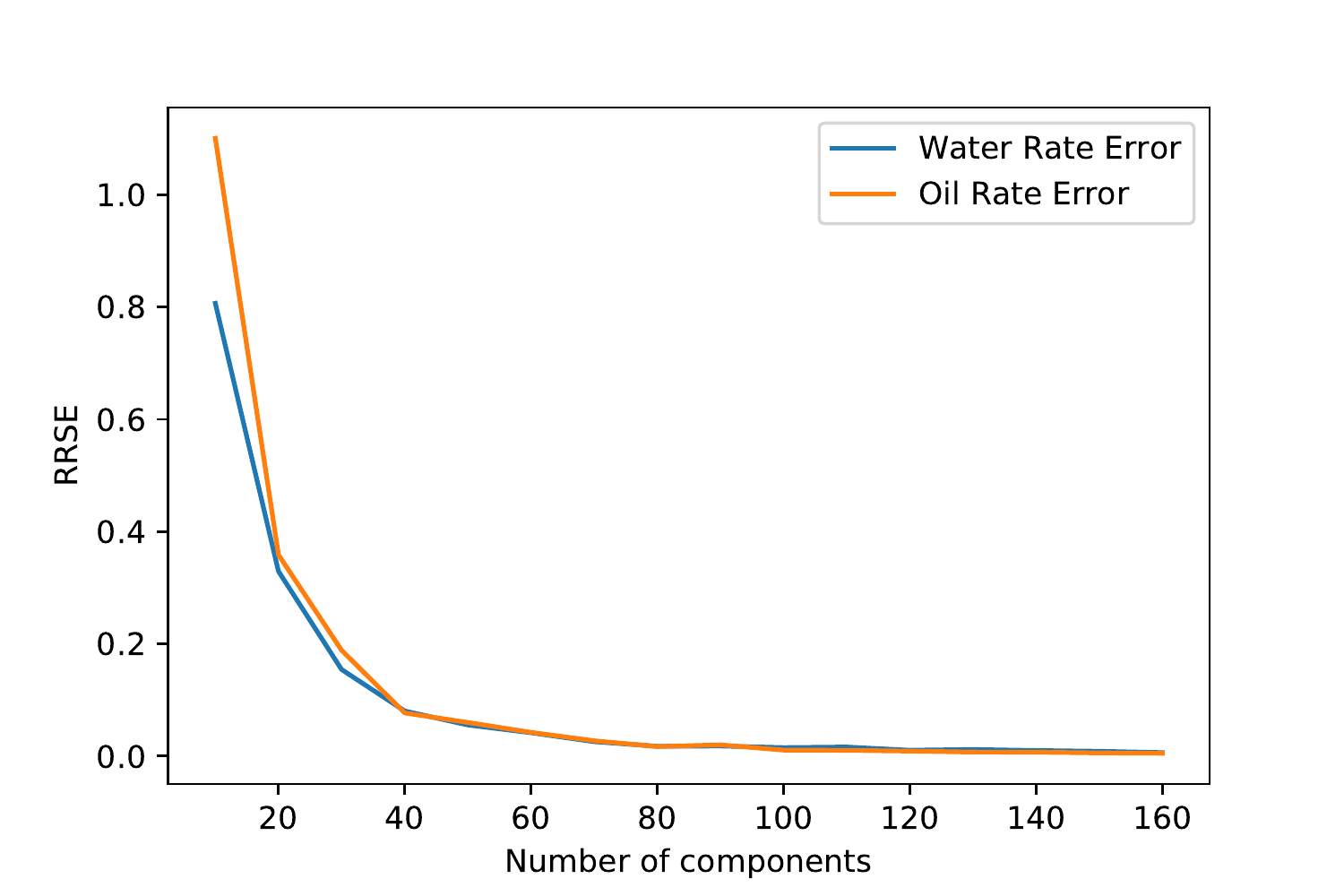}
    \end{center}
    \caption{Root relative squared error of the simulated fluid production rates as a function of number of the POD basis components.}
    \label{fig: error_components_universal_basis}
\end{figure}

Root relative squared errors (RRSE) \eqref{eq: rrse} were calculated for the deviation of the
solutions obtained with different numbers of the reduced basis components with respect to the full model solution.
RRSE is designed to be relative to an error of a simple predictor (constant mean value) \cite{ian_h_witten_data_2011}.

\begin{equation}\label{eq: rrse}
    E = \sqrt{\frac{\sum\limits_{i=1}^n(\hat{q}_i - q_i)^2}{\sum\limits_{i=1}^n(q_i - \overline{q}_i)^2}}\;,
\end{equation}
where $\hat{q}$ - production rate calculated by ROM, $q$ - reference production rate ( from full-scale model),
$\overline{q}$ - mean value of reference production rate. Errors were calculated separately for oil and water rates.
RRSE for oil and water production rates as a function of number of the POD basis components are
shown in Fig. \ref{fig: error_components_universal_basis}.

One can observe that a POD basis containing at least 80 components should be used in order
to accurately simulate the production rates.

\subsection{Adaptive POD basis}
\subsubsection{Local POD basis}
\begin{figure}[htbp]
    \begin{center}
        \includegraphics[width=100mm]{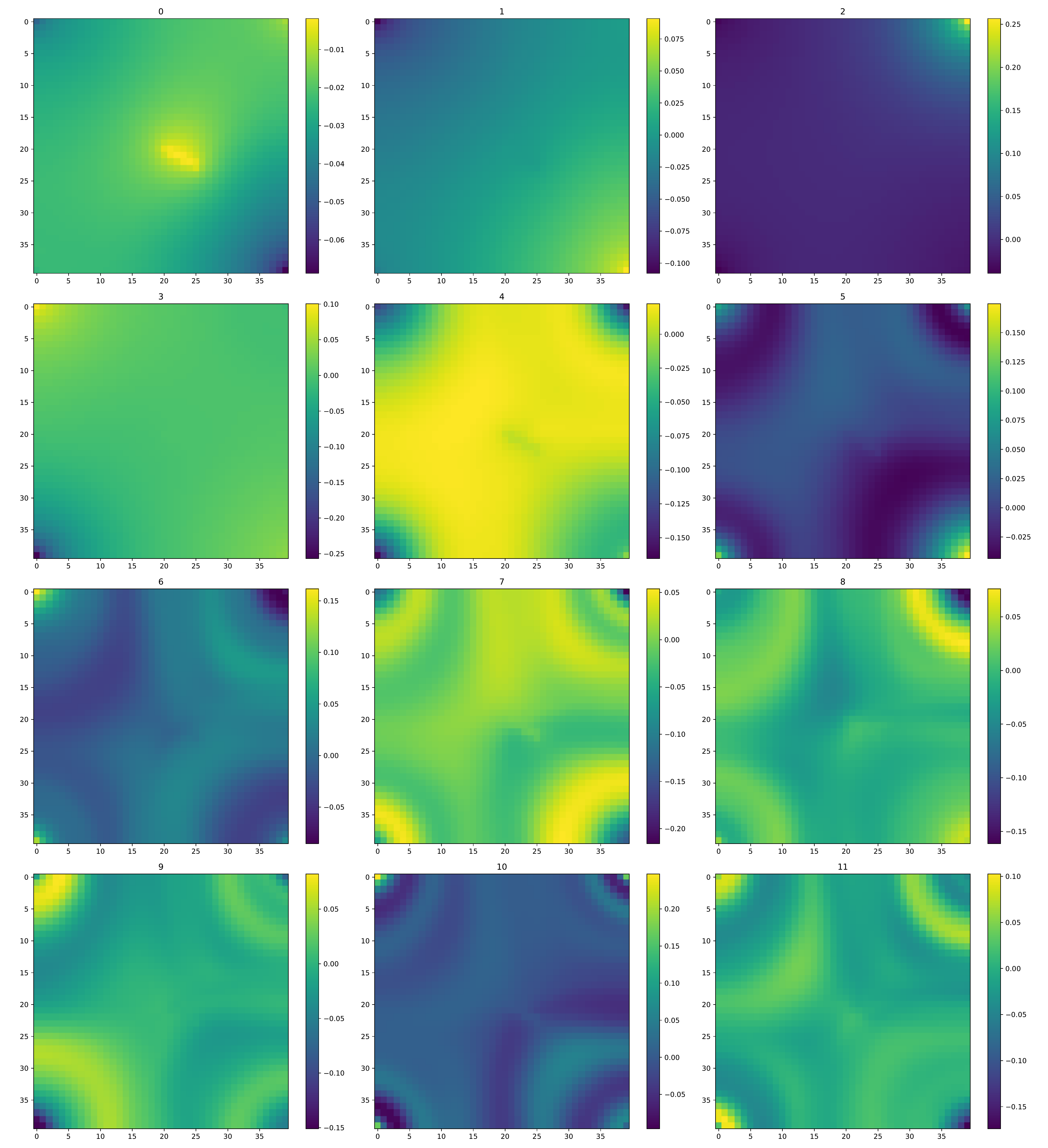}
    \end{center}
    \caption{The first 12 principal components of the local POD basis.}
    \label{fig: local_63_basis_components}
\end{figure}

Let us compare the universal basis approach with the construction of the POD-basis for a specific
well configuration. The latter is referred in this work
as the local basis approach. In this approach, a series of snapshots is generated in the same way as in the universal basis case,
however the geometry and the locations of the wells are fixed while the regimes of the injector wells are randomly changed.
These snapshot are flattened and stacked into a matrix, after which SVD is applied to that matrix, and the reduced basis is obtained.
It should be noted that such a local POD basis is only suitable for simulating scenarios with the specific well configuration for which it was constructed.
In Fig.
\ref{fig: local_63_basis_components}, the first 12 principal components of such a local POD-basis
for a model with the
producer oriented 63 degrees clockwise
from the horizontal axis are presented.

\begin{figure}[htbp]
    \begin{subfigure}{0.5\textwidth}
        \begin{center}
            \includegraphics[width=0.9\textwidth]{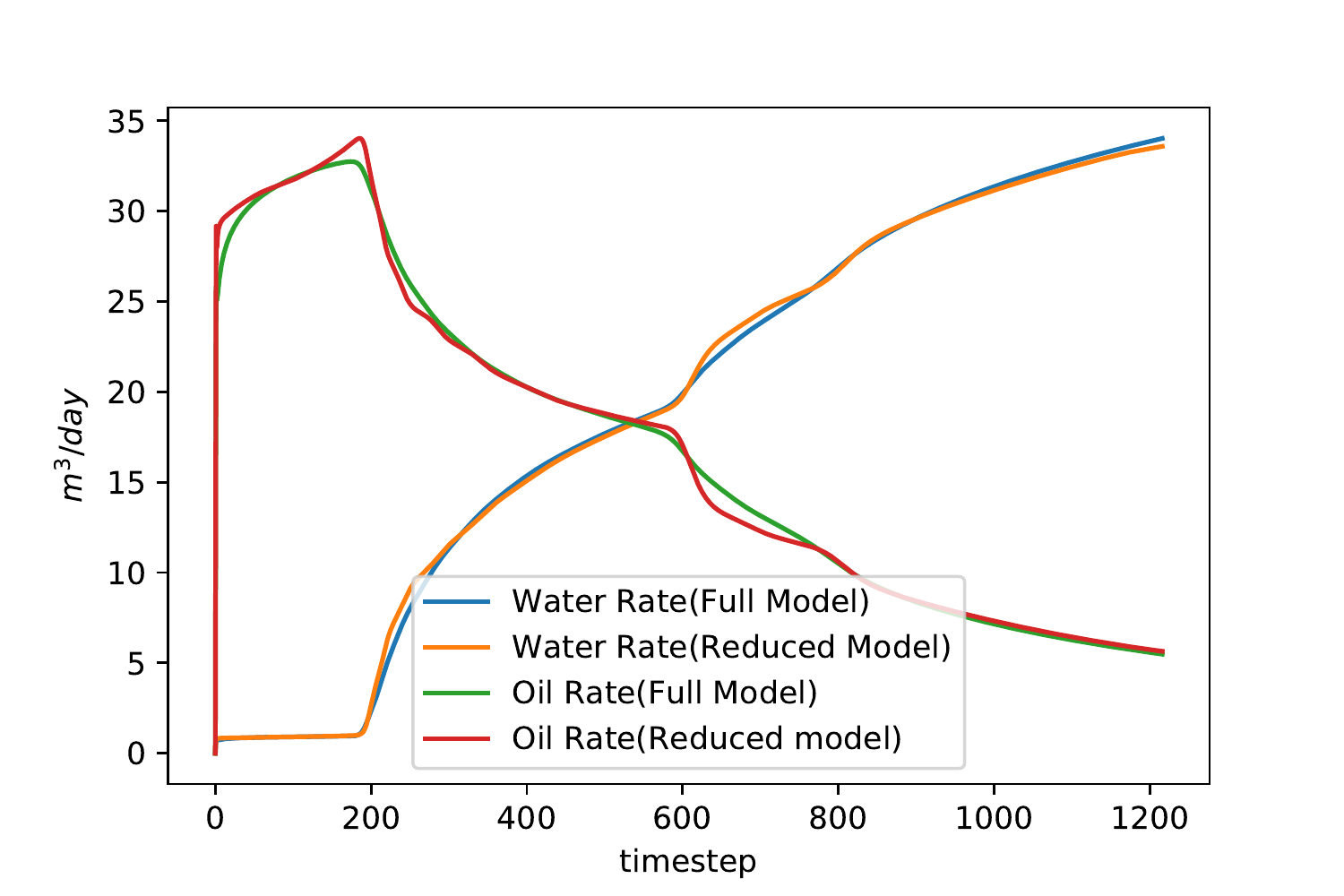}
        \end{center}
        \caption{}
        \label{fig: rates_63_10_components_local}
    \end{subfigure}
    \begin{subfigure}{0.5\textwidth}
        \begin{center}
            \includegraphics[width=0.9\textwidth]{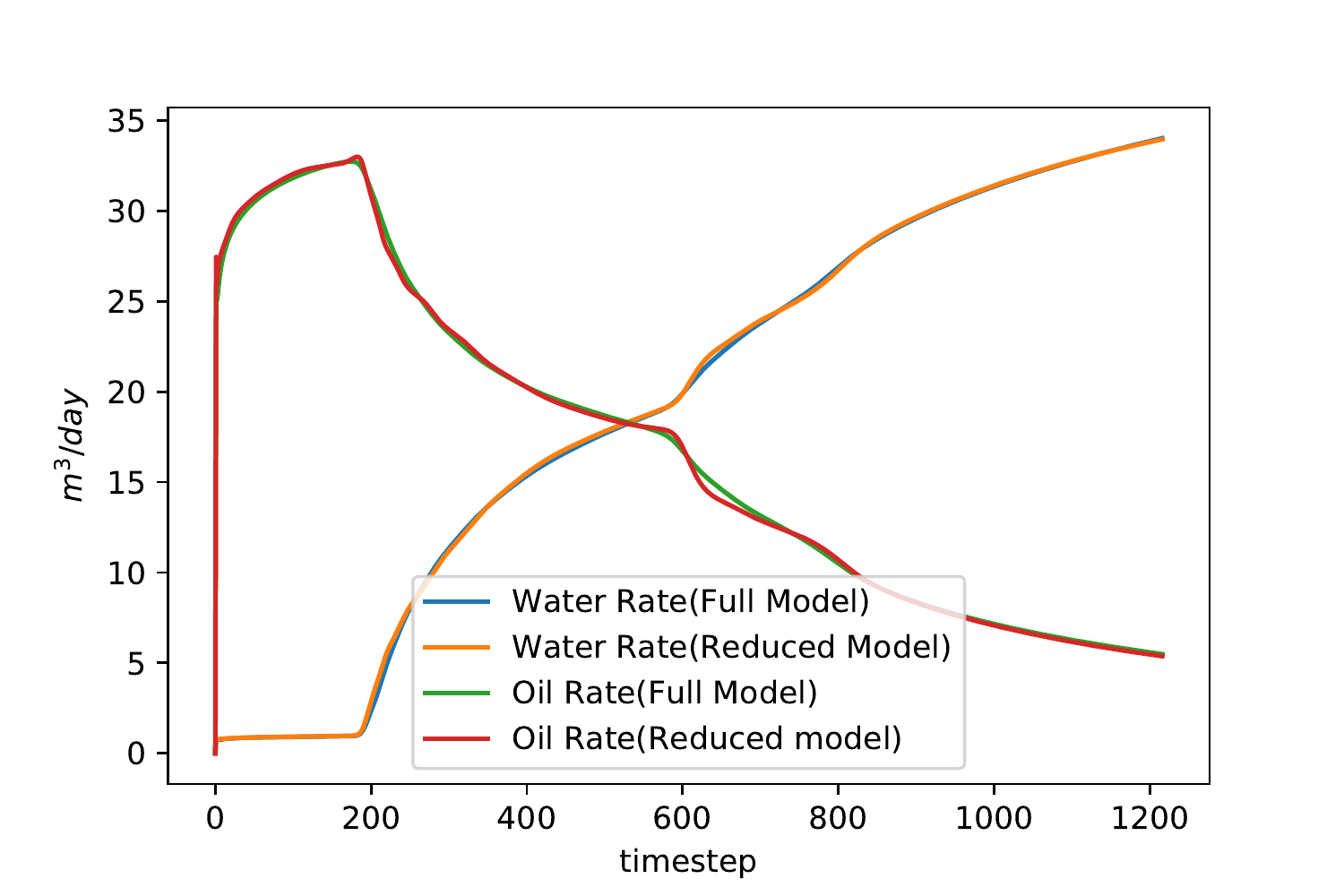}
        \end{center}
        \caption{}
        \label{fig: rates_63_20_components_local}
    \end{subfigure}
    \begin{subfigure}{0.5\textwidth}
        \begin{center}
            \includegraphics[width=0.9\textwidth]{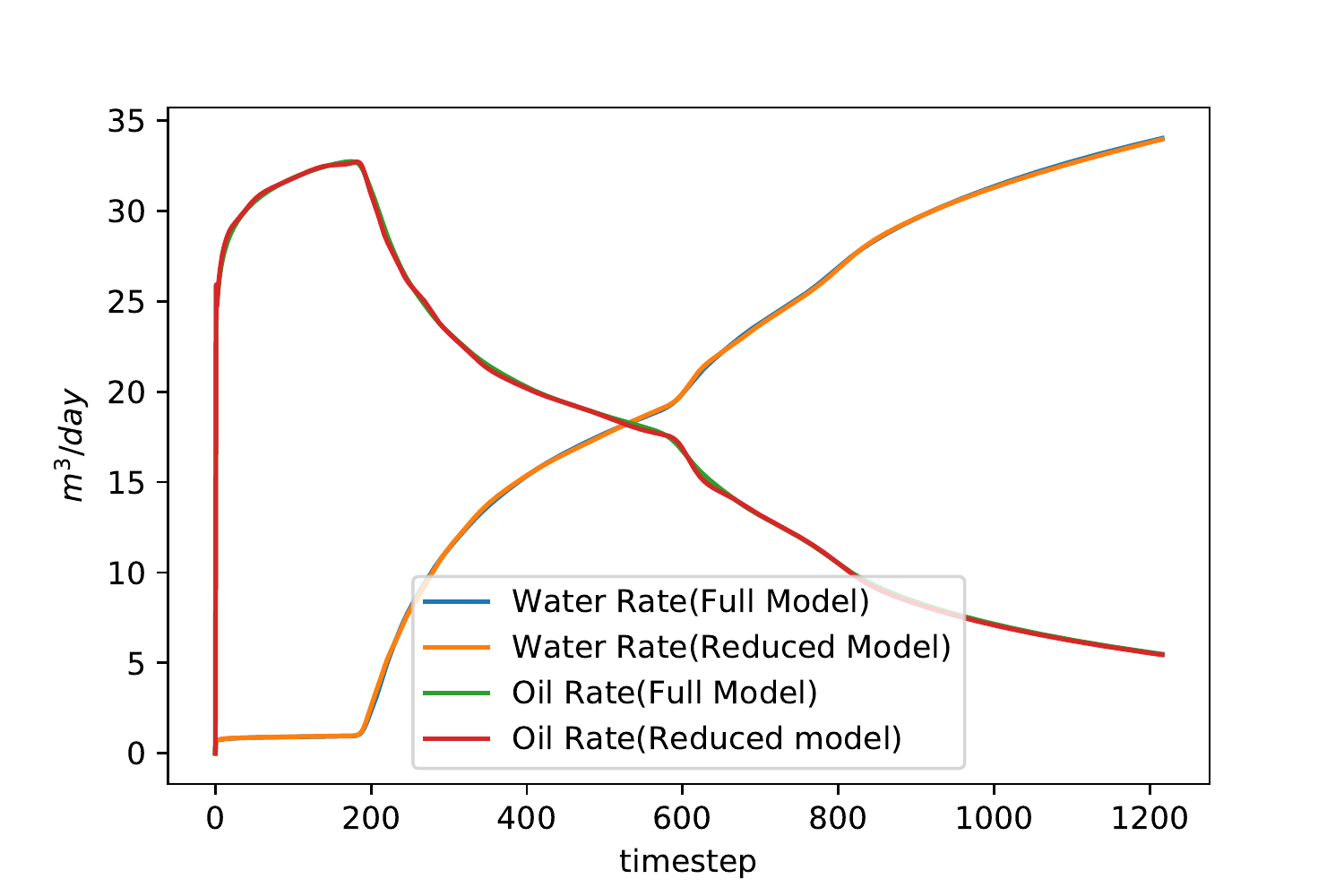}
        \end{center}
        \caption{}
        \label{fig: rates_63_30_components_local}
    \end{subfigure}
    \begin{subfigure}{0.5\textwidth}
        \begin{center}
            \includegraphics[width=0.9\textwidth]{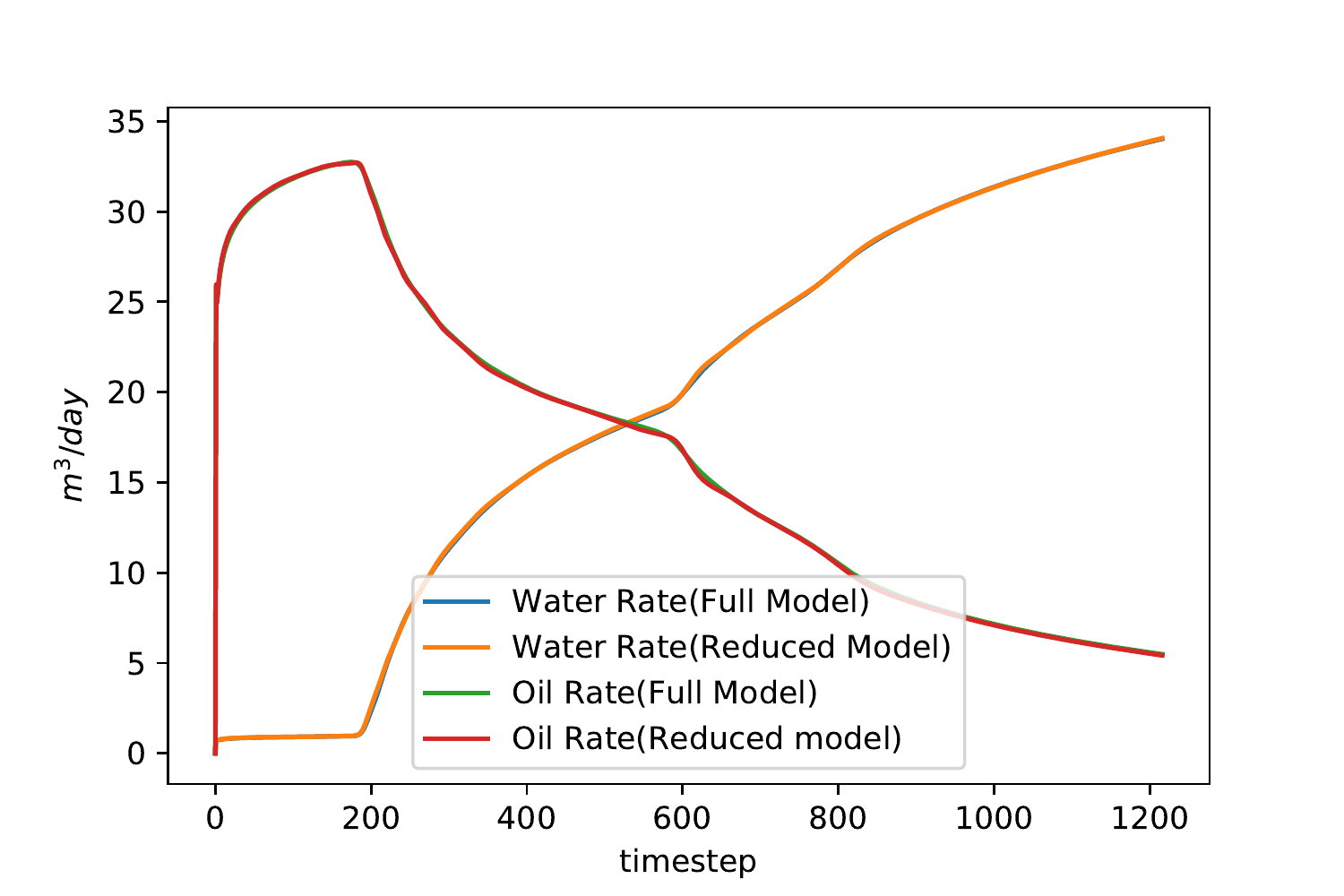}
        \end{center}
        \caption{}
        \label{fig: rates_63_35_components_local}
    \end{subfigure}
    \caption{Production rates simulated using the local POD basis with different number of components.
    \ref{sub@fig: rates_63_10_components_local} - 10 components, \ref{sub@fig: rates_63_20_components_local} -
     20 components,
    \ref{sub@fig: rates_63_30_components_local} - 30 components,
    \ref{sub@fig: rates_63_35_components_local} - 35 components.}
    \label{fig: rates_various_components_local}
\end{figure}

In Fig. \ref{fig: rates_various_components_local} the simulated production rates obtained
using different number of components of the local POD basis are shown.

\begin{figure}[htbp]
    \begin{subfigure}{0.5\textwidth}
        \begin{center}
            \includegraphics[width=0.9\textwidth]{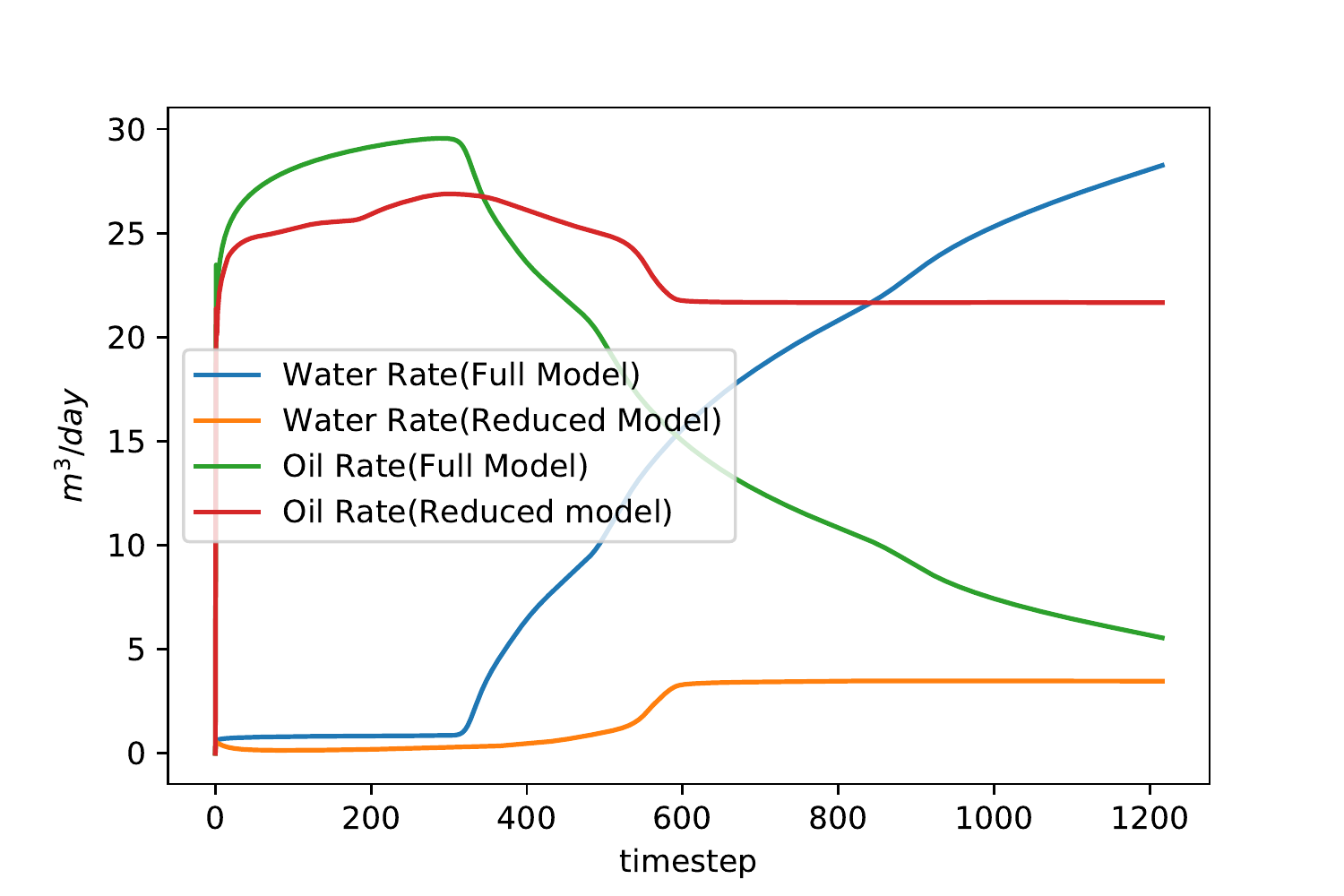}
        \end{center}
        \caption{}
        \label{fig: rates_63_175_20_components_local}
    \end{subfigure}
    \begin{subfigure}{0.5\textwidth}
        \begin{center}
            \includegraphics[width=0.9\textwidth]{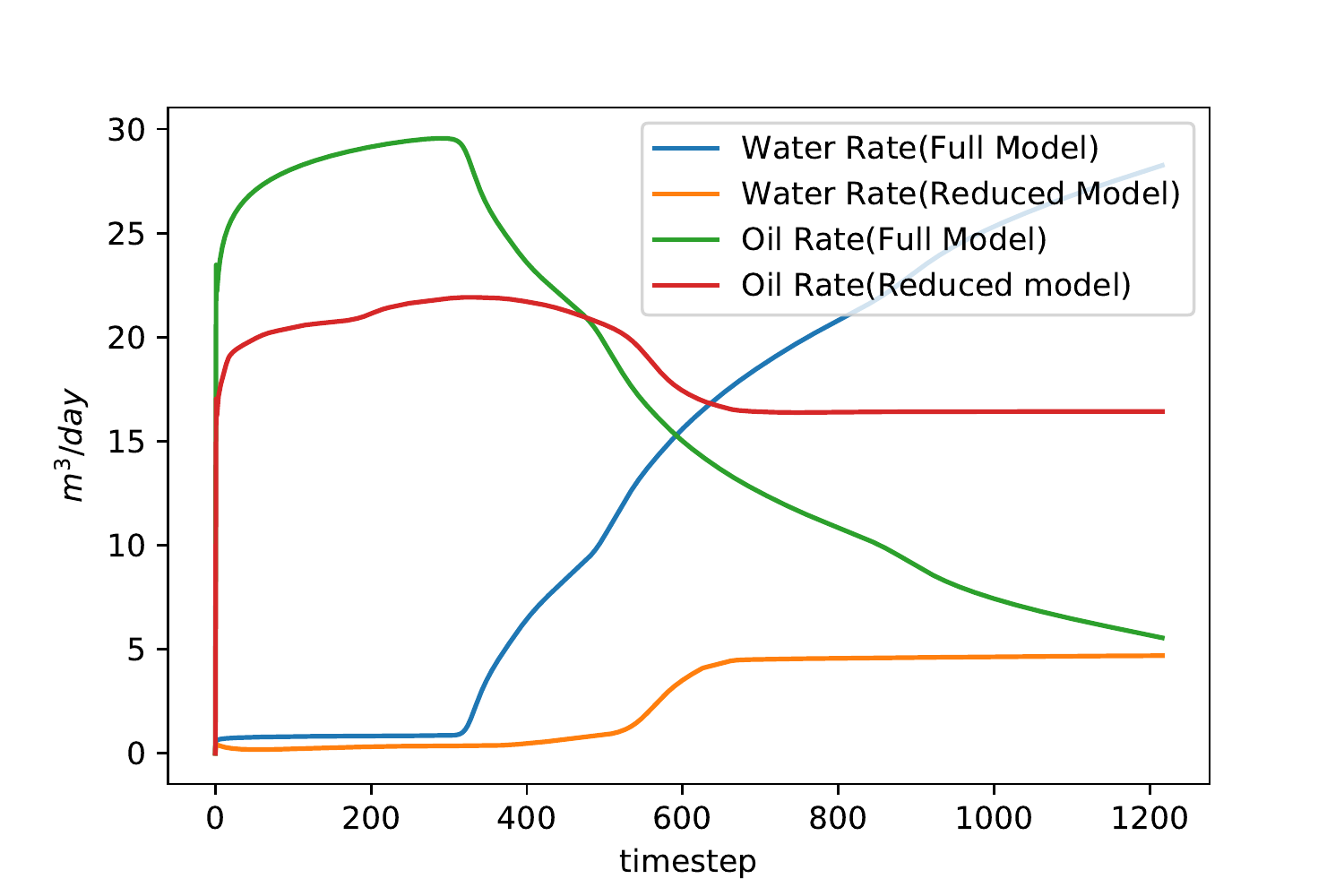}
        \end{center}
        \caption{}
        \label{fig: rates_63_175_40_components_local}
    \end{subfigure}
    \begin{subfigure}{0.5\textwidth}
        \begin{center}
            \includegraphics[width=0.9\textwidth]{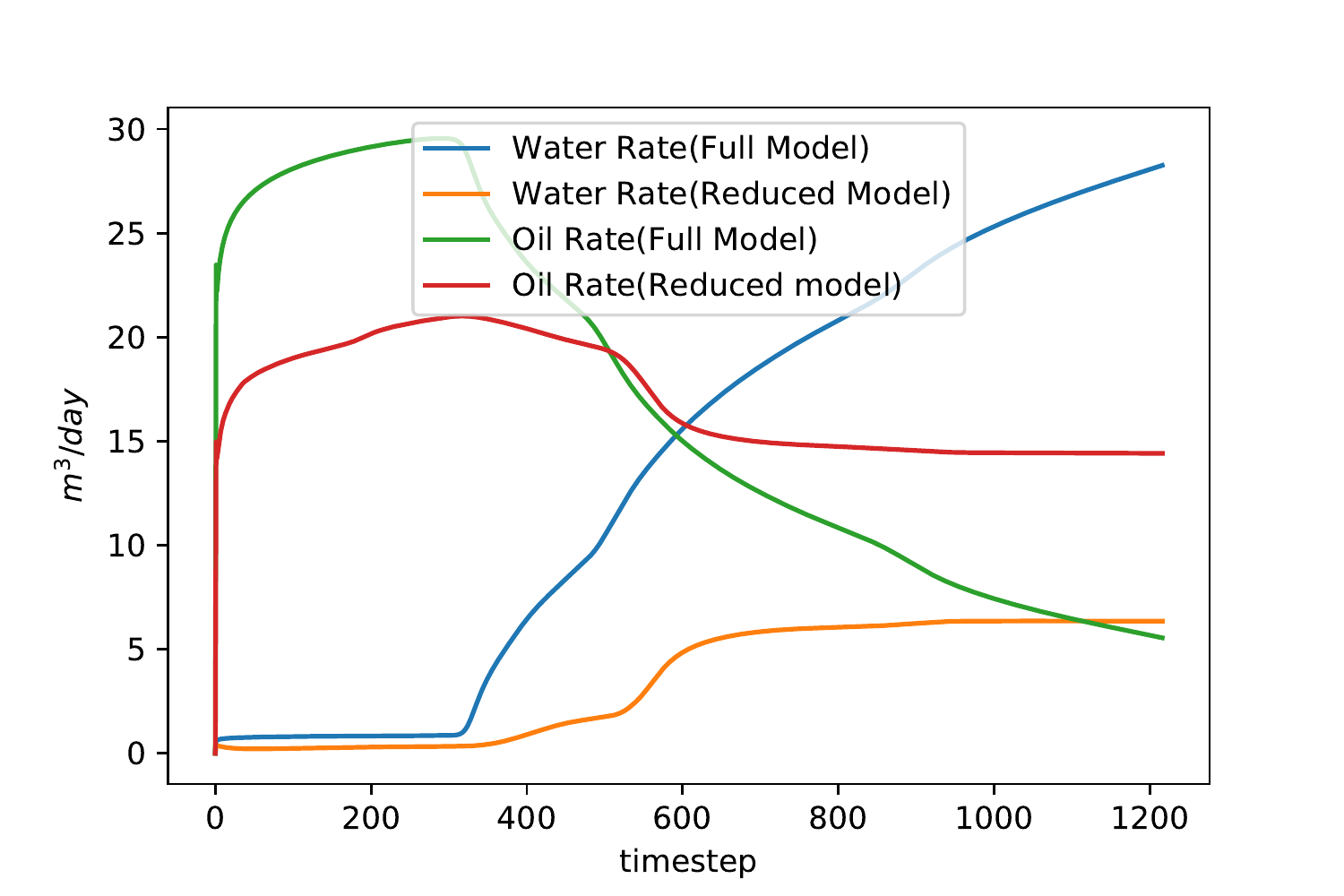}
        \end{center}
        \caption{}
        \label{fig: rates_63_175_100_components_local}
    \end{subfigure}
    \begin{subfigure}{0.5\textwidth}
        \begin{center}
            \includegraphics[width=0.9\textwidth]{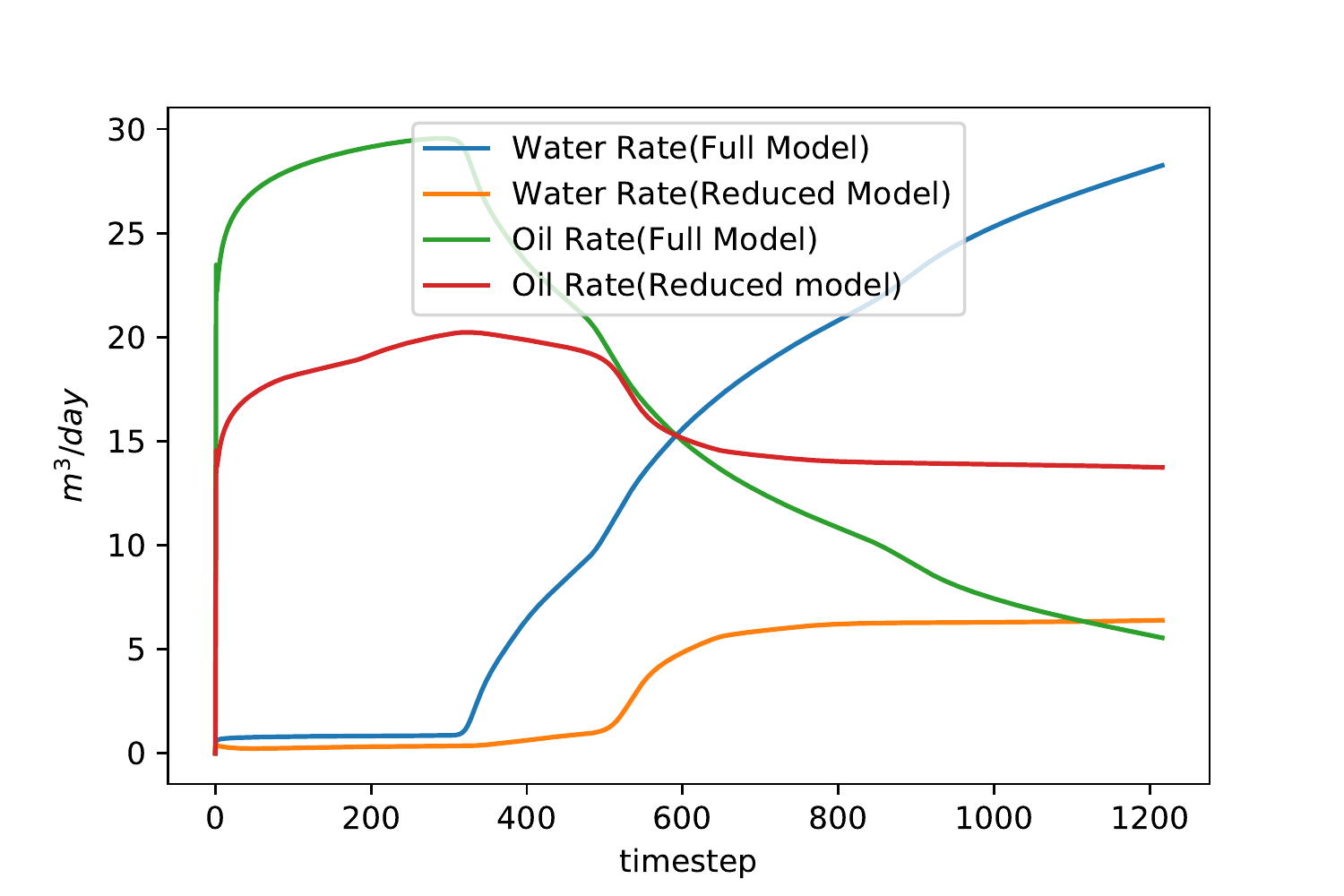}
        \end{center}
        \caption{}
        \label{fig: rates_63_175_200_components_local}
    \end{subfigure}
    \caption{Simulated production rates of the well oriented 175 degrees clockwise
    from the horizontal axis using the local POD basis constructed for a mismatching producer orientation (63 degrees) and a different number of components of the reduced basis:
    \ref{sub@fig: rates_63_175_20_components_local} - 20 components, \ref{sub@fig: rates_63_175_40_components_local} -
     40 components,
    \ref{sub@fig: rates_63_175_100_components_local} - 100 components,
    \ref{sub@fig: rates_63_175_200_components_local} - 200 components.}
    \label{fig: rates_various_components_local_63_175}
\end{figure}

One can observe that with the local POD-basis, a similar accuracy of the simulations is obtained
with fewer POD basis components.
For example, simulations with 20 components of the local POD-basis (Fig. \ref{fig: rates_63_20_components_local}) give practically the same accuracy as those with 100 components
of the universal POD basis
(Fig. \ref{fig: rates_63_100_components}) at a significantly lower computational cost.
Another advantage of the local basis is that it is much easier to construct since it
requires significantly less additional snapshots compared to the original full basis.
However, POD-Galerkin model with a local basis is only capable of simulating scenarios with
one specific well location and geometry, and if one needs to simulate a new well configuration
it is necessary to build a new basis corresponding to that configuration.
In Fig. \ref{fig: rates_various_components_local_63_175} the simulated production rates
for a well oriented at 175 degrees obtained
with the POD basis constructed for a mismatching well orientation (63 degrees clockwise
from the horizontal axis) and with different number of the basis components are presented.

One can see that in this case, the simulated production curves do not match the full model solution,
and the quality of the simulation remains poor even when using up to 200 POD basis components.

\subsubsection{Adaptive POD basis}
Although using the local POD basis constructed for a well configuration that does not match the simulated
well configuration yields quite poor results, such a mismatching basis still contains some useful information
about the simulated problem. In this section, a new approach based on utilization of this information is described which allows us
to update the POD basis and make it properly handle simulations with
new well configurations.
In order to make the old POD basis applicable to a new problem, one needs to update it
with additional components while keeping the information about the generic features
of the model.
Let us suppose that one has the local POD basis constructed for a particular well configuration
and needs to simulate production scenarios for different well configurations.
Building a new local POD basis from scratch may require a lot of additional computations since one needs to generate a
new training data set that typically consists of thousands of simulated snapshots of the full resolution
model. Additional computing resource required to generate such a training data set may entirely offset the gains achieved due to POD model reduction and thus make its use meaningless.
The proposed approach is based on updating the existing POD basis with a few new components obtained from a limited amount
of new snapshots, and allows us to produce accurate simulations with
new well configurations using the updated POD basis.

Let us consider a reduced POD basis $\mathbf{U}^r_o$ constructed for a given well configuration,
and calculate a few additional snapshots corresponding to the new well location ${\mathbf{s}_p}_i$.
The information which is lost by projecting these snapshots onto the reduced subspace defined by
the basis $\mathbf{U}^r_o$ can be expressed as
\begin{equation}\label{eq: residual_snapshots}
    {\mathbf{s}_p}_i^{res} = {\mathbf{s}_p}_i - \mathbf{U}^r_o {\mathbf{s}_p}_i {\mathbf{U}^r_o}^\top ,
\end{equation}
or in the matrix form:
\begin{equation}\label{eq: residual_snapshots_matrix}
    {\mathbf{S}_p}^{res} = {\mathbf{S}_p} - \mathbf{U}^r_o {\mathbf{S}_p} {\mathbf{U}^r_o}^\top .
\end{equation}
One can then apply SVD to this residual snapshot
matrix in order to obtain a residual basis $\mathbf{U}^r_{res}$. By construction, this basis
will be orthogonal to $\mathbf{U}^r_o$, so that one can use a combination of components from
this basis in order to build an updated basis $\widetilde{\mathbf{U}^r}$ that can be used
in the POD-Galerkin method. The suggested method complements $\mathbf{U}^r_o$
by a few components from $\mathbf{U}^r_{res}$ in order to construct $\widetilde{\mathbf{U}^r}$
which can be used in a reduced POD-Galerkin model for the new well configuration.
The workflow of the proposed method can be summarized as follows:
\begin{enumerate}
    \item{calculate $n$ snapshots with the new well configuration;}
    \item{compose the snapshot matrix $\mathbf{S}$;}
    \item{calculate the residual snapshot matrix $\mathbf{{S}_p}^{res}$; }
    \item{perform SVD and take the first $r^{res}$ components ($r^{res}\ll r$);}
    \item{update the existing POD basis $\mathbf{U}^r_o$ with these components and obtain $\widetilde{\mathbf{U}^r}$.}
    \item{Use $\widetilde{\mathbf{U}^r}$ to formulate the updated reduced POD-Galerkin model.}
\end{enumerate}

\begin{figure}[htbp]
    \begin{subfigure}{0.3\textwidth}
        \begin{center}
            \includegraphics[width=\textwidth]{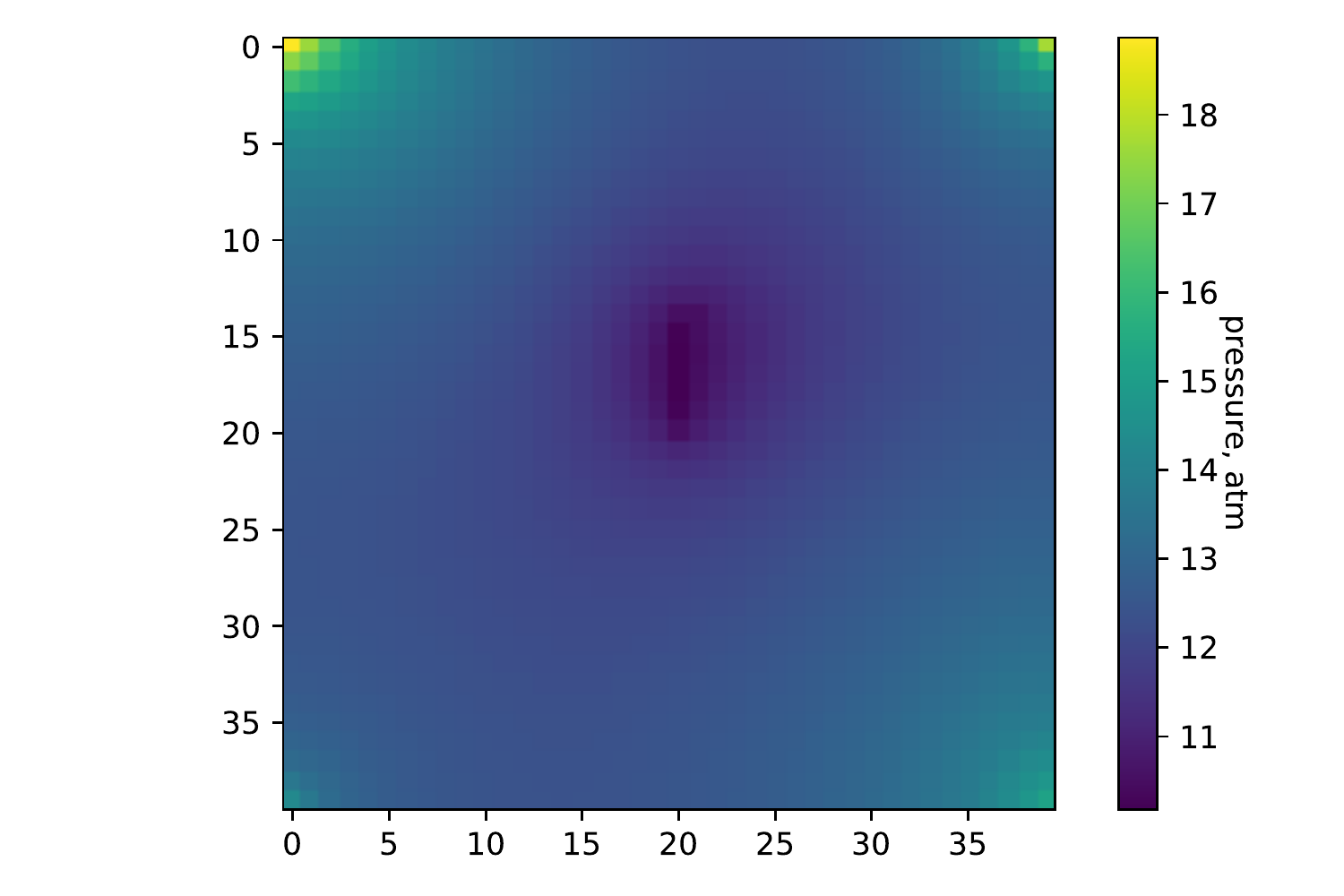}
        \end{center}
        \caption{}
        \label{fig: snapshot_projection: original_snapshot}
    \end{subfigure}
    \begin{subfigure}{0.3\textwidth}
        \begin{center}
            \includegraphics[width=\textwidth]{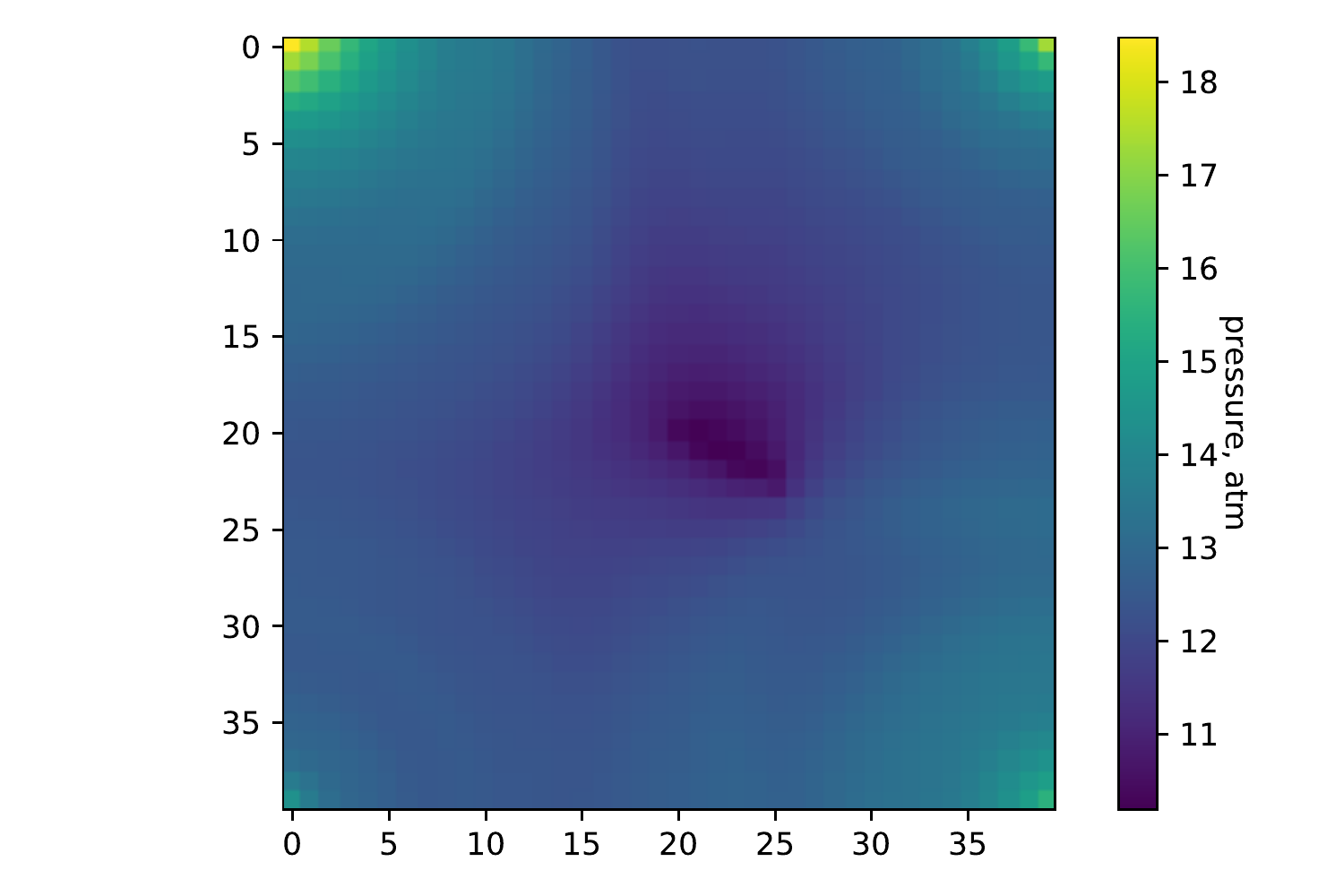}
        \end{center}
        \caption{}
        \label{fig: snapshot_projection: reduced_snapshot}
    \end{subfigure}
    \begin{subfigure}{0.3\textwidth}
        \begin{center}
            \includegraphics[width=\textwidth]{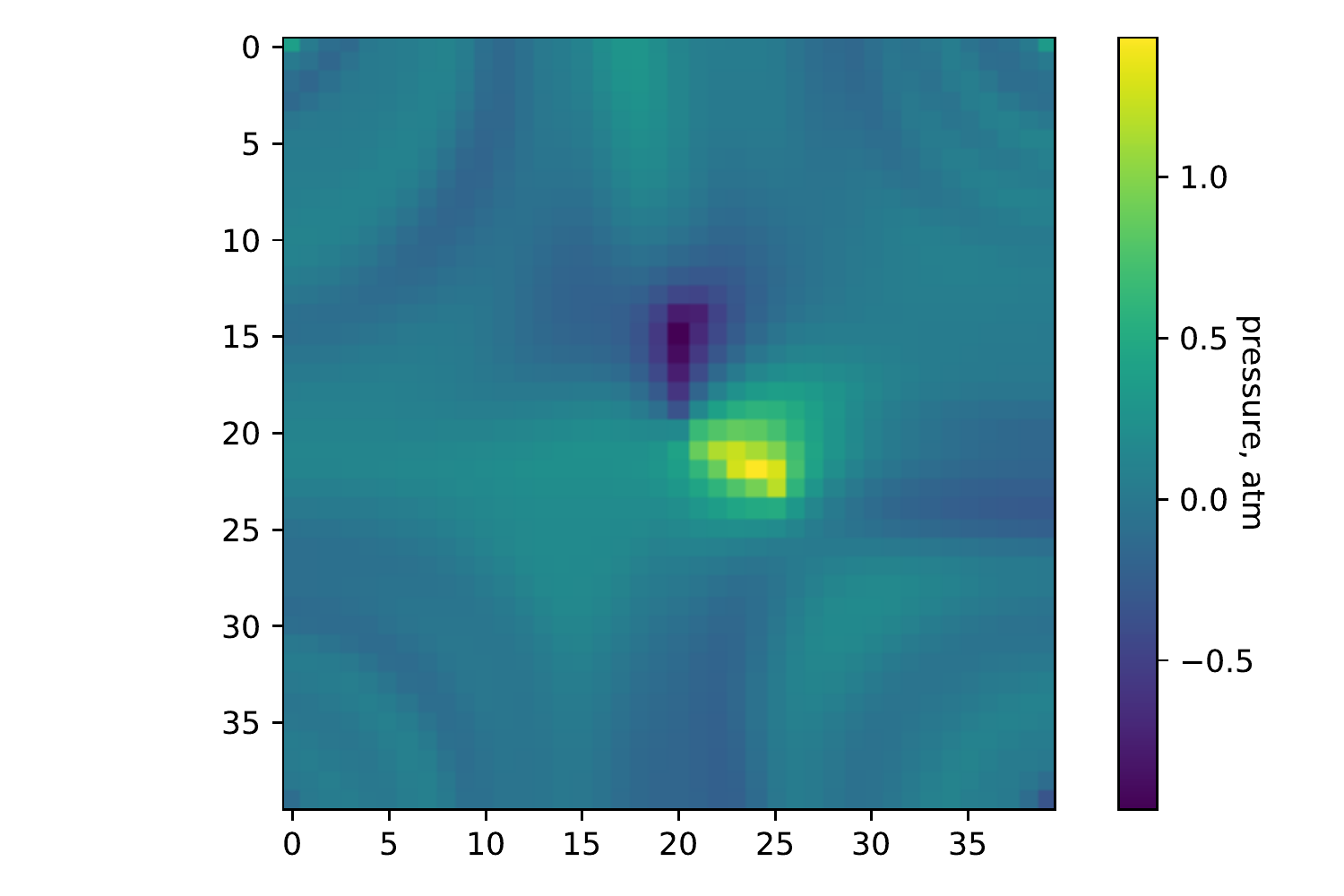}
        \end{center}
        \caption{}
        \label{fig: snapshot_projection: residual_snapshot}
    \end{subfigure}
    \caption{\ref{sub@fig: snapshot_projection: original_snapshot} - original snapshot  (${\mathbf{s}_p}_i$),
             \ref{sub@fig: snapshot_projection: reduced_snapshot} -reduced snapshot
             ($\mathbf{U}^r_o {\mathbf{s}_p}_i {\mathbf{U}^r_o}^\top$),
             \ref{sub@fig: snapshot_projection: residual_snapshot} -residual snapshot ($\mathbf{{S}_p}^{res}$).}
    \label{fig: snapshot_projection}
\end{figure}

\begin{figure}[htbp]
    \begin{center}
        \includegraphics[width=100mm]{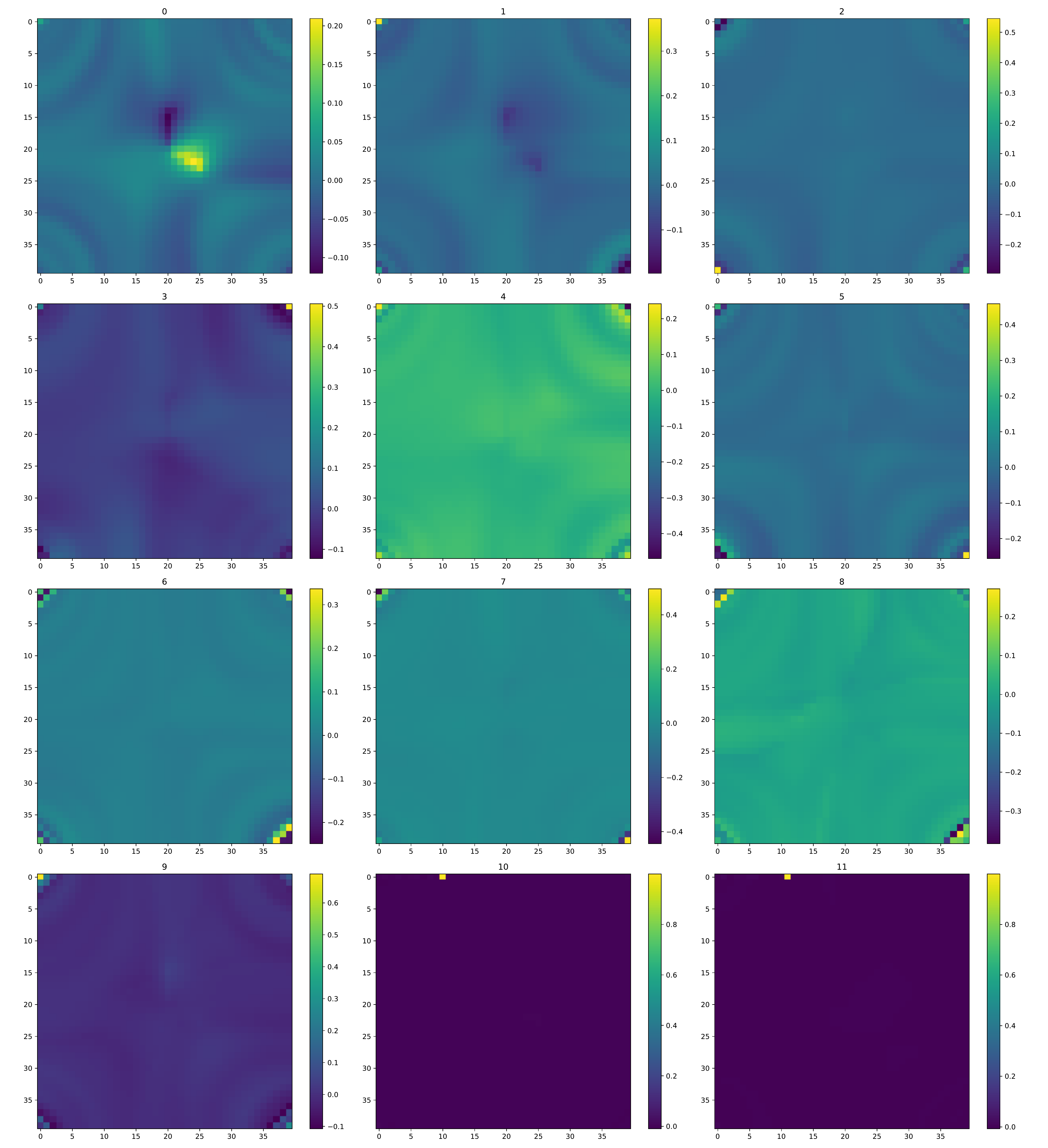}
    \end{center}
    \caption{The first 12 principal components of the residual snapshot matrix decomposition ($\mathbf{U}^r_{res}$).}
    \label{fig: residual_basis_components}
\end{figure}

\begin{figure}[htbp]
    \begin{subfigure}{0.3\textwidth}
        \begin{center}
            \includegraphics[width=\textwidth]{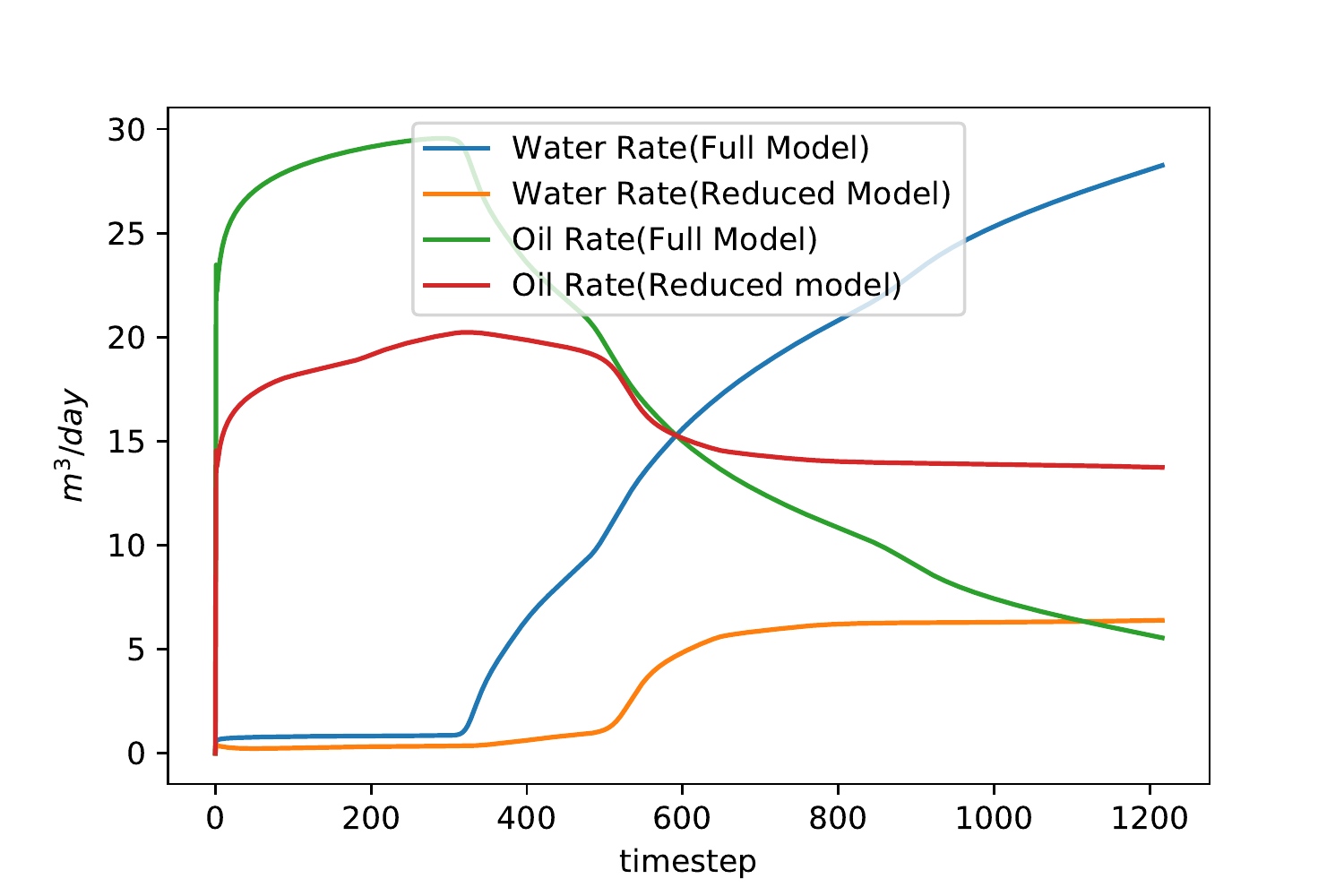}
        \end{center}
        \caption{}
        \label{fig: adaptive_basis_rates: 63_basis}
    \end{subfigure}
    \begin{subfigure}{0.3\textwidth}
        \begin{center}
            \includegraphics[width=\textwidth]{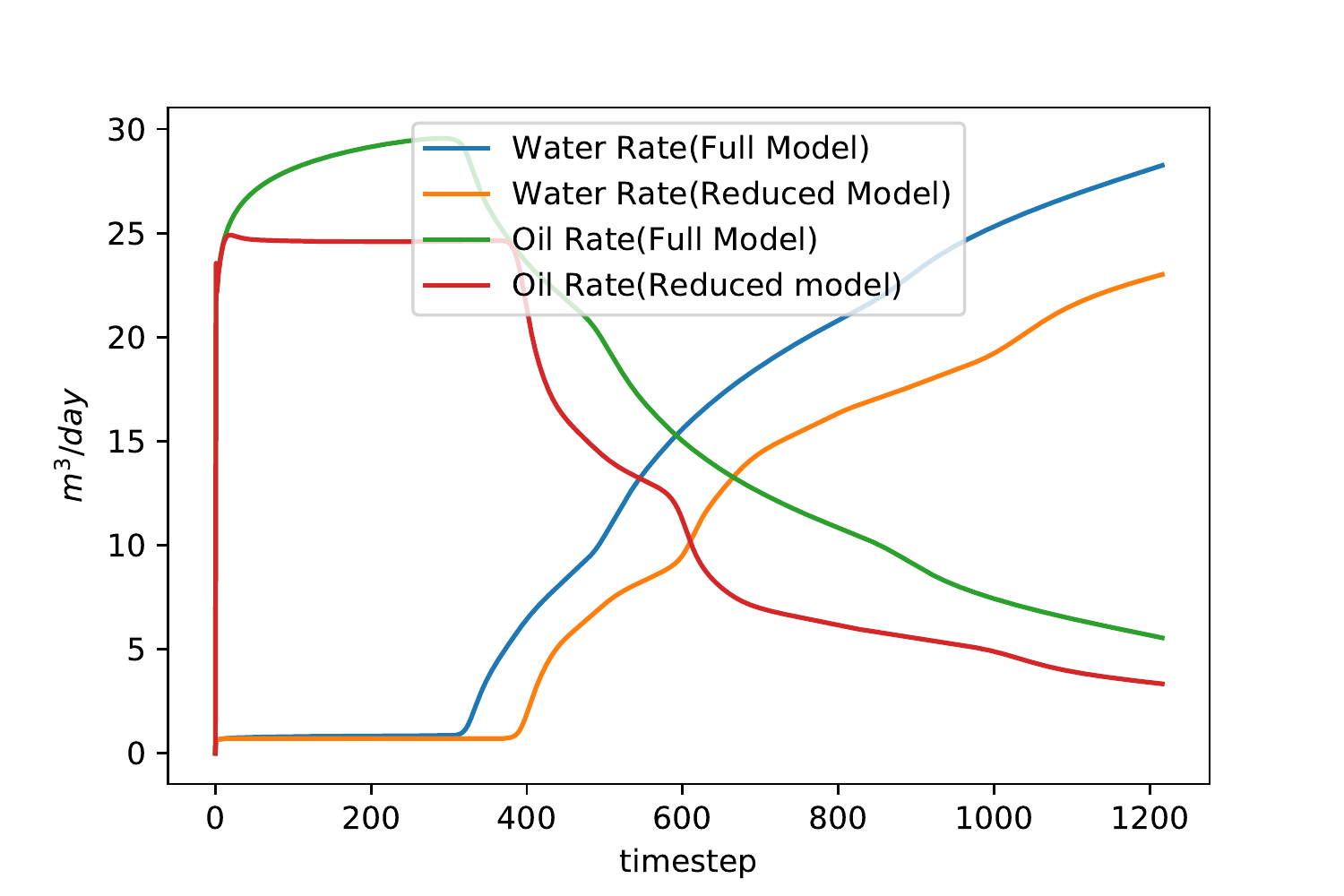}
        \end{center}
        \caption{}
        \label{fig: adaptive_basis_rates: local_basis_10_components}
    \end{subfigure}
    \begin{subfigure}{0.3\textwidth}
        \begin{center}
            \includegraphics[width=\textwidth]{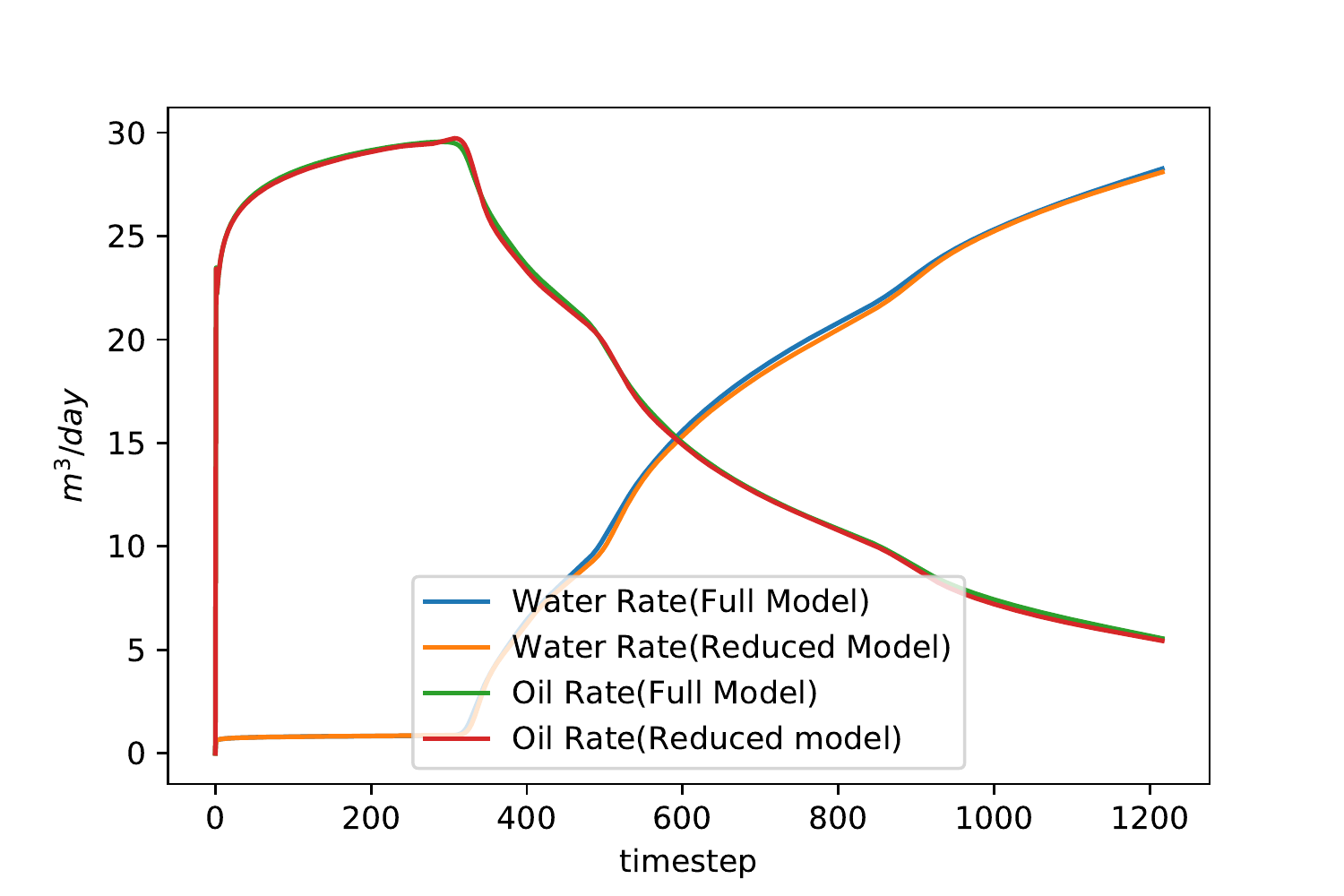}
        \end{center}
        \caption{}
        \label{fig: adaptive_basis_rates: adaptive_basis}
    \end{subfigure}
    \caption{Simulated fluid production rates:
             \ref{sub@fig: adaptive_basis_rates: 63_basis} - using the local basis for a mismatching well
             orientation;
             \ref{sub@fig: adaptive_basis_rates: local_basis_10_components} - using the updated POD
             basis constructed from 10 additional snapshots;
             \ref{sub@fig: adaptive_basis_rates: adaptive_basis} - using the adaptive POD basis.}
    \label{fig: adaptive_basis_rates}
\end{figure}

Now let us consider an example of the application of the proposed method. As the initial basis, $\mathbf{U}^r_o$
a local basis consisting of 20 components and corresponding to the
well oriented
63 degrees clockwise from the horizontal axis (Fig. \ref{fig: local_63_basis_components}) is taken.
Only 10 additional snapshots for the new well direction
(175 degrees clockwise from the horizontal axis) are simulated. As a reminder:
construction from scratch of a new POD basis for that problem would require generating about a thousand snapshots.
The snapshot matrix ${\mathbf{S}_p}$ is then composed, and the residual snapshot matrix $\mathbf{{S}_p}^{res}$ is obtained
using \eqref{eq: residual_snapshots_matrix}.
In Fig. \ref{fig: snapshot_projection} one of such snapshots,
the corresponding reduced snapshot ($\mathbf{U}^r_o {\mathbf{s}_p}_i {\mathbf{U}^r_o}^\top$),
and the residual snapshot  \eqref{eq: residual_snapshots} are shown.
SVD is then performed on the residual snapshot matrix, and the additional components are obtained.
The first 12 of the resulting additional components are presented in the Fig. \ref{fig: residual_basis_components}.
Since only 10 additional snapshots were used, most of the higher order residual components
are quite noisy and contain little structural information about the model.
Consequently, only the first three residual components will be used
to build the updated basis $\widetilde{\mathbf{U}^r}$ by adding them
to the original basis $\mathbf{U}^r_o$. One can then use the updated basis to formulate the
POD-Galerkin problem for the new well configuration.
The simulated production rates obtained for the new well configuration with the updated POD basis
are shown in Fig. \ref{fig: adaptive_basis_rates}.
One can observe that the proposed method of the local basis adaptation
significantly increases the accuracy of the simulations. In Fig.
\ref{fig: adaptive_basis_rates: local_basis_10_components}, the results of
simulation of the production rates using the local basis with 23 components
(the same number as in the adaptive basis) built using the same 10 additional snapshots, are presented.
One can conclude from these simulations that 10 additional snapshots of the model are not sufficient to construct a new basis without using the information contained in original basis. However, the
adaptive POD-basis scheme with the same 10 additional snapshots provides quite a
satisfactory result.

\begin{figure}[htbp]
    \begin{subfigure}{0.5\textwidth}
        \begin{center}
            \includegraphics[width=0.9\textwidth]{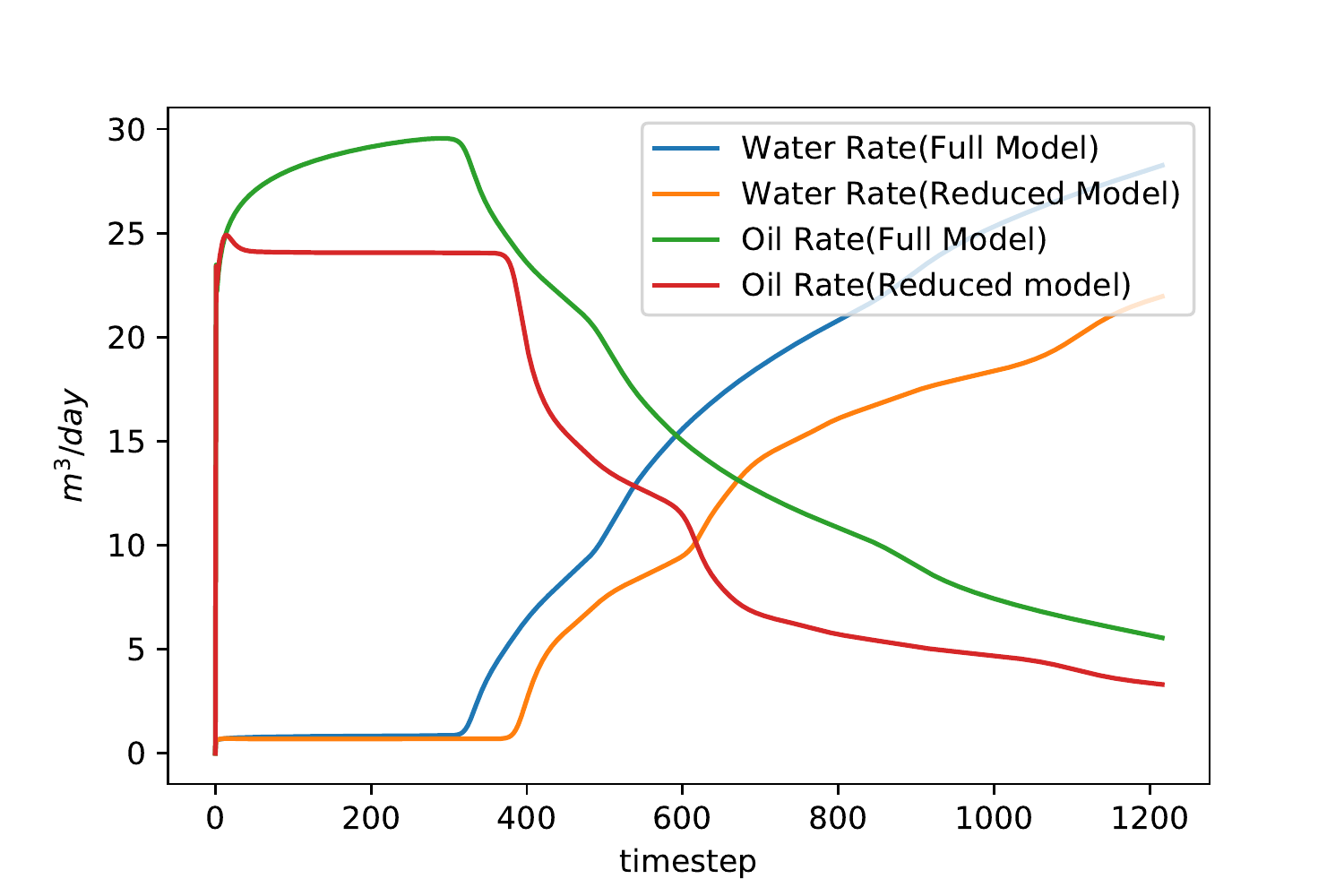}
        \end{center}
        \caption{}
        \label{fig: snapshot_number: 10_snapshots}
    \end{subfigure}
    \begin{subfigure}{0.5\textwidth}
        \begin{center}
            \includegraphics[width=0.9\textwidth]{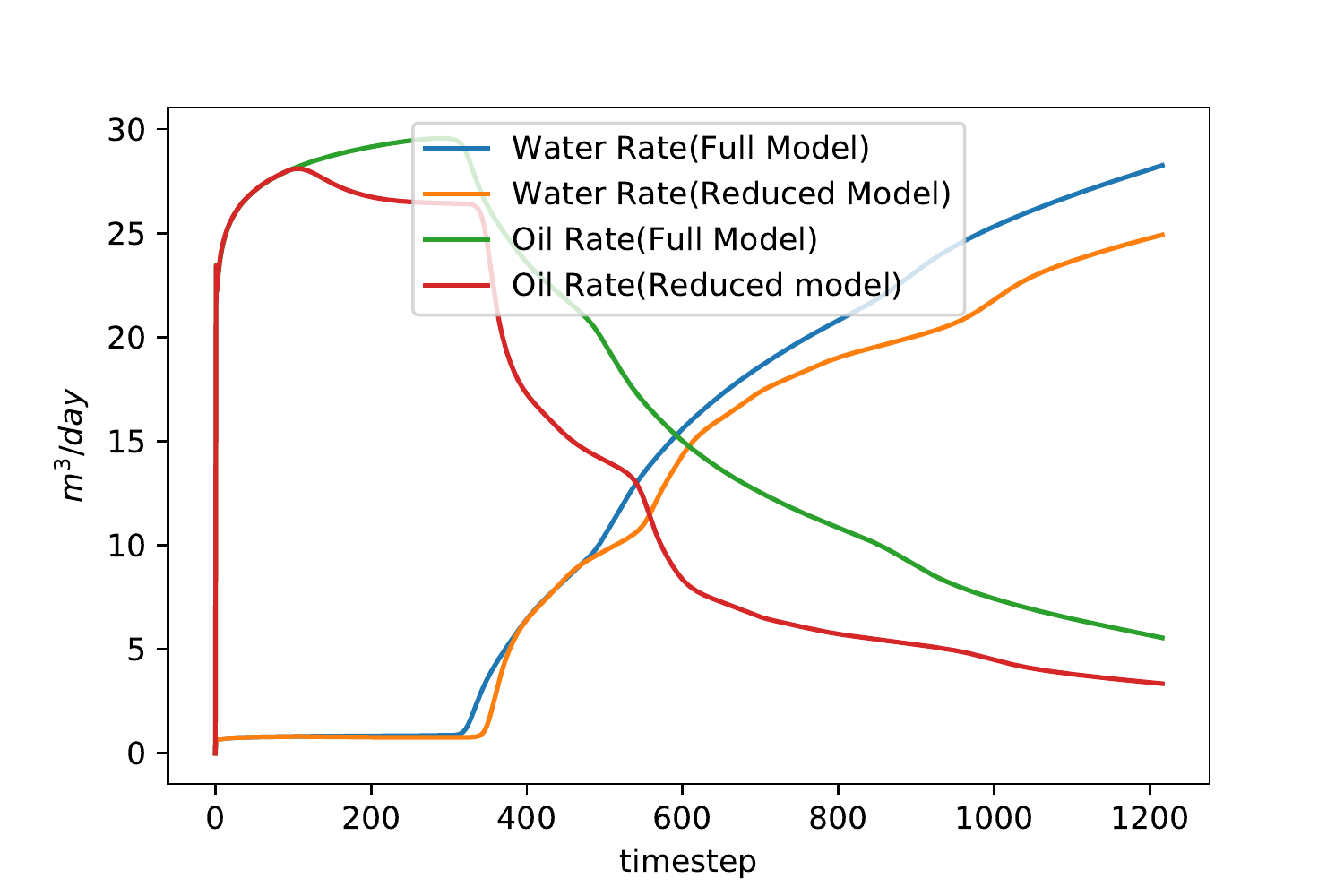}
        \end{center}
        \caption{}
        \label{fig: snapshot_number: 100_snapshots}
    \end{subfigure}
    \begin{subfigure}{0.5\textwidth}
        \begin{center}
            \includegraphics[width=0.9\textwidth]{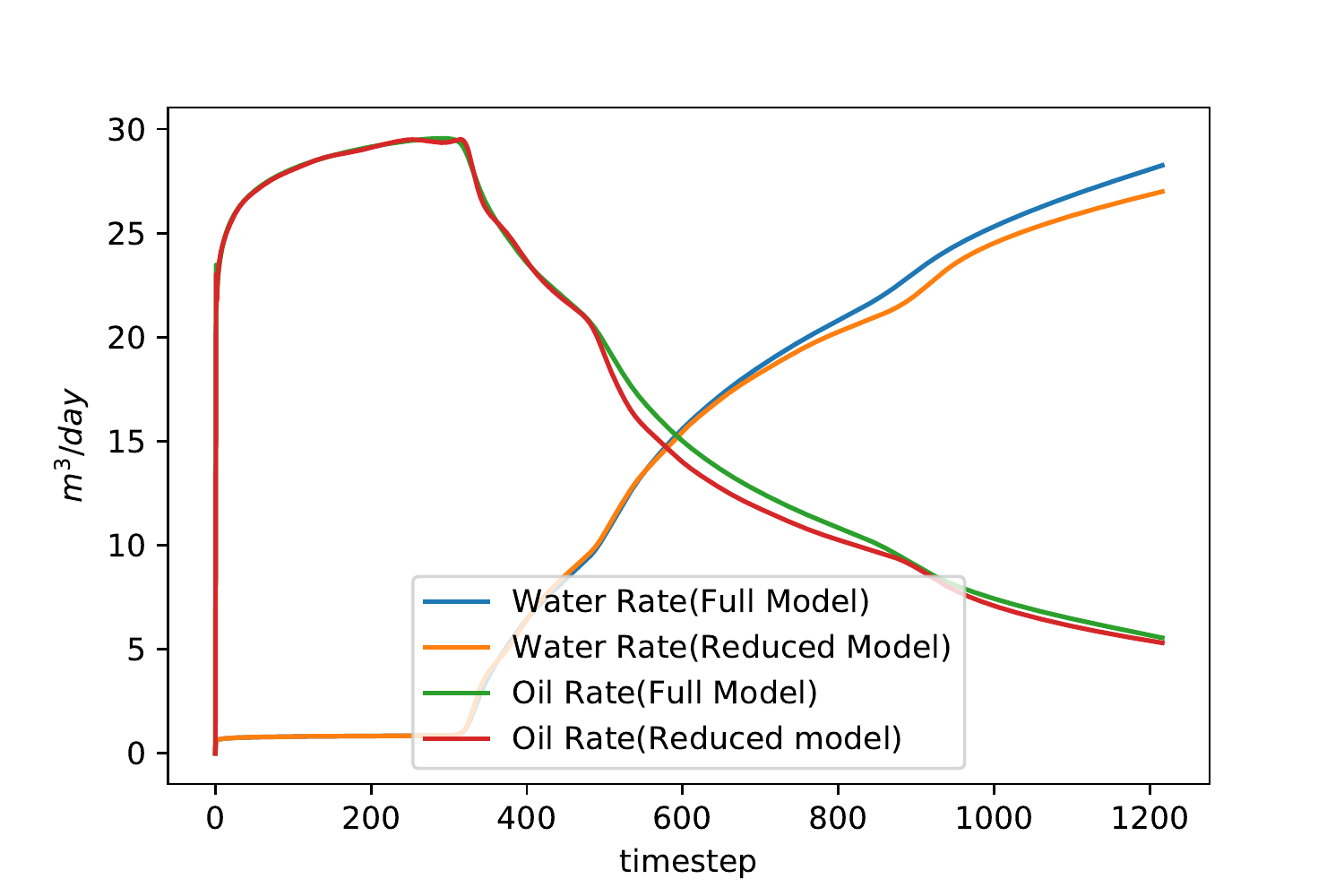}
        \end{center}
        \caption{}
        \label{fig: snapshot_number: 500_snapshots}
    \end{subfigure}
    \begin{subfigure}{0.5\textwidth}
        \begin{center}
            \includegraphics[width=0.9\textwidth]{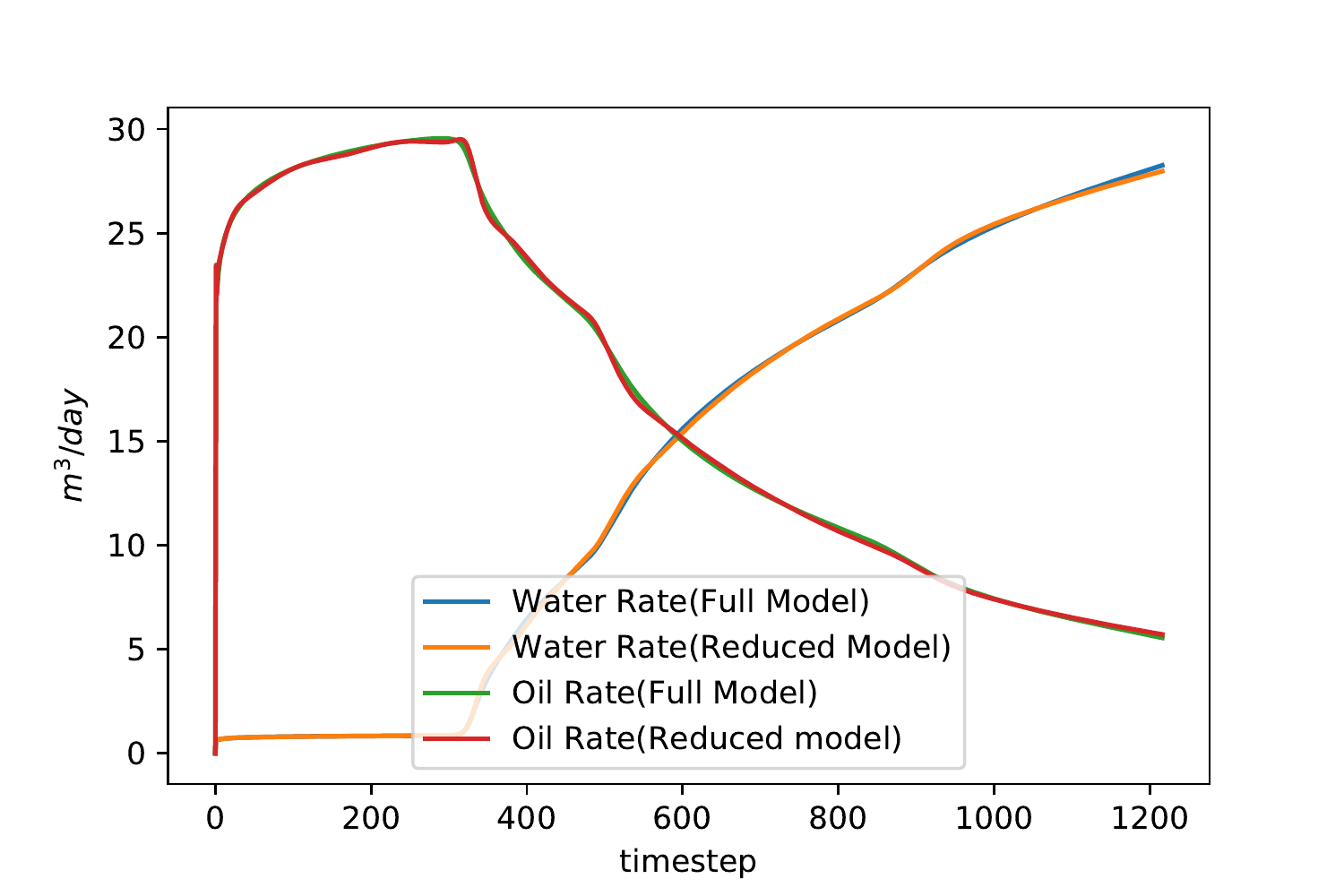}
        \end{center}
        \caption{}
        \label{fig: snapshot_number: 1000_snapshots}
    \end{subfigure}
    \caption{Production rates simulated using POD-Galerkin models with
    local bases constructed from different number of snapshots:
    \ref{sub@fig: snapshot_number: 10_snapshots} - 10 snapshots;
    \ref{sub@fig: snapshot_number: 100_snapshots} - 100 snapshots;
    \ref{sub@fig: snapshot_number: 500_snapshots} - 500 snapshots;
    \ref{sub@fig: snapshot_number: 1000_snapshots} - 1000 snapshots.}
    \label{fig: snapshot_number}
\end{figure}

In order to estimate the number of snapshots required to build the local basis from scratch,
a number of simulations with POD-Galerkin models using bases constructed from
varying numbers of snapshots were performed.
The corresponding simulated production rates are shown in Fig. \ref{fig: snapshot_number}.
One can observe that in order to construct a proper local basis, about 1000 snapshots are required,
in contrast with just about 10 additional snapshots needed for the adaptation of the existing basis.

\subsubsection{Adaptive POD bases for models with variations of well location and well length}
\begin{figure}[htbp]
    \begin{center}
        \includegraphics[width=0.6\textwidth]{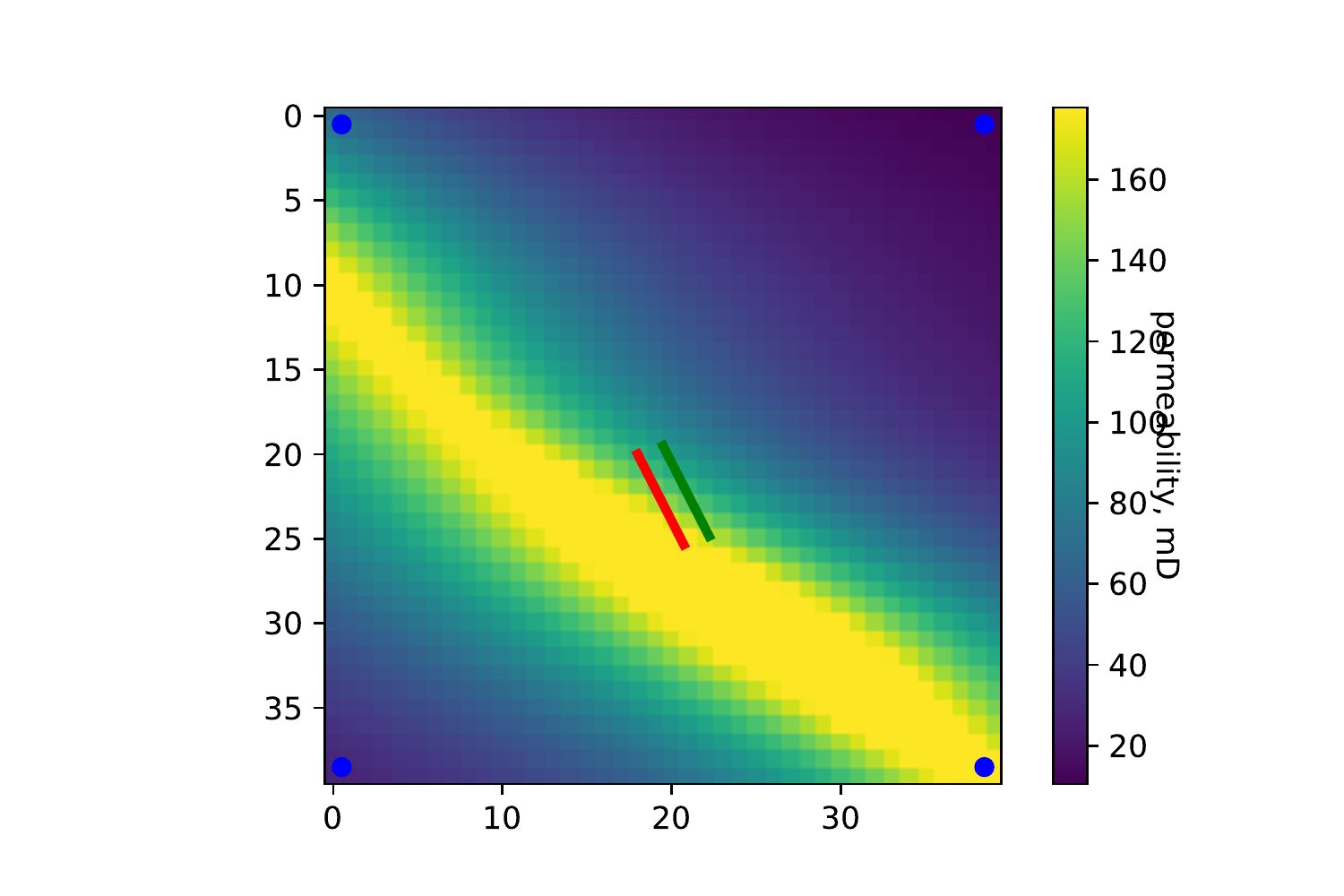}
    \end{center}
    \caption{Well placement scheme: blue circles represent the injector wells; the green line
    represents the original position of the producer; the red line represents the new
    position of the producer.}
    \label{well_map_18_20}
\end{figure}

\begin{figure}[htbp]
    \begin{subfigure}{0.49\textwidth}
        \begin{center}
            \includegraphics[width=\textwidth]{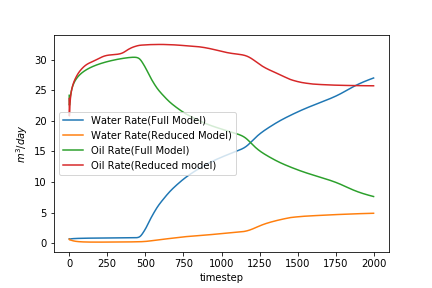}
        \end{center}
        \caption{}
        \label{fig: position_rates_calculation: original_basis}
    \end{subfigure}
    \begin{subfigure}{0.49\textwidth}
        \begin{center}
            \includegraphics[width=\textwidth]{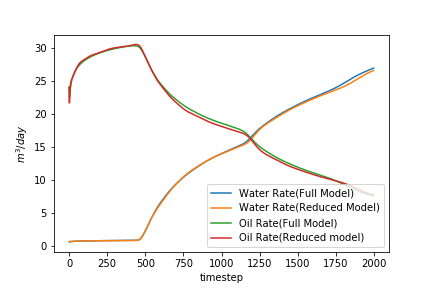}
        \end{center}
        \caption{}
        \label{fig: position_rates_calculation: adaptive_basis}
    \end{subfigure}
    \caption{Simulated fluid production rates: \ref{sub@fig: position_rates_calculation: original_basis} -
    the original POD basis, \ref{sub@fig: position_rates_calculation: adaptive_basis} - the adaptive POD basis constructed using 50 new snapshots and 1 additional component.}
    \label{fig: position_rates_calculation}
\end{figure}

In order to further test the applicability of the proposed technique to problems with changing well configurations,  it was used to build the adaptive bases for models with variations
of well location and well length.
In Fig. \ref{well_map_18_20}, the well location
for which the original POD basis was constructed, as well as the new
location for which the adaptation of the original basis will be performed are presented.
In Fig. \ref{fig: position_rates_calculation}, a comparison of the simulated fluid production
rates for the new producer location is shown for models using the original POD basis
and the adaptive POD basis.

\begin{figure}[htbp]
    \begin{center}
        \begin{subfigure}{0.49\textwidth}
            \includegraphics[width=0.9\textwidth]{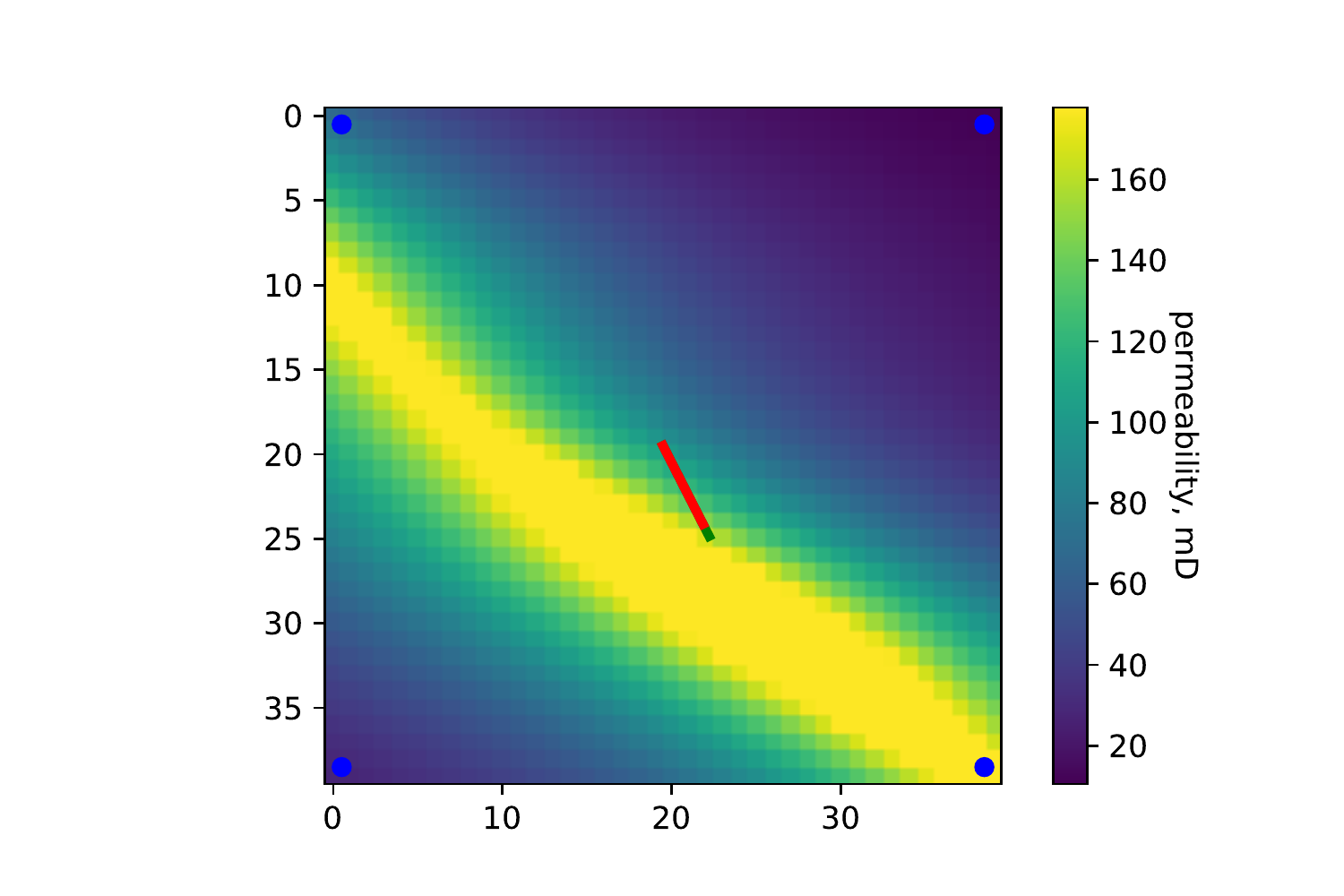}
            \caption{}
            \label{fig: well_map_130_150: map}
        \end{subfigure}
        \begin{subfigure}{0.49\textwidth}
            \includegraphics[width=0.9\textwidth]{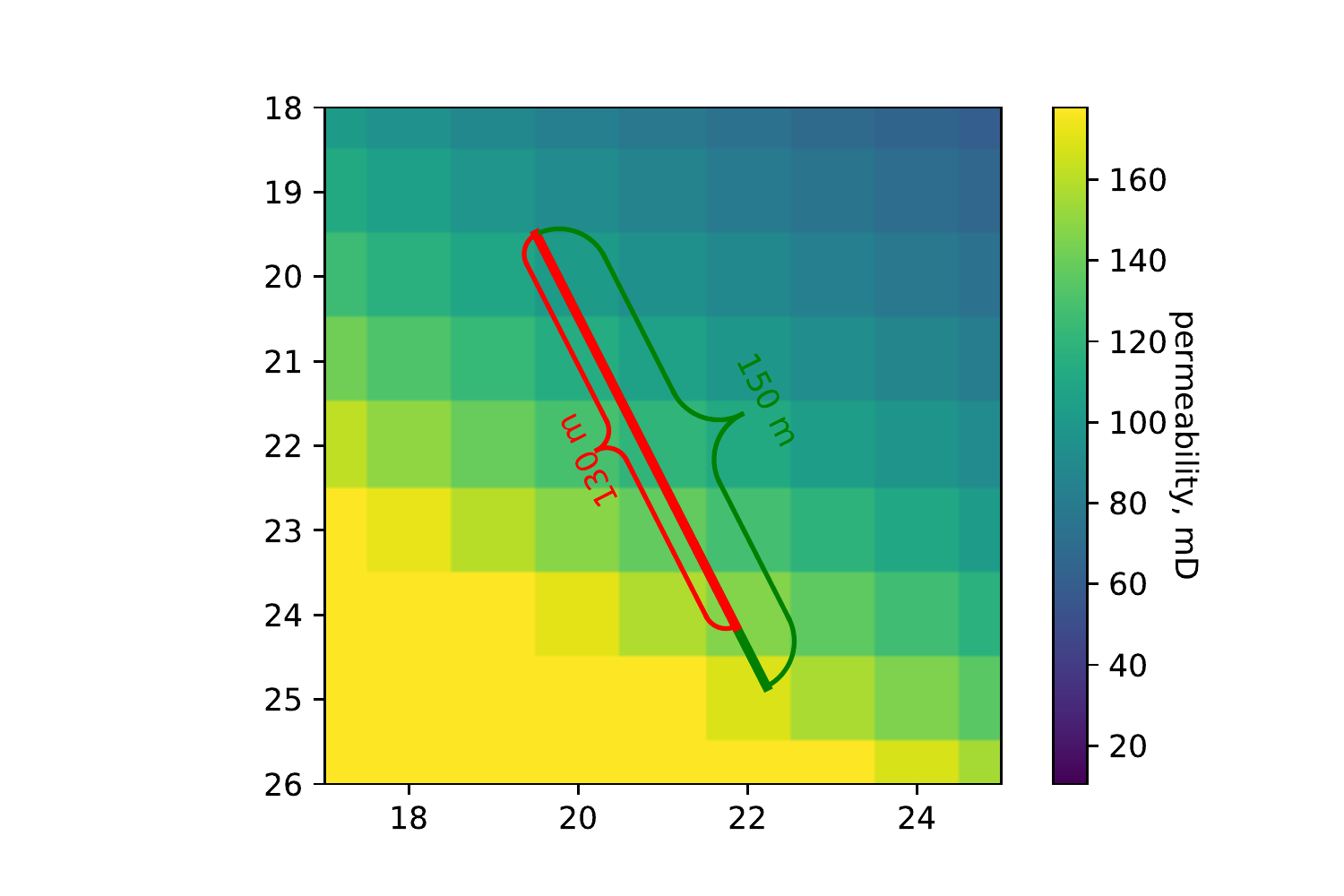}
            \caption{}
            \label{fig: well_map_130_150: zoom}
        \end{subfigure}
        \caption{\ref{sub@fig: well_map_130_150: map} - Well placement scheme. Blue circles represent the injector wells; the green line
            shows the original position and length of the producer; the red line indicates
            the new length of the producer. \ref{sub@fig: well_map_130_150: zoom} - zoomed in part of the model with the producer well.
            }
        \label{fig: well_map_130_150}
    \end{center}
\end{figure}

\begin{figure}[htbp]
    \begin{subfigure}{0.49\textwidth}
        \begin{center}
            \includegraphics[width=\textwidth]{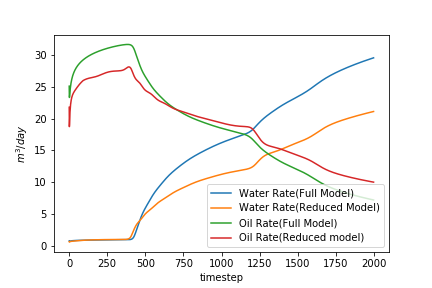}
        \end{center}
        \caption{}
        \label{fig: length_rates_calculation: original_basis}
    \end{subfigure}
    \begin{subfigure}{0.49\textwidth}
        \begin{center}
            \includegraphics[width=\textwidth]{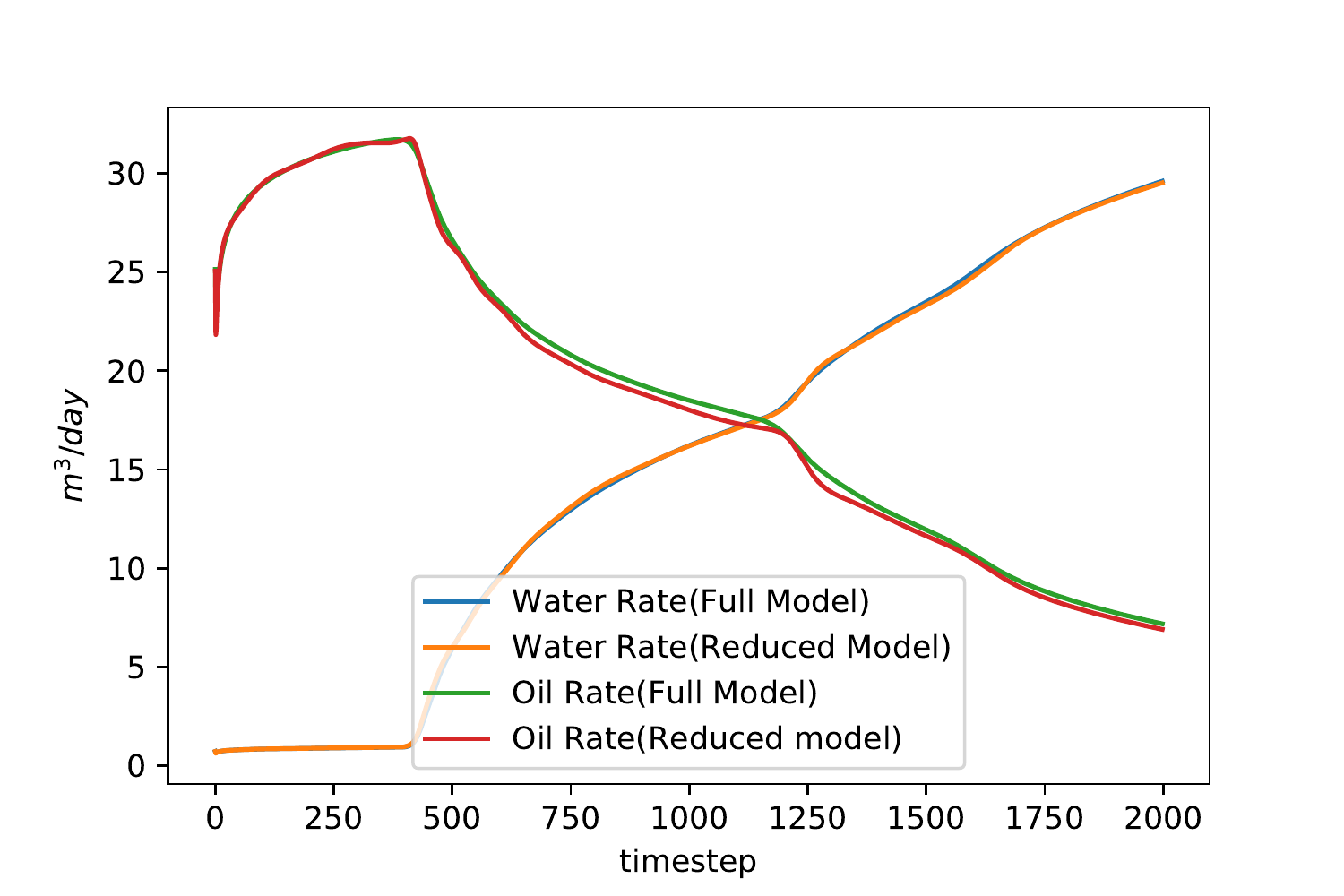}
        \end{center}
        \caption{}
        \label{fig: length_rates_calculation: adaptive_basis}
    \end{subfigure}
    \caption{Simulated fluid production rates: \ref{sub@fig: length_rates_calculation: original_basis} -
    the original POD basis, \ref{sub@fig: length_rates_calculation: adaptive_basis} - the adaptive POD basis constructed using 50 new snapshots and 1 additional component.}
    \label{fig: length_rates_calculation}
\end{figure}

In Fig. \ref{fig: well_map_130_150}, a well placement scheme with the variation of the length
of the producing section is shown.
In Fig. \ref{fig: length_rates_calculation}, a comparison of the simulated fluid production rates for this model obtained with the use of the original POD basis and with the adaptive POD basis
is presented.

One can see that the proposed approach allows us to adapt the existing ROM to a wide range of new well configurations,
including varying well orientations, length, and position at the expense of a relatively small number of
additional snapshots. In contrast, if one tries to build an universal POD-based ROM capable of accurately simulating production rates for models with variable well positions
and geometries, a set of snapshots that scales exponentially with
the number of varying parameters (well orientation, well length, etc.) is required. In statistics and in ML applications this
problem is known as the curse of dimensionality \cite{bishop_pattern_2006}. Another problem that one would
face in this case, is an increasing complexity of the snapshots space that would in turn lead to a
significant increase of the number of POD basis components
required to obtain adequate simulation accuracy and stability. This makes such an approach impractical
for problems with multiple optimization parameters.

\subsubsection{Sensitivity analysis}
\begin{figure}[htbp]
    \begin{subfigure}{0.49\textwidth}
        \begin{center}
            \includegraphics[width=\textwidth]{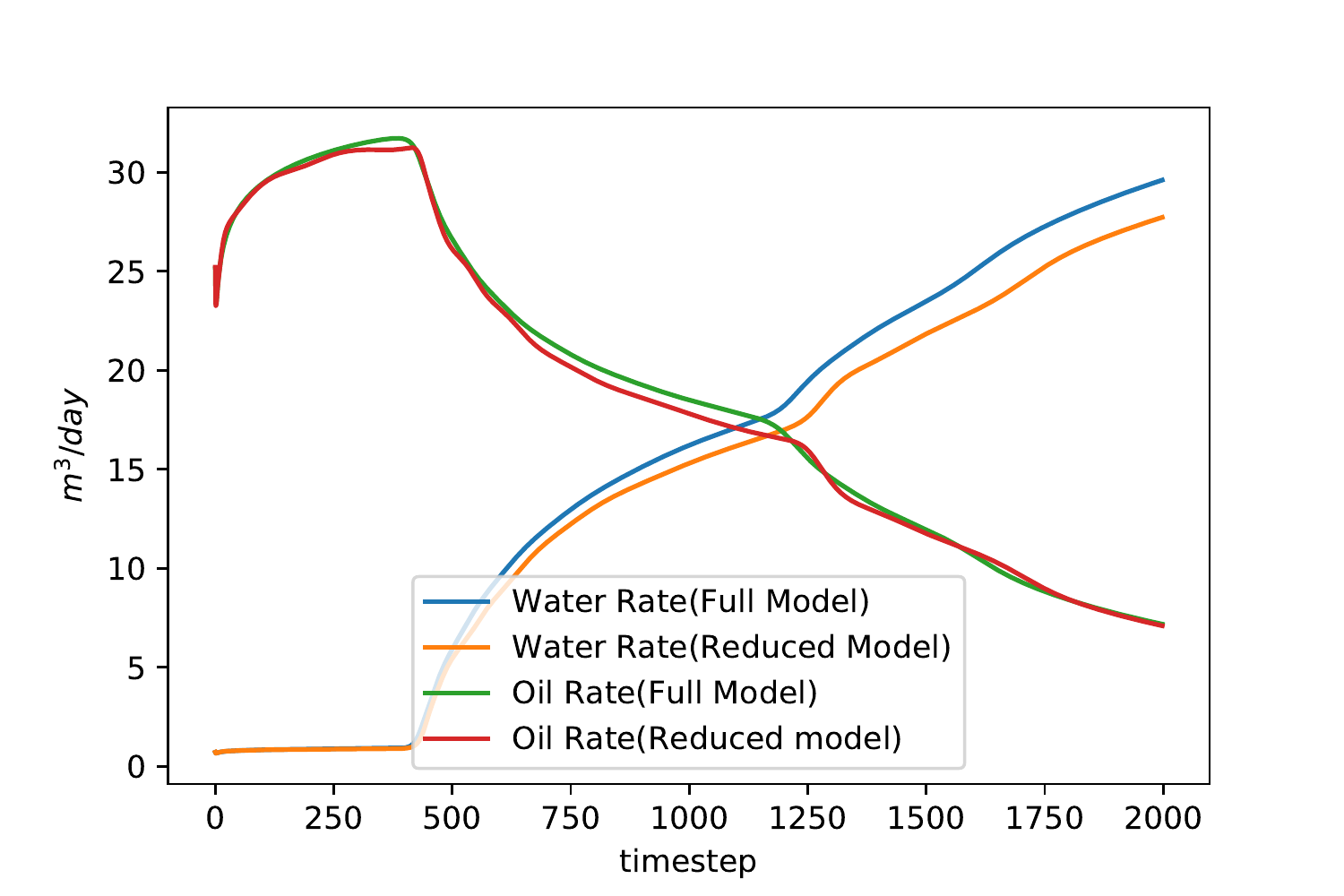}
        \end{center}
        \caption{}
        \label{fig: snapshots_number_adaptive_basis: 10_snapshots}
    \end{subfigure}
    \begin{subfigure}{0.49\textwidth}
        \begin{center}
            \includegraphics[width=\textwidth]{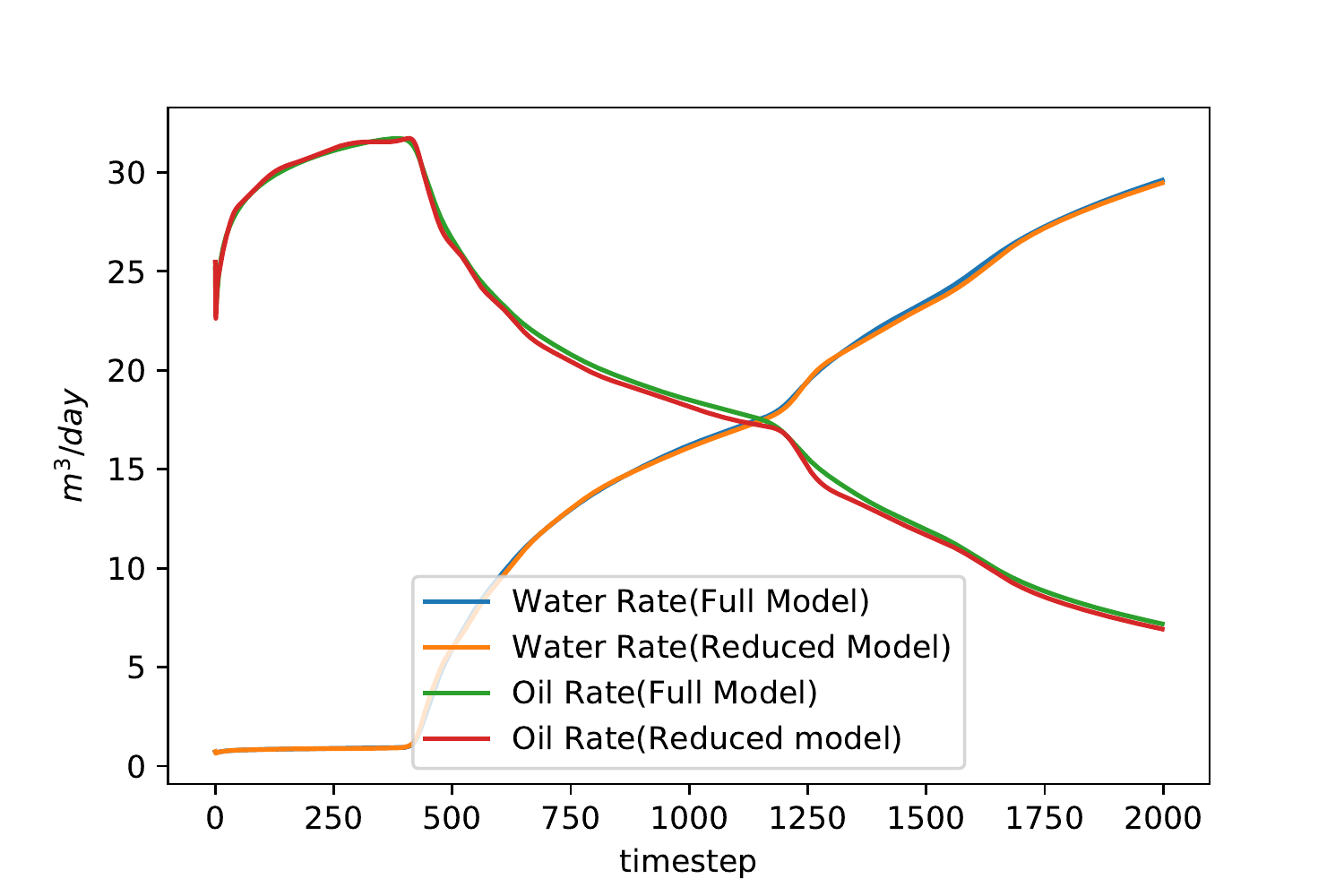}
        \end{center}
        \caption{}
        \label{fig: snapshots_number_adaptive_basis: 30_snapshots}
    \end{subfigure}
    \begin{subfigure}{0.49\textwidth}
        \begin{center}
            \includegraphics[width=\textwidth]{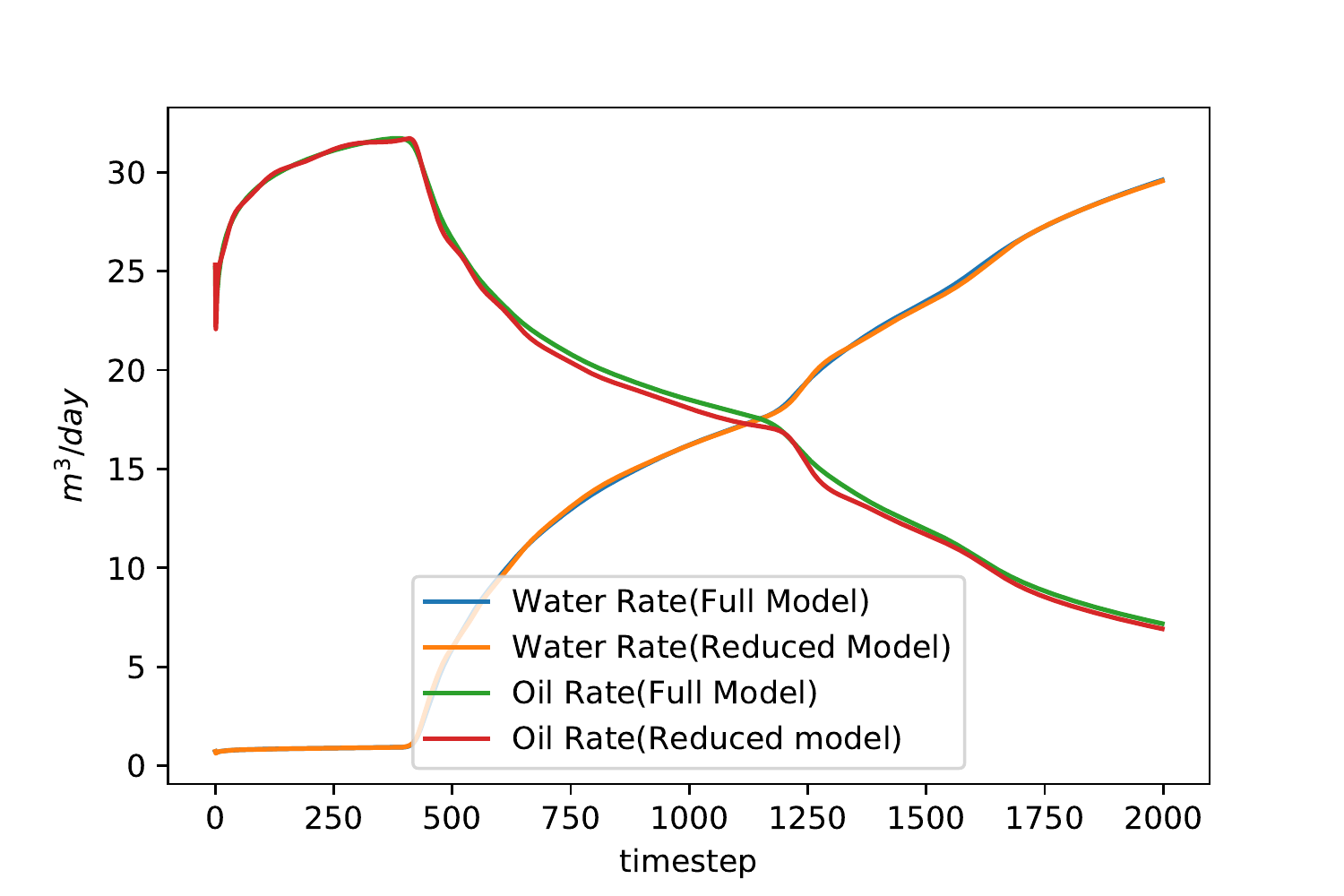}
        \end{center}
        \caption{}
        \label{fig: snapshots_number_adaptive_basis: 50_snapshots}
    \end{subfigure}
    \begin{subfigure}{0.49\textwidth}
        \begin{center}
            \includegraphics[width=\textwidth]{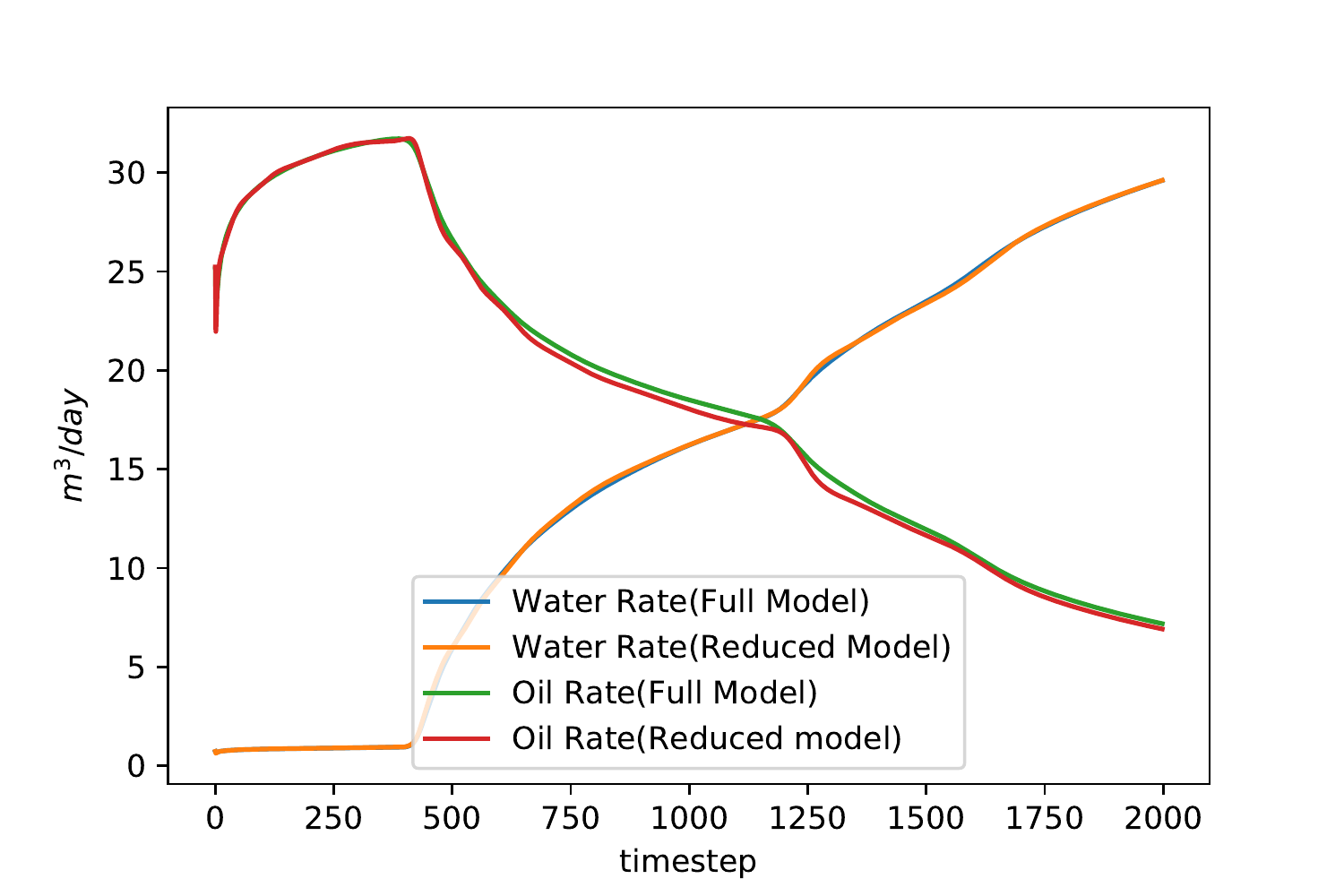}
        \end{center}
        \caption{}
        \label{fig: snapshots_number_adaptive_basis: 100_snapshots}
    \end{subfigure}
    \caption{
        Simulated fluid production rates for the model with a variable length of the producer, obtained with the adaptive POD basis using different number of snapshots:
        \ref{sub@fig: snapshots_number_adaptive_basis: 10_snapshots} - 10 snapshots,
        \ref{sub@fig: snapshots_number_adaptive_basis: 30_snapshots} - 30 snapshots,
        \ref{sub@fig: snapshots_number_adaptive_basis: 50_snapshots} - 50 snapshots,
        \ref{sub@fig: snapshots_number_adaptive_basis: 100_snapshots} - 100 snapshots.
    }
    \label{fig: snapshots_number}
\end{figure}

In order to estimate the impact of the number of additional snapshots used in the adaptive POD basis on
the accuracy of the flow rate simulations, several variants of the model with
different lengths of the producer and a different number of additional snapshots for each of the variants
(Fig. \ref{fig: well_map_130_150}) were simulated.
The results of these simulations are presented in Fig. \ref{fig: snapshots_number}.

\begin{figure}[htbp]
    \begin{center}
        \includegraphics[width=0.6\textwidth]{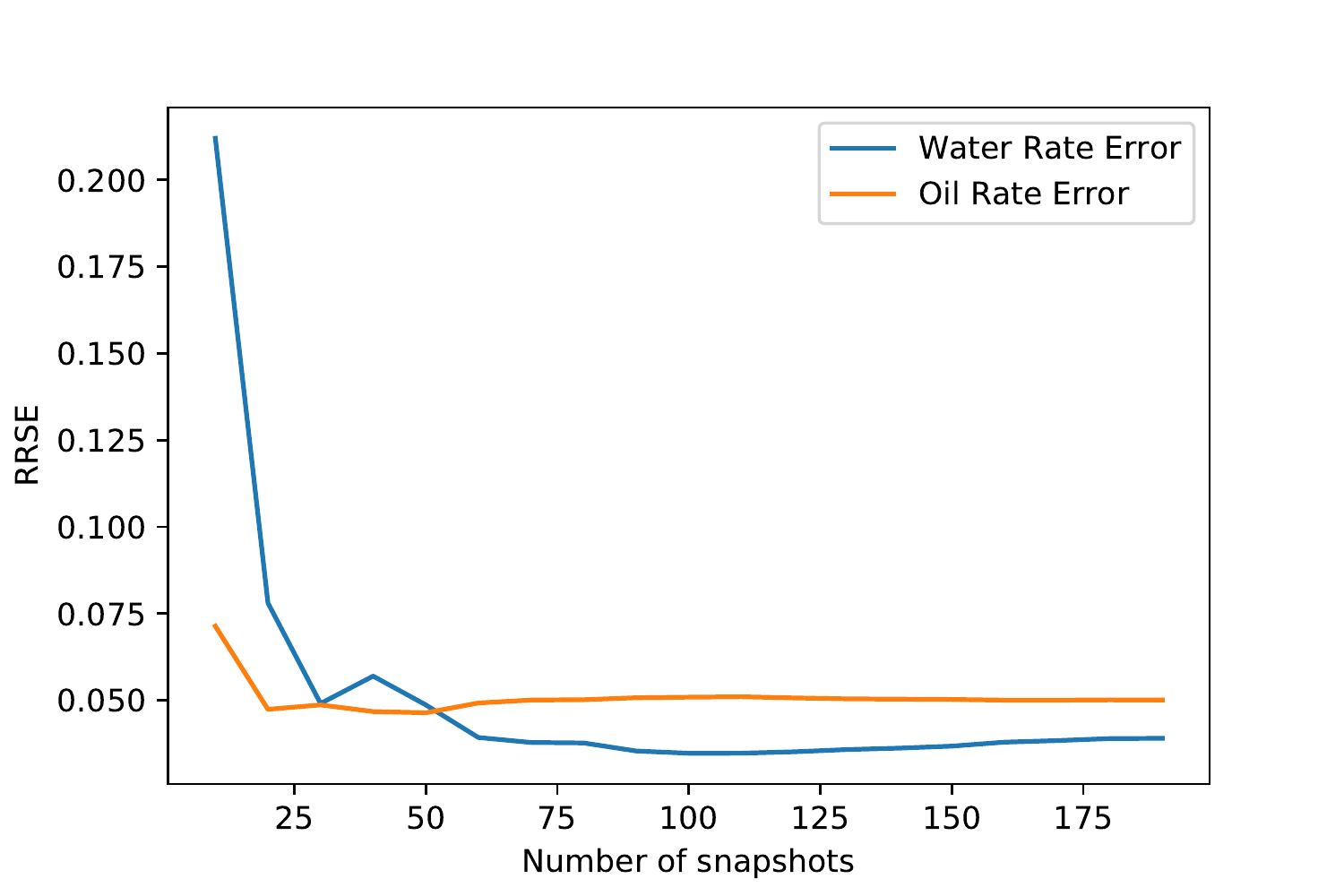}
    \end{center}
    \caption{Root relative squared error of simulated fluid production rates with respect to the full model as a function of number of additional snapshots used in the construction of the adaptive POD basis.}
    \label{fig: error_snapshots}
\end{figure}

Deviations of the simulated fluid production rates with respect to the full model were calculated for the models with the adaptive POD basis constructed using different number of additional snapshots.
The corresponding RRSE graphs are presented in Fig. \ref{fig: error_snapshots}.

Not surprisingly, the increase of the number of additional snapshots up to 30-50 generally improves the accuracy of the flow rate simulations. However, further
increase of the number of additional snapshots beyond that range does not seem to improve the accuracy.

\begin{figure}[htbp]
    \begin{subfigure}{0.49\textwidth}
        \begin{center}
            \includegraphics[width=\textwidth]{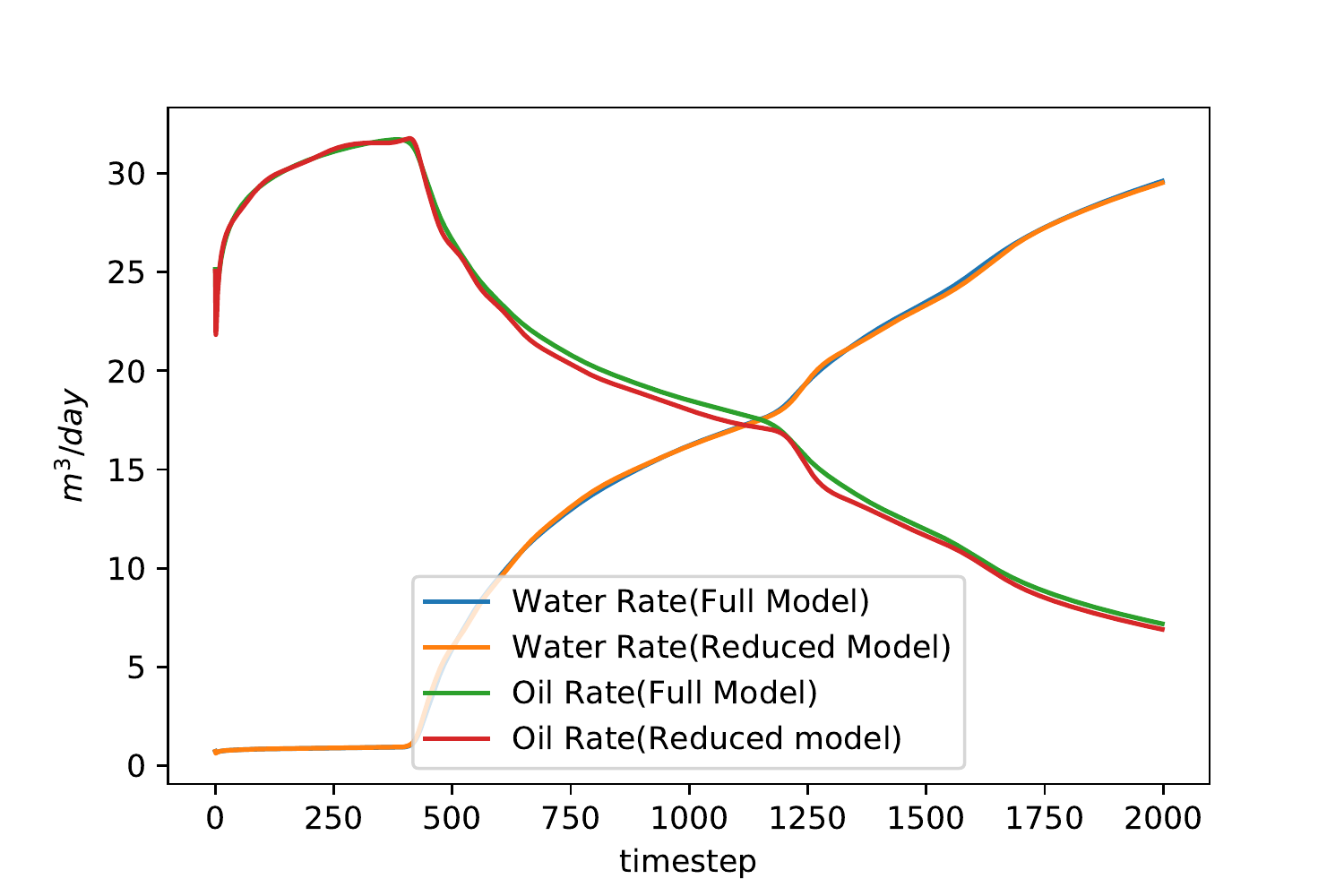}
        \end{center}
        \caption{}
        \label{fig: componet_number_adaptive_basis: 1_component}
    \end{subfigure}
    \begin{subfigure}{0.49\textwidth}
        \begin{center}
            \includegraphics[width=\textwidth]{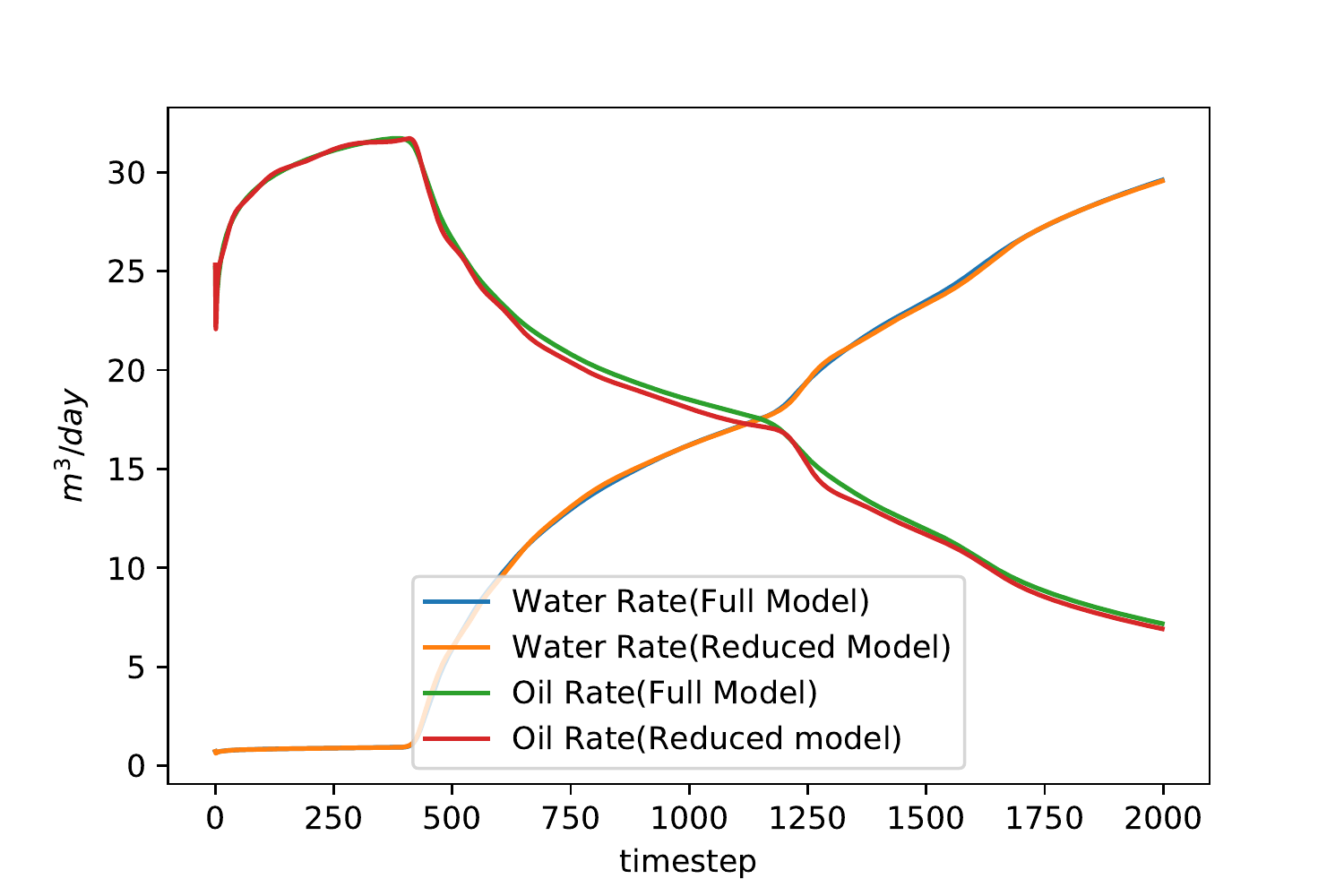}
        \end{center}
        \caption{}
        \label{fig: componet_number_adaptive_basis: 3_component}
    \end{subfigure}
    \begin{subfigure}{0.49\textwidth}
        \includegraphics[width=\textwidth]{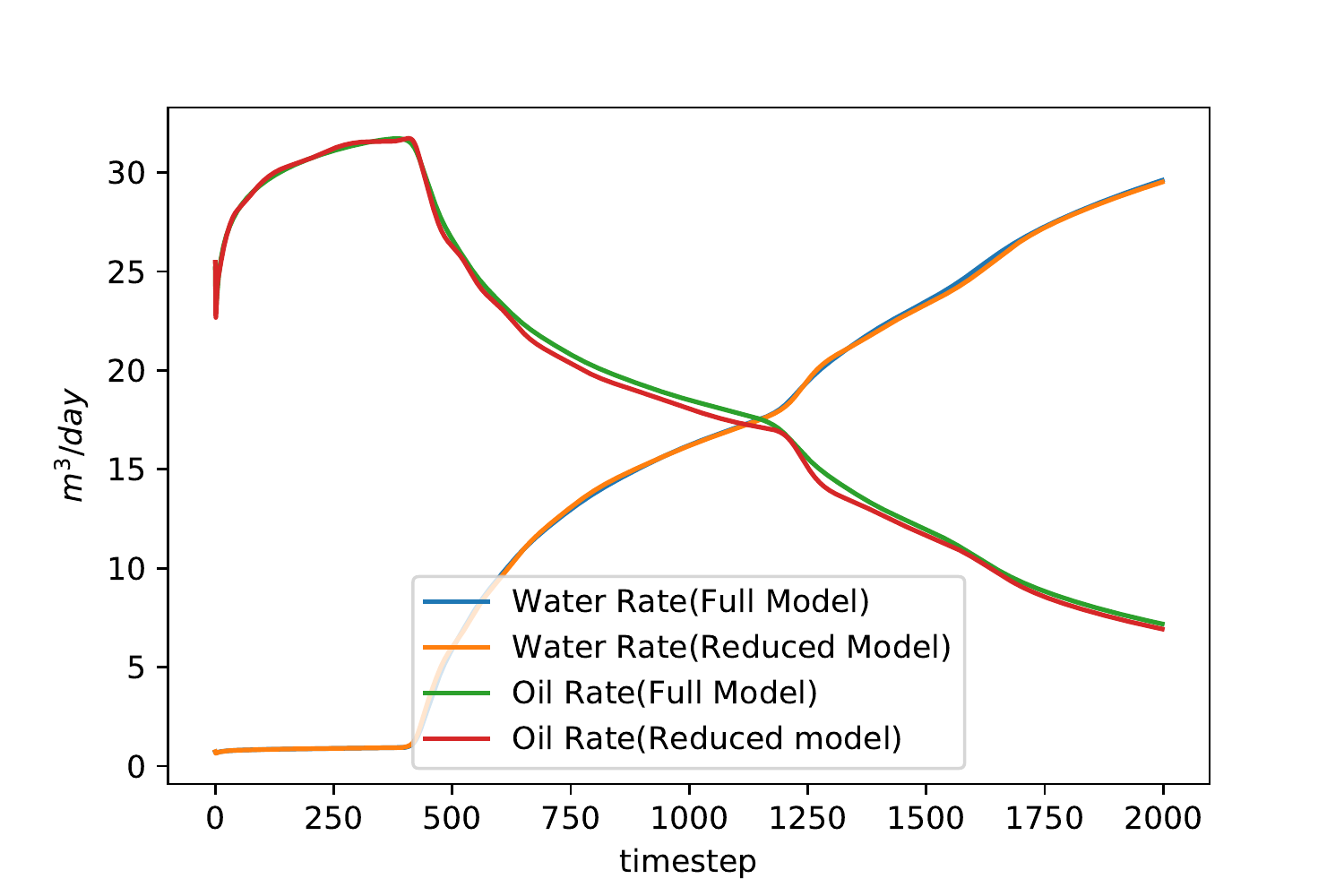}
        \caption{}
        \label{fig: componet_number_adaptive_basis: 5_component}
    \end{subfigure}
    \begin{subfigure}{0.49\textwidth}
        \begin{center}
            \includegraphics[width=\textwidth]{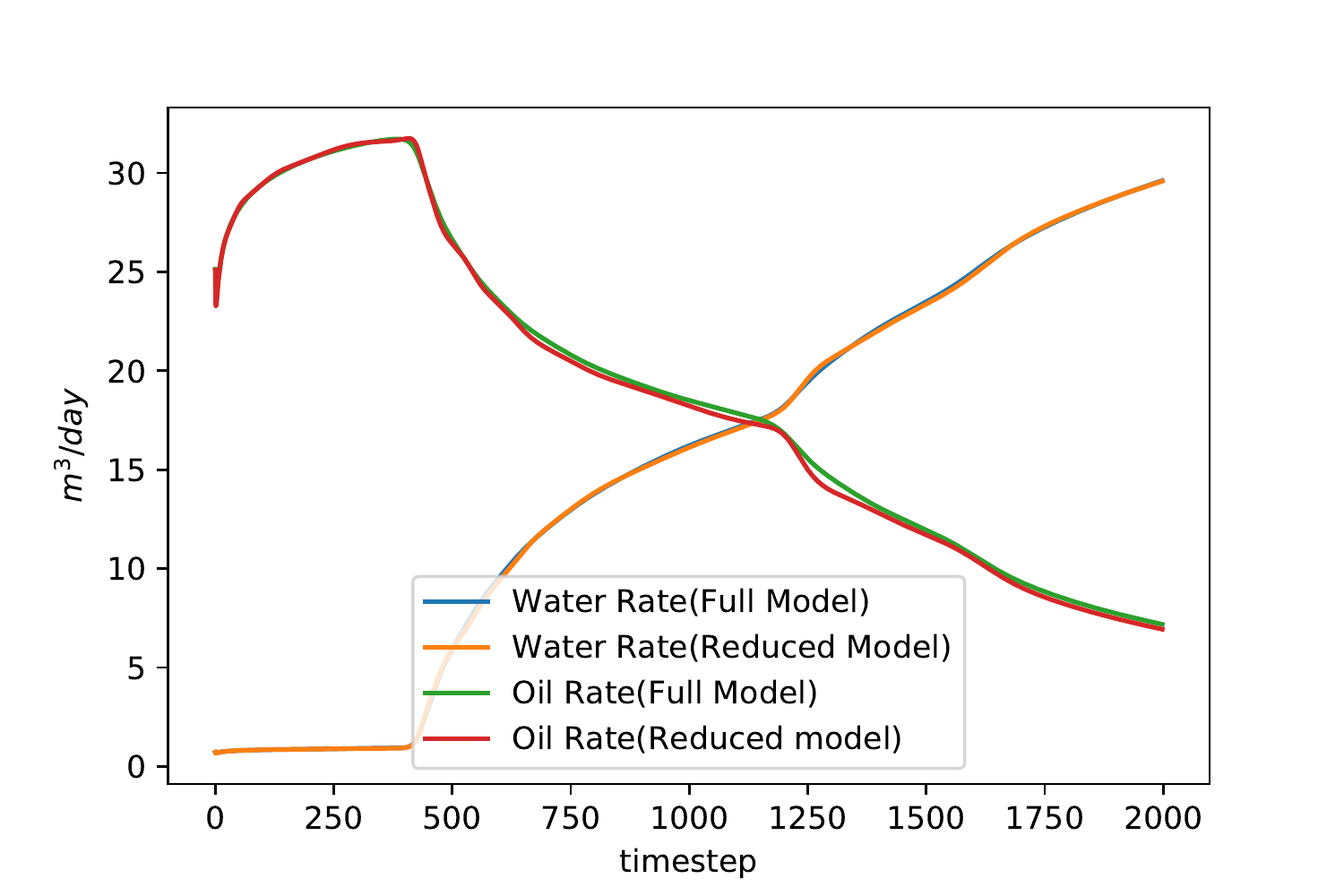}
        \end{center}
        \caption{}
        \label{fig: componet_number_adaptive_basis: 10_component}
    \end{subfigure}
    \caption{
        Simulated fluid production rates for models with the adaptive POD basis with different
        number of additional components:
        \ref{sub@fig: componet_number_adaptive_basis: 1_component} - 1 component,
        \ref{sub@fig: componet_number_adaptive_basis: 3_component} - 3 components,
        \ref{sub@fig: componet_number_adaptive_basis: 5_component} - 5 components,
        \ref{sub@fig: componet_number_adaptive_basis: 10_component} - 10 components.
    }
    \label{fig: components_number}
\end{figure}

In Fig. \ref{fig: components_number}, the simulated fluid production rates are shown for models using adaptive POD basis with different number
of additional components. In all the presented cases, the POD basis adaptation was performed
using 50 additional snapshots.

\begin{figure}[htbp]
    \begin{center}
        \includegraphics[width=0.6\textwidth]{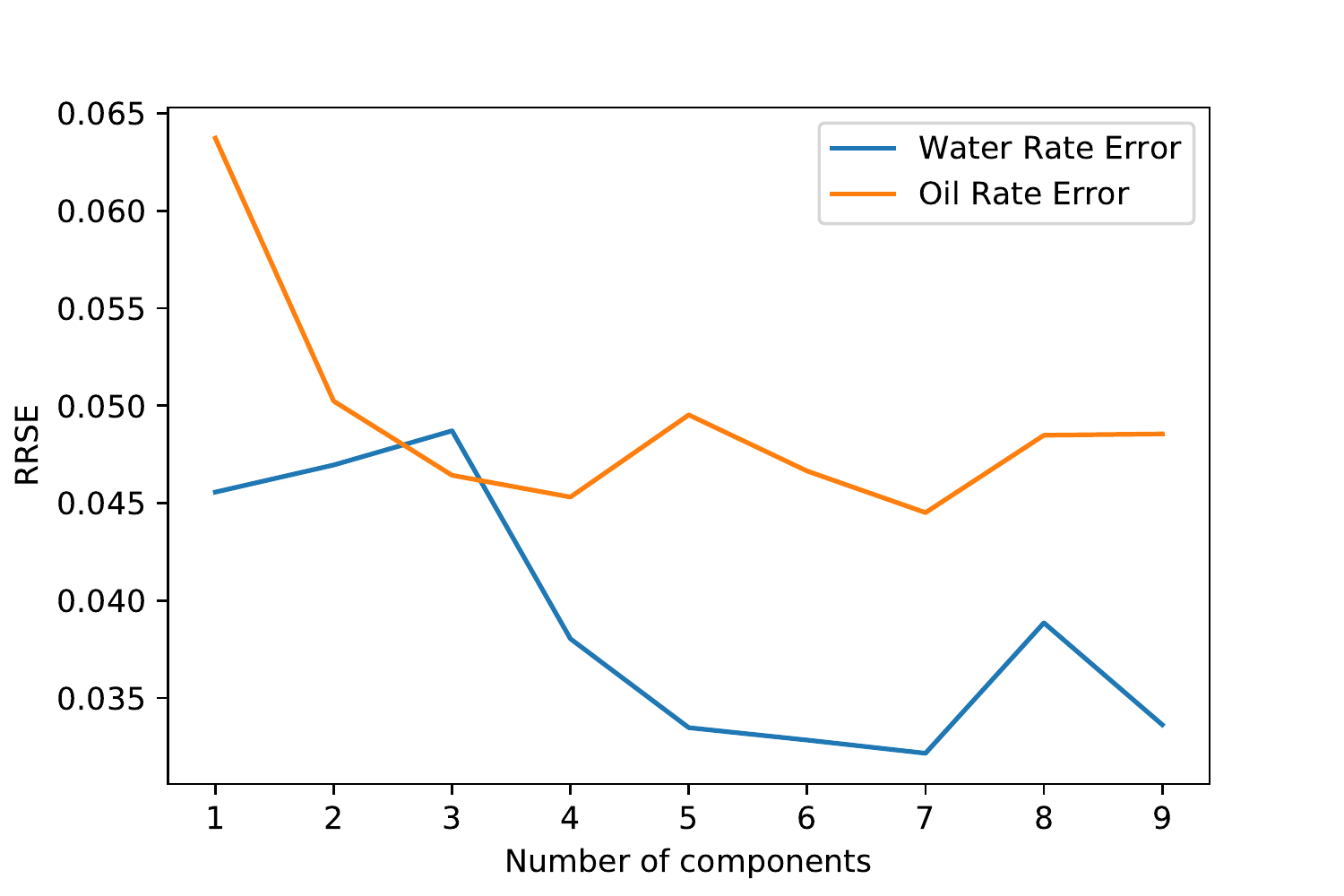}
    \end{center}
    \caption{Root relative squared error of simulated fluid production rates with respect to the full model
    as a function of the number of additional components of the adaptive POD basis.}
    \label{fig: error_components}
\end{figure}

Deviations of the simulated fluid production rates with respect to the full model were calculated for the cases with the adaptive POD basis constructed using different number of additional components.
The corresponding RRSE graphs are presented in Fig. \ref{fig: error_components}.

One can observe that the influence of the number of additional components of the adaptive basis
on the production rate simulations is relatively small, and typically 1-3
additional components are sufficient for a satisfactory basis adaptation.

\section{Summary and Further Work}\label{sec: summary}

In this work, different approaches for an efficient use of POD-Galerkin ROMs
for the reservoir simulation were explored for problems where changes of the boundary
conditions, such as well location and geometry, are essential.
In the universal basis approach, a training data set is generated such that it contains snapshots
of the solutions related to all the considered well geometries.
This approach allows us to use the same POD basis for simulating scenarios with different well
geometries. However, that universality comes at the expense of a relatively large number
of the POD basis components required in order to obtain a reasonable accuracy of the simulations.
Another drawback of this approach is the requirement of a very large training data set that should
reflect all the range of possible scenarios. This approach may be relevant
in cases where one needs to explore a relatively small range of possible well geometries.

A new approach based of the adaptation of the POD basis to varying
well configurations was introduced. This approach allows us to use the POD basis
constructed for a specific
model configuration in order to model a different problem setting. Such an adaptation of the POD
basis is achieved at a relatively low computational cost, after updating the basis with a few
additional components obtained from the snapshots of the new model configuration.
It was found that the number of such additional snapshots is substantially smaller compared to that
required for building of a new POD basis from scratch. The new adaptive method was validated
on several test cases where a POD-Galerkin ROM was used to simulate immiscible oil displacement
by water injection in the model settings where the horizontal producer had variable location, length,
and orientation.

So far, the POD reduction technique was applied only to the pressure field.
An extension of the proposed adaptive approach to reduce the saturation related part of the model
is in the scope of the future work.
Another envisaged direction of the work is related to the application of the proposed adaptive POD
approach in order to build efficient ROM simulators that are capable to address
realistic reservoir optimization problems, such as the optimal well design, completion,
and hydraulic fracturing.

\newpage
\printbibliography
\end{document}